\documentclass[11pt]{article}

\usepackage{latexsym}
\usepackage{amssymb}
\usepackage{amsthm}
\usepackage{amscd}
\usepackage{amsmath}
\usepackage[all]{xy}
\usepackage{graphicx}

\newtheorem{theorem}{Theorem}[section]

\newtheorem{lemma}{Lemma}[section]

\newtheorem{proposition}{Proposition}[section]
\newtheorem{corollary}{Corollary}[section]

\xyoption{dvips}

\numberwithin{equation}{section}

\newcommand{\FF}{\mathbb{F}}    \newcommand{\CC}{\mathbb{C}} 
   
  \def\Hom{\mathrm{Hom}} 
 \def\ZZ{\mathbb{Z}} 
 \def\X{\mathfrak{X}}  \def\Ind{\mathrm{Ind}} \def\GL{\mathrm{GL}}    \def\dd{\displaystyle}  \def\ch{\mathrm{ch}} \def\spanning{\textnormal{-span}}
\def\vphi{\varphi}  \def\P{\mathcal{P}}
  \def\wt{\mathrm{wt}}\def\Indf{\mathrm{Indf}}
  \def\sh{\mathrm{sh}}\def\hgt{\mathrm{ht}}
  \def\sss{\scriptscriptstyle} \def\ss{\scriptstyle}
\newcommand{\cT}{\mathcal{T}}\newcommand{\cD}{\mathcal{D}}  

\newcommand{\cX}{\mathcal{X}}
\newcommand{\cS}{\mathcal{S}}
\newcommand{\cB}{\mathcal{B}}

\def\cP{\mathcal{P}}
\newcommand{\bmu}{\boldsymbol{\mu}}
\newcommand{\blam}{\boldsymbol{\lambda}}
\newcommand{\bnu}{\boldsymbol{\nu}}
\newcommand{\bgamma}{\boldsymbol{\gamma}}

\newcommand{\btau}{\boldsymbol{\tau}}
\newcommand{\bdeta}{\boldsymbol{\eta}}

\newcommand{\core}{\mathrm{core}}
\newcommand{\quot}{\mathrm{quot}}

\makeatletter
\renewcommand{\@makefnmark}{\mbox{\textsuperscript{}}}
\makeatother

\CompileMatrices

\def\adots{\mathinner{\mkern2mu\raise0pt\hbox{.}  % antidiagonal dots
\mkern2mu\raise4pt\hbox{.}\mkern1mu
\raise7pt\vbox{\kern7pt\hbox{.}}\mkern1mu}}

\allowdisplaybreaks[1]

\begin{document}

\bibliographystyle{amsplain}

\title{Gelfand-Graev characters of the finite unitary groups}
\author{Nathaniel Thiem\footnote{Stanford University: \textsf{thiem@math.stanford.edu}}{ } and C. Ryan Vinroot\footnote{University of Arizona: \textsf{vinroot@math.arizona.edu}}}
\date{}

\maketitle

%%%

%%%INTRODUCTION

%%%

\begin{abstract}
Gelfand-Graev characters and their degenerate counterparts have an important role in the representation theory of finite groups of Lie type.  Using a characteristic map to translate the character theory of the finite unitary groups into the language of symmetric functions, we study degenerate Gelfand-Graev characters of the finite unitary group from a combinatorial point of view.  In particular, we give the values of Gelfand-Graev characters at arbitrary elements, recover the decomposition multiplicities of degenerate Gelfand-Graev characters in terms of tableau combinatorics, and conclude with some multiplicity consequences.
\end{abstract}

\section{Introduction\protect\footnote{MSC 2000: 20C33 (05E05)}\protect\footnote{Keywords: finite unitary group, degenerate Gelfand-Graev, multiplicity, domino tableaux, characteristic map}}

Gelfand-Graev modules have played an important role in the representation theory of finite groups of Lie type \cite{dignemichel,GG62,St67}.  In particular, if $G$ is a finite group of Lie type, then Gelfand-Graev modules of $G$ both contain cuspidal representations of $G$ as submodules, and have a multiplicity free decomposition into irreducible $G$-modules.  Thus, Gelfand-Graev modules can give constructions for cuspidal $G$-modules.  This paper uses a combinatorial correspondence between characters and symmetric functions (as described in \cite{TV05}) to examine the Gelfand-Graev character and its degenerate relatives for the finite unitary group.

Let $B^<$ be a maximal unipotent subgroup of a finite group of Lie type $G$.  Then the Gelfand-Graev character $\Gamma$ of $G$ is the character obtained by inducing a generic linear character from $B^<$ to $G$.  The {\em degenerate} Gelfand-Graev characters of $G$ are obtained by inducing arbitrary linear characters.  In the case $\GL(n,\FF_q)$, Zelevinsky \cite{zel} described the multiplicities of irreducible characters in degenerate Gelfand-Graev characters by counting multi-tableaux of specified shape and weight, and obtained the result that every irreducible appears with multiplicity one in some degenerate Gelfand-Graev character.  It is known that this multiplicity one result is not true in a general finite group of Lie type, and in fact there are characters which do not appear in any degenerate Gelfand-Graev character in general.  This follows, for example, by work of Kotlar \cite{kot} which gives a geometric description of the irreducible characters appear in some degenerate Gelfand-Graev character in general type.  In the finite unitary case, we give a combinatorial description of which irreducible characters appear in some degenerate Gelfand-Graev character, as well as a combinatorial description of a large family of characters which appear with multiplicity one.

Section 2 reviews the necessary combinatorics and results from \cite{TV05}.   Section 3 examines the Gelfand-Graev character, and uses a remarkable formula for the character values of the Gelfand-Graev character of $\GL(n,\FF_q)$ (for an elementary proof see \cite{howzwor}) to obtain the corresponding formula for ${\rm U}(n,\FF_{q^2})$.  The following is the main result of this section. 

\vspace{.15cm}

\noindent \textbf{I.} (Corollary \ref{UGGFormula})    If $\Gamma_{(n)}$ is the Gelfand-Graev character of ${\rm U}(n,\FF_{q^2})$, and $g\in {\rm U}(n,\FF_{q^2})$, then
$$\Gamma_{(n)}(g)=\biggl\{\begin{array}{@{}ll} (-1)^{\lfloor n/2\rfloor+ {\ell\choose 2}}(q^{\ell}-(-1)^{\ell})\cdots (q-1) &\text{if $g$ is unipotent, block type $(\mu_1,\mu_2,\ldots, \mu_\ell)$,}\\ 0 & \text{otherwise.}\end{array}$$
Section 4 computes the decomposition of degenerate Gelfand-Graev characters in a fashion analogous to \cite{zel}, with the main result as follows.

\vspace{.15cm}

\noindent\textbf{II.} (Theorem \ref{decomp})  The degenerate Gelfand-Graev character $\Gamma_{(k,\nu)}$ decomposes as
$$\Gamma_{(k,\nu)}=\sum_{\blam}m_{\blam} \chi^{\blam},$$
where $\blam$ is a multi-partition and $m_{\blam}$ is a nonnegative integer obtained by counting ``battery tableaux" of a given weight and shape. 

\vspace{.15cm}

\noindent In the process of proving Theorem \ref{decomp}, we obtain some combinatorial Pieri-type formulas (Lemma \ref{PieriLemma}) and  decompositions of induced characters from $\GL(n,\FF_{q^2})$ to ${\rm U}(2n,\FF_{q^2})$ (Theorem \ref{DGGBaseCase} and Theorem \ref{InducedIrreducibles}).

Section 5 concludes with a discussion of the multiplicity implications of Section 4.  In particular, we improve a multiplicity one result by Ohmori \cite{ohm}.

\vspace{.15cm}

\noindent\textbf{III.} (Theorem \ref{MultiplicityOne}) We give combinatorial conditions on multipartitions $\blam$ that guarantee that the irreducible character $\chi^{\blam}$ appears with multiplicity one in some degenerate Gelfand-Graev character.

\vspace{.15cm}

Another question one might ask is how the {\em generalized} Gelfand-Graev representations of the finite unitary group decompose.  Generalized Gelfand-Graev representations, which were defined by Kawanaka in \cite{kawan}, are obtained by inducing certain irreducible representations (not necessarily one dimensional) from a unipotent subgroup.  Rainbolt studies the generalized Gelfand-Graev representations of ${\rm U}(3, \FF_{q^2})$ in \cite{JR98}, but in the general case they seem to be significantly more complicated than the degenerate Gelfand-Graev representations.

\vspace{.25cm}

\noindent\textbf{Acknowledgements.}  We would like to thank G. Malle for suggesting the questions that led to the results in Section 5,  S. Assaf for a helpful discussion regarding Section 5.1, and T. Lam helping us connect Lemma \ref{PieriLemma} to the literature.

%%%

%%%PRELIMINARIES

%%%

\section{Preliminaries}

\subsection{Partitions}

Let 
$$\cP=\bigcup_{n\geq 0} \cP_n,\qquad \text{where} \qquad \cP_n=\{\text{partitions of $n$}\}.$$
For $\nu = (\nu_1, \nu_2, \ldots, \nu_l) \in \cP_n$, the {\em length} $\ell(\nu)$ of $\nu$ is the number of parts $l$, and the {\em size} $|\nu|$ of $\nu$ is the sum of the parts $n$.  Let $\nu'$ denote the conjugate of the partition $\nu$.    We also write
$$\nu=(1^{m_1(\nu)} 2^{m_2(\nu)} \cdots),\qquad\text{where}\qquad m_i(\nu) = |\{ j\in \ZZ_{\geq 1} \mid \nu_j = i \}|.$$
Define $n(\nu)$ to be
$$n(\nu)=\sum_i (i - 1) \nu_i.$$
If $\mu, \nu \in \cP$, we define $\mu \cup \nu \in \cP$ to be the partition of size $|\mu| + |\nu|$ whose set of parts is the union of the parts of $\mu$ and $\nu$.  For $k \in \ZZ_{\geq 1}$, let $k\nu=(k\nu_1, k\nu_2, \ldots)$, and if every part of $\nu$ is divisible by $k$, then we let $\nu/k=(\nu_1/k, \nu_2/k, \ldots)$.  A partition $\nu$ is {\em even} if $\nu_i$ is even for $1 \leq i \leq \ell(\nu)$.

\subsection{The ring of symmetric functions}

Let $X=\{X_1,X_2,\ldots\}$ be an infinite set of variables and let 
$$\Lambda(X)=\CC[p_1(X),p_2(X),\ldots],\qquad \text{where}\qquad p_k(X)=X_1^k+X_2^k+\cdots,$$
be the graded $\CC$-algebra of symmetric functions in the variables $\{X_1,X_2,\ldots\}$.  For $\nu=(\nu_1,\nu_2,\ldots,\nu_\ell)\in \cP$, the \emph{power-sum symmetric function} $p_\nu(X)$ is
$$p_\nu(X)=p_{\nu_1}(X)p_{\nu_2}(X)\cdots p_{\nu_\ell}(X).$$

The irreducible characters $\omega^\lambda$ of $S_n$ are indexed by $\lambda\in\cP_n$.  Let $\omega^\lambda(\nu)$ be the value of $\omega^\lambda$ on a permutation with cycle type $\nu$.  

The \emph{Schur function} $s_\lambda(X)$ is given by
\begin{equation} \label{SchurToPower} s_\lambda(X)=\sum_{\nu\in \cP_{|\lambda|}}  \omega^\lambda(\nu)z_{\nu}^{-1} p_\nu(X),\quad\text{where}\quad  z_{\nu} = \prod_{i \geq 1} i^{m_i} m_i! \end{equation}
is the order of the centralizer in $S_n$ of the conjugacy class corresponding to $\nu=(1^{m_1}2^{m_2}\cdots)\in \cP$.

Fix $t\in \CC^\times$.  For $\mu\in \cP$, the \emph{Hall-Littlewood symmetric function} $P_\mu(X;t)$ is given by
\begin{equation}\label{SchurToHall}
s_\lambda(X)=\sum_{\mu\in \cP_{|\lambda|}} K_{\lambda\mu}(t) P_\mu(X;t),
\end{equation}
where $K_{\lambda\mu}(t)$ is the Kostka-Foulkes polynomial (as in \cite[III.6]{Mac}). For $\nu,\mu\in \cP_n$, the \emph{classical Green function} $Q_\nu^\mu(t)$ is given by
\begin{equation}\label{PowerToHall} p_\nu(X)=\sum_{\mu\in \cP_{|\nu|}} Q_\nu^\mu(t^{-1})t^{n(\mu)} P_\mu(X;t).\end{equation}

As a graded ring,
\begin{align*} \Lambda(X) &=\CC\spanning\{p_\nu(X)\ \mid\ \nu\in \cP\}\\
					 &=\CC\spanning\{s_\lambda(X)\ \mid\ \lambda\in \cP\}\\
					 &=\CC\spanning\{P_\mu(X;t)\ \mid\ \mu\in \cP\},
\end{align*}
with change of bases given in (\ref{SchurToPower}), (\ref{SchurToHall}), and (\ref{PowerToHall}).

We will also use several products formulas in the ring of symmetric functions.   The usual product on Schur functions
\begin{equation} \label{LRCoefficients} s_\nu s_\mu=\sum_{\lambda\in \cP} c_{\nu\mu}^\lambda s_\lambda\end{equation}
gives us the Littlewood-Richardson coefficients $c_{\nu\mu}^\lambda$.  The \emph{plethysm} of $p_\nu$ with $p_k$ is
$$p_\nu \circ p_k= p_{k\nu}.$$
Thus, we can consider the nonnegative integers $c_\lambda^\gamma$ given by
\begin{equation} s_\lambda\circ p_k = \sum_{\nu\in \P_{|\lambda|}} \frac{\omega^\lambda(\nu)}{z_\nu} p_{k\nu} 
=\sum_{\gamma\in\cP_{k|\lambda|}} c_\lambda^\gamma s_\gamma.\label{PlethysmCoefficients}
\end{equation}
Chen, Garsia, and Remmel \cite{cgr} give a combinatorial algorithm for computing the coefficients $c_\lambda^\gamma$.  We will use the case $k=2$ in Section \ref{SectionBaseCase}.  

\vspace{.5cm}

\noindent\textbf{Remark.}  The unipotent characters $\chi^{\tilde\lambda}$ of $\GL(n,\FF_{q^2})$ are indexed by partitions $\tilde\lambda$ of $n$ and the unipotent characters $\chi^\gamma$ of $\mathrm{U}(2n,\FF_{q^2})$ are indexed by partitions $\gamma$ of $2n$.  It will follow from Theorem \ref{InducedIrreducibles} that
$$R_{\GL(n,\FF_{q^2})}^{\mathrm{U}(2n,\FF_{q^2})}(\chi^{\tilde\lambda})=\sum_{|\gamma|=2|\tilde\lambda|} c_{\tilde\lambda}^\gamma \chi^\gamma,$$
where $R_H^G$ is Harish-Chandra induction.

\subsection{The finite unitary groups}

Let $\bar{G}_n=\GL(n,\bar\FF_q)$ be the general linear group with entries in the algebraic closure of the finite field $\FF_q$ with $q$ elements.    

For the Frobenius automorphisms $\tilde{F},F,F':\bar{G}_n\rightarrow \bar{G}_n$ given by
\begin{align}
\tilde{F}((a_{ij})) & = (a_{ij}^q),& %\text{for $(a_{ij})\in \bar{G}_n$,}
\notag\\
F((a_{ij})) & = (a_{ji}^q)^{-1},& %\text{for $(a_{ij})\in \bar{G}_n$,}
\label{FrobeniusDefinitions}\\
F'((a_{ij})) & =  (a_{n-j,n-i}^q)^{-1},& \text{where $(a_{ij})\in \bar{G}_n$,}\notag
\end{align}
let
\begin{equation}\label{GroupDefinitions}
\begin{split}
G_n & =\bar{G}_n^{\tilde{F}}=\{a\in \bar{G}_n\ \mid\ \tilde{F}(a)=a\},\\
U_n & =\bar{G}_n^F=\{a\in \bar{G}_n\ \mid\ F(a)=a\},\\
U_n' & =\bar{G}_n^{F'}=\{a\in \bar{G}_n\ \mid\ F'(a)=a\}.
\end{split}
\end{equation}
Then $G_n= \GL(n,\FF_{q})$ and $U_n'\cong U_n$ are isomorphic to the finite unitary group $\mathrm{U}(n,\FF_{q^2})$.  In fact, if follows from the Lang-Steinberg theorem that $U_n'$ and $U_n$ are conjugate subgroups of $\bar{G}_n$.  

For $k\in \ZZ_{\geq 0}$, let
$$\tilde{T}_{(k)}=\bar{G}_1^{\tilde{F}^k}\cong\FF_{q^k}^\times\qquad \text{and}\qquad T_{(k)}=\bar{G}_1^{F^k}\cong\left\{\begin{array}{ll} \FF_{q^k}^\times & \text{if $k$ is even}, \\ \{t\in \FF_{q}\ \mid\ t^{q^k+1}=1\} &\text{if $k$ is odd.}\end{array}\right.$$
For every partition $\eta=(\eta_1,\eta_2,\ldots,\eta_\ell)\in \cP_n$ let
\begin{align*}
T_\eta &=T_{(\eta_1)}\times T_{(\eta_2)}\times\cdots \times T_{(\eta_\ell)}\\
\tilde{T}_\eta &= \tilde{T}_{(\eta_1)}\times \tilde{T}_{(\eta_2)}\times\cdots \times \tilde{T}_{(\eta_\ell)}.
\end{align*}
Every maximal torus of $G_n$ is isomorphic to $\tilde{T}_\eta$ for some $\eta\in \cP_n$, and every maximal torus of $U_n$ is isomorphic to $T_\eta$ for some $\eta\in \cP_n$.

\subsection{Multipartitions}\label{multipartitions}

Let $F:\bar{G}_n\rightarrow \bar{G}_n$ be as in (\ref{FrobeniusDefinitions}), and let $\bar\FF_q^*=\{\vphi:\bar\FF_q^\times\rightarrow \CC^\times\}$ be the group of multiplicative complex-valued characters of $\bar\FF_q^\times$. Consider
\begin{align*}
\Phi=\{F\text{-orbits of $\bar{\FF}_q^\times$}\} \qquad \text{and}\qquad  \Theta  = \{F\text{-orbits of $\bar{\FF}_q^*$}\}.
\end{align*}

For $\cX\in \{\Phi,\Theta\}$, an \emph{$\cX$-partition} $\blam=(\blam^{(x_1)}, \blam^{(x_2)},\ldots)$ is a sequence of partitions indexed by $\cX$.   The \emph{size} of $\blam$ is
$$|\blam|=\sum_{x\in \cX} |x||\blam^{(x)}|,$$
where $|x|$ is the size of the orbit $x$.

Let
$$\cP^\cX=\bigcup_{n\geq 0} \cP_n^\cX,\qquad\text{where}\qquad \cP_n^\cX=\{\text{$\cX$-partitions of size $n$}\}.$$
For $\blam\in \cP^\cX$, let
$$\ell(\blam) = \sum_{x\in \cX} \ell(\blam^{(x)})\qquad \text{and}\qquad n(\blam)=\sum_{x\in \cX} |x| n(\blam^{(x)}).$$
The {\em conjugate} of $\blam \in \cP^\cX$ is the $\cX$-partition $\blam'$ defined by $\blam'^{(x)} = (\blam^{(x)})'$, and if $\bmu, \blam \in \cP^\cX$, then $\bmu \cup \blam \in \cP^\cX$ is defined by $(\bmu \cup \blam)^{(x)} = \bmu^{(x)} \cup \blam^{(x)}$.

The \emph{semisimple part} $\blam_s$ of an $\cX$-partition $\blam$ is the $\cX$-partition given by
\begin{equation} \label{muSemisimple}\blam_s^{(x)}=(1^{|\blam(x)|}), \qquad\text{for $x\in \cX$}.\end{equation}
For $\blam \in \cP^{\cX}$, define the set $\cP^{\blam}_s$ by
$$ \cP^{\blam}_s = \{ \bmu \in \cP^{\cX} \mid \bmu_s = \blam_s \}.$$
The \emph{unipotent part} $\blam_u$ of $\blam$ is the $\cX$-partition given by 
\begin{equation}\label{muUnipotent}\blam_u^{(\{1\})} \quad \text{has parts} \quad \{|x|\blam^{(x)}_i \ \mid \ x \in \mathcal{X}, i = 1, \ldots, \ell(\blam(x)) \},
\end{equation}
where $\{1\}$ is the orbit containing $1$ in $\Phi$ or the trivial character in $\Theta$, and $\blam_u^{(x)} = \emptyset$ when $x \neq \{1\}$.

Note that we can think of ``normal" partitions as $\cX$-partitions $\blam$ that satisfy $\blam_u=\blam$.  By a slight abuse of notation, we will sometimes interchange the mulipartition $\blam_u$ and the partition $\blam_u^{(\{1\})}$.  For example, $T_{\blam_u}$ will denote the torus corresponding to the partition $\blam_u^{(\{1\})}$.

Given the torus $T_\eta$, $\eta=(\eta_1,\eta_2,\ldots, \eta_\ell)\in \cP_n$, there is a natural surjection
\begin{equation}\label{TorusToTheta}
\begin{array}{rccc}  \btau_\Theta: &  \{ \theta=\theta_1\otimes\theta_2\otimes\cdots\otimes\theta_\ell\in \Hom(T_\eta,\CC^\times)\}& \longrightarrow & \{\bnu\in \cP^\Theta\ \mid\ \bnu_u^{(\{1\})}=\eta\}\\ &  \theta=\theta_1\otimes\theta_2\otimes \cdots \otimes \theta_\ell & \mapsto  & \btau_\Theta(\theta),\end{array}
\end{equation}
where
$$\btau_\Theta(\theta)^{(\vphi)}=(\eta_{i_1}/|\vphi|,\eta_{i_2}/|\vphi|,\ldots, \eta_{i_r}/|\vphi|),\qquad \text{where}\qquad \theta_{i_1},\theta_{i_2},\ldots, \theta_{i_r}\in \vphi.$$
It follows from a short calculation that if $\bnu\in \cP^\Theta$ has support $\{\vphi_1,\vphi_2,\ldots, \vphi_r\}$, then the preimage $\btau_{\Theta}^{-1}(\nu)$ has size
\begin{align} 
&\prod_{j=1}^{r} |\vphi_j|^{\ell(\bnu^{(\vphi_j)})} \prod_{i\geq 1}\binom{m_i(\bnu_u^{(\{1\})})}{m_{i/|\vphi_1|}(\bnu^{(\vphi_1)}),\ m_{i/|\vphi_2|}(\bnu^{(\vphi_2)}),\ \cdots,\ m_{i/|\vphi_r|}(\bnu^{(\vphi_r)})}\notag \\
& = \prod_{\vphi \in \Theta} |\vphi|^{\ell(\bnu^{(\vphi)})} \prod_{i \ge 1} \frac{\big(m_{i}(\bnu_u^{(\{1\})}) \big)!}{\prod_{\vphi \in \Theta}(m_{i/|\vphi|}(\bnu^{(\vphi)}))!} \label{TorusMultiplicities}
\end{align} 

The conjugacy classes $K^{\bmu}$ of $U_n$ are parameterized by $\bmu \in \cP_n^{\Phi}$, a fact on which we elaborate in Section \ref{UnitaryCharacteristic}.   We have another natural surjection,
\begin{equation}\label{TorusToPhi}
\begin{array}{rccc} \btau_\Phi: &T_\eta & \rightarrow & \{\bnu\in\cP^\Phi\ \mid\ \bnu_u^{(\{1\})}=\eta\}\\ & t=(t_1,t_2,\ldots,t_\ell) & \mapsto & \btau_\Phi(t_1)\cup \btau_\Phi(t_2)\cup\cdots\cup \btau_\Phi(t_\ell),\end{array}
\end{equation}
where 
$$\btau_\Phi(t_i)=\bmu', \qquad\text{if $t_i\in K^{\bmu}$ in $U_{\eta_i}$.}$$

\subsection{The characteristic map}\label{UnitaryCharacteristic}

For every $f\in \Phi$, let $X^{(f)}=\{X^{(f)}_1,X^{(f)}_2,\ldots\}$ be an infinite set of variables, and  
for every $\vphi\in \Theta$, let $Y^{(\vphi)}=\{Y^{(\vphi)}_1,Y^{(\vphi)}_2,\ldots\}$ be an infinite set of characters.  We relate symmetric functions in the variables $X^{(f)}$ to those in the variables $Y^{(\vphi)}$ through the transform
$$p_k(Y^{(\vphi)})=(-1)^{k|\phi|-1}\sum_{x\in T_{k|\vphi|}} \xi(x) p_{k|\vphi|/|f_x|}(X^{(f_x)}),
\qquad\text{where $\xi\in \vphi$, $x\in f_x$.}$$  
The \emph{ring of symmetric functions} $\Lambda$ is
\begin{equation*}
\Lambda=\bigotimes_{f\in \Phi} \Lambda(X^{(f)})=\bigotimes_{\vphi\in \Theta} \Lambda(Y^{(\vphi)}).
\end{equation*}
For $\bmu\in \cP^\Phi$, the Hall-Littlewood polynomial $P_{\bmu}$ is
$$P_{\bmu}=(-q)^{-n(\bmu)}\prod_{f\in \Phi} P_{\bmu^{(f)}}(X^{(f)};(-q)^{-|f|}),$$
and for $\blam\in \cP^\Theta$, the power-sum symmetric function $p_{\blam}$ and the Schur function $s_{\blam}$ are
$$p_{\blam}= \prod_{\vphi\in \Theta} p_{\blam^{(\vphi)}}(Y^{(\vphi)}) \qquad\text{and}\qquad 
s_{\blam}= \prod_{\vphi\in \Theta} s_{\blam^{(\vphi)}}(Y^{(\vphi)}).$$
For $\bmu, \bnu \in \cP^\Phi$, the Green function is
$$Q_{\bnu}^{\bmu}(-q) = \prod_{f \in \Phi_{\bmu}} Q_{\bnu^{(f)}}^{\bmu^{(f)}}\big((-q)^{|f|}\big),$$
where $\Phi_{\bmu} = \{ f \in \Phi \; \mid \; \bmu^{(f)} \neq \emptyset \}$.
As a graded rings,
\begin{align*}
 \Lambda  &=\CC\spanning\{p_{\bnu}\ \mid\ \bnu\in \cP^\Theta\}\\
			 	      &=\CC\spanning\{s_{\blam}\ \mid\ \blam\in \cP^\Theta\}\\
				      &=\CC\spanning\{P_{\bmu}\ \mid\ \bmu\in \cP^\Phi\}.
\end{align*}

The conjugacy classes $K^{\bmu}$ of $U_n$ are indexed by $\bmu\in \cP^\Phi_n$ and the irreducible characters $\chi^{\blam}$ of $U_n$ are indexed by $\blam\in \cP^\Theta_n$ \cite{ennolaconj, ennola}.  Thus, the ring of class functions $C_n$ of $U_n$ is given by
\begin{align*}
C_n&=\CC\spanning\{\chi^{\blam}\ \mid\ \blam\in \cP^\Theta\}\\
 &=\CC\spanning\{\kappa^{\bmu}\ \mid\ \bmu\in \cP^\Theta\},
\end{align*}
where $\kappa^{\bmu}:U_n\rightarrow \CC$ is given by
$$\kappa^{\bmu}(g)=\left\{\begin{array}{ll} 1 & \text{if $g\in K^{\bmu}$}\\ 0 & \text{otherwise.}\end{array}\right.$$
We let $\chi^{\blam}(\bmu)$ denote the value of the character $\chi^{\blam}$ on any element in the conjugacy $K^{\bmu}$. 

For $\bnu\in \cP^\Theta_n$, let the Deligne-Lusztig character $R_{\bnu} = R_{\bnu}^{U_n}$ be given by
$$R_{\bnu}=R_{T_{\bnu_u}}^{U_n}(\theta)$$
where  $\theta\in \Hom(T_{\bnu_u},\CC^\times)$ is any homomorphism such that $\btau_{\Theta}(\theta)=\bnu$ (see (\ref{TorusToTheta})).

Let $C=\bigoplus_{n\geq 1} C_n$ so that
\begin{align*}
C &=\CC\spanning\{\chi^{\blam}\ \mid\ \blam\in \cP^\Theta\}\\
&=\CC\spanning\{\kappa^{\bmu}\ \mid\ \bmu\in \cP^\Phi\}\\
&=\CC\spanning\{R_{\bnu}\ \mid\ \bnu\in \cP^\Theta\}\\
\end{align*}
is a ring with multiplication given by
$$R_{\blam} R_{\bdeta}= R_{\blam \cup \bdeta}.$$

The next theorem follows from the results of \cite{dignemichel,ennola,hotspr,  kawan, luszsrin, TV05}.  A summary of the relevant results in these papers and how they imply the following theorem is given in \cite{TV05}.
\begin{theorem}[Characteristic Map] \label{CharacteristicMap}
The map
$$\begin{array}{rccc} \ch: & C & \rightarrow & \Lambda\\ & \chi^{\blam} & \mapsto & (-1)^{\lfloor |\blam|/2\rfloor+n(\blam)} s_{\blam}\\  &\kappa^{\bmu} & \mapsto & P_{\bmu}\\
& R_{\bnu} & \mapsto & (-1)^{|\bnu|-\ell(\bnu)} p_{\bnu}\end{array}$$
is an isometric ring isomorphism with respect to the natural inner products
$$\langle\chi^{\blam},\chi^{\bdeta}\rangle=\delta_{\blam\bdeta} \qquad \text{and}\qquad \langle s_{\blam},s_{\bdeta}\rangle=\delta_{\blam\bdeta}.$$
\end{theorem}

In the following change of basis equations, (\ref{SCHURtoHALL}) follows from Theorem \ref{CharacteristicMap}, (\ref{SCHURtoPOWER}) follows from (\ref{SchurToPower}), and (\ref{POWERtoHALL}) follows from \cite[Theorem 4.2]{TV05}.  
\begin{align}
(-1)^{\lfloor k/2\rfloor+n(\blam)}s_{\blam} &= \sum_{\bmu\in \cP^\Phi_k}  \chi^{\blam}(\bmu) P_{\bmu} & & \text{for $\blam\in \cP_k^\Theta$,}\label{SCHURtoHALL}\\
s_{\blam} &= \sum_{\bnu\in \cP^\Theta_k\atop \blam_s=\bnu_s}\biggl(\prod_{\vphi\in \Theta} \frac{\omega^{\blam^{(\vphi)}}(\bnu^{(\vphi)})}{z_{\bnu^{(\vphi)}}}\biggr) p_{\bnu} & & \text{for $\blam\in \cP_k^\Theta$,}\label{SCHURtoPOWER}\\
(-1)^{k-\ell(\bnu)}p_{\bnu} &= \sum_{\bmu\in \cP^\Phi_k}\biggl(\sum_{t\in T_{\bnu_u}\atop \btau_\Phi(t)_s=\bmu_s}  \theta(t)Q_{\btau_\Phi(t)}^{\bmu}(-q)\biggr) P_{\bmu} & & \text{for $\bnu\in \cP_k^\Theta$, $\btau_{\Theta}(\theta) = \bnu$}.\label{POWERtoHALL}
\end{align}

\section{Gelfand-Graev characters on arbitrary elements}\label{SectionGelfandGraev}

\subsection{$G_n=\GL(n,\FF_q)$ notation}

In this Section \ref{SectionGelfandGraev}, let
$$\tilde\Phi=\{\tilde{F}\text{-orbits in $\bar{\FF}_q$}\} \qquad \text{and}\qquad  \tilde\Theta  = \{\tilde{F}\text{-orbits in $\bar{\FF}_q^*$}\}.$$
The same set-up of sections \ref{multipartitions} and \ref{UnitaryCharacteristic} gives a characteristic map for $G_n=\GL(n,\FF_q)$ by replacing $\Phi$ by $\tilde\Phi$, $\Theta$ by $\tilde\Theta$, $-q$ by $q$, $T_{(k)}$ by $\tilde{T}_{(k)}$, and $(-1)^{\lfloor n/2\rfloor+n(\blam)} s_{\blam}$ by $s_{\blam}$.  With the exception of the Deligne-Lusztig characters (which follows from the parallel argument of \cite[Theorem 4.2]{TV05}), this can be found in \cite[Chapter IV]{Mac}.

\subsection{The Gelfand-Graev character}

We will use $U_n'=\GL(n,\bar{\FF}_q)^{F'}$ (see (\ref{GroupDefinitions})) to give an explicit description of the Gelfand-Graev character.  For a more general description see \cite{dignemichel}, for example. 

For $1\leq i< j\leq n$ and $t\in \FF_q$, let $x_{ij}(t)$ denote the matrix with ones on the diagonal, $t$ in the $i$th row and $j$th column, and zeroes elsewhere.  Let
\begin{align*} 
u_{ij}(t) &= x_{ij}(t)x_{n+1-j,n+1-i}(-t^q) & & \text{for $1\leq i<j\leq \lfloor n/2\rfloor$, $t\in \FF_{q^2}$,}\\
u_{i,n+1-j}(t) &=  x_{i,n+1-j}(t)x_{j,n+1-i}(-t^q) & & \text{for $1\leq i<j\leq \lfloor n/2\rfloor$, $t\in \FF_{q^2}$,}
\end{align*}
and for $1\leq k\leq \lfloor n/2\rfloor$, and $t,a,b\in \FF_{q^2}$, let

\begin{align*}
u_k(a)&=x_{k,n+1-k}(a) & & \text{for $n$ even, and $a^q+a=0$,}\\
u_k(a,b)&=x_{\lceil n/2\rceil,n+1-k}(-a^q)x_{k,n+1-k}(b) x_{k,\lceil n/2\rceil}(a) & & \text{for $n$ odd, and $a^{q+1}+b+b^q=0$.}
\end{align*}

\noindent\textbf{Examples.}  In $U_4'$, we have
$$u_{12}(t)=\left(\begin{smallmatrix} 1 & t & 0 & 0\\ 0 & 1 & 0 & 0\\ 0 & 0 & 1 & -t^q\\ 0 & 0 & 0 & 1\end{smallmatrix}\right), \quad u_{13}(t)=\left(\begin{smallmatrix} 1 & 0 & t & 0\\ 0 & 1 & 0 & -t^q\\ 0 & 0 & 1 & 0\\ 0 & 0 & 0 & 1\end{smallmatrix}\right),\quad u_1(a)=\left(\begin{smallmatrix} 1 & 0 & 0 & a\\ 0 & 1 & 0 & 0\\ 0 & 0 & 1 & 0\\ 0 & 0 & 0 & 1\end{smallmatrix}\right), \quad u_2(a)=\left(\begin{smallmatrix} 1 & 0 & 0 & 0\\ 0 & 1 & a & 0\\ 0 & 0 & 1 & 0\\ 0 & 0 & 0 & 1\end{smallmatrix}\right),$$
where $a^q+a=0$.  In $U_5'$, we have
$$u_{12}(t)=\left(\begin{smallmatrix} 1 & t & 0 & 0 & 0\\ 0 & 1 & 0 & 0 & 0\\ 0 & 0 & 1 & 0 & 0\\ 0 & 0 & 0 & 1 & -t^q\\ 0 & 0 & 0 & 0 & 1 \end{smallmatrix}\right),\qquad u_{14}(t)=\left(\begin{smallmatrix} 1 & 0 & 0 & t & 0\\ 0 & 1 & 0 & 0 & -t^q\\ 0 & 0 & 1 & 0 & 0\\ 0 & 0 & 0 & 1 & 0\\ 0 & 0 & 0 & 0 & 1 \end{smallmatrix}\right),$$
$$u_1(a,b)=\left(\begin{smallmatrix} 1 & 0 & a & 0 & b\\ 0 & 1 & 0 & 0 & 0\\ 0 & 0 & 1 & 0 & -a^q\\ 0 & 0 & 0 & 1 & 0\\ 0 & 0 & 0 & 0 & 1 \end{smallmatrix}\right),\qquad u_2(a,b)=\left(\begin{smallmatrix} 1 & 0 & 0 & 0 & 0\\ 0 & 1 & a & b & 0\\ 0 & 0 & 1 & -a^q & 0\\ 0 & 0 & 0 & 1 & 0\\ 0 & 0 & 0 & 0 & 1 \end{smallmatrix}\right),$$
where $a^{q+1}+b+b^q=0$.

For $1\leq i<j\leq \lfloor n/2\rfloor$, and $1\leq k\leq \lfloor n/2\rfloor$, define the one-parameter subgroups
\begin{align*}
\X_{ij}&=\{u_{ij}(t)\ \mid\ t\in \FF_{q^2}\}\cong \FF_{q^2}^+,\\
\X_{i,n+1-j}&=\{u_{i,n+1-j}(t)\ \mid\ t\in \FF_{q^2}\}\cong \FF_{q^2}^+,\\
\X_k &=\left\{\begin{array}{ll}\{ u_k(t)\ \mid\ t\in \FF_{q^2}, t^q+q=0\}&\text{if $n$ is even,}\\ \{u_k(a,b)\ \mid\ a,b\in \FF_{q^2}, a^{q+1}+b+b^q=0\} & \text{if $n$ is odd.}\end{array}\right.
\end{align*}
so that
$$B_n^<=\langle \X_{ij}, \X_{i,n-j}, \X_k\ \mid\ 1\leq i<j\leq \lfloor n/2\rfloor, 1\leq k\leq \lfloor n/2\rfloor\rangle\subseteq U_n'$$
is the subgroup of $U_n'$ of upper-triangular matrices with ones on the diagonal.  Noting that 
$$\X_k/[\X_k,\X_k]\cong \FF_q^+,$$
a direct calculation gives 
$$U_n'/[U_n',U_n']\cong \X_{12}\times \X_{23}\times\ldots\times \X_{\lfloor n/2\rfloor-1, \lfloor n/2\rfloor}\times \X_{\lfloor n/2\rfloor}\cong (\FF_{q^2}^+)^{\lfloor n/2\rfloor-1}\times \FF_{q}^+.$$

Similarly, let
$$\tilde{B}_n^<=\langle x_{ij}(t)\ \mid\ 1\leq i<j\leq \ell,t\in \FF_{q}\rangle\subseteq G_n$$
be the subgroup of unipotent upper-triangular matrices in $G_n$.

Fix a homomorphism $\psi:\FF_{q^2}^+\rightarrow \CC^\times$ of the additive group of the field such that for all $1\leq k\leq \lfloor n/2\rfloor$, $\psi$ is nontrivial on $\X_k/[\X_k,\X_k]$.  Define the homomorphism $\psi_{(n)}: B_n^<\rightarrow \CC$ by
$$\psi_{(n)}\bigg|_{\X_\alpha/[\X_\alpha,\X_\alpha]}=\left\{\begin{array}{ll} \psi & \text{if $\alpha=(i,i+1)$, $1\leq i< \lfloor n/2\rfloor$, or if $\alpha = \lfloor n/2 \rfloor$,}\\ 1 &\text{otherwise.}\end{array}\right.  $$
The \emph{Gelfand-Graev character} of $U_n'$ is
$$\Gamma_n'=\Ind_{B_n^<}^{U_n'}(\psi_{(n)}).$$
Recall that $U_n'$ is conjugate to $U_n$ in $\bar{G}_n$.  If $U_n'=yU_ny^{-1}$, then let
$$\Gamma_n=\Ind_{y^{-1}B_n^< y}^{U_n}(y^{-1} \psi_{(n)}y).$$

Similarly, the Gelfand-Graev character $\tilde\Gamma_{(n)}$ of $G_n$ is
$$\tilde\Gamma_{(n)}=\Ind_{\tilde{B}_n^<}^{G_n}(\tilde\psi_{(n)}).$$
where $\tilde\psi_{(n)}:\tilde{B}_n^<\rightarrow\CC$ is given by
$$\tilde\psi_{(n)}(x_{ij}(t))=\delta_{j,i+1} \psi(t).$$

It is well-known that the Gelfand-Graev character has a multiplicity free decomposition into irreducible characters \cite{St67,Yo68,Yo69}.  The following explicit decompositions essentially follow from \cite{dellusz}.  Specific proofs are given in \cite{zel} in the $G_n$ case and in \cite{ohm} in the $U_n$ case.

\begin{theorem} \label{GGDecomposition}
Let $\hgt(\blam) = \mathrm{max} \{ \ell(\blam^{(\varphi)})\}.$  Then
$$\Gamma_{(n)} = \sum_{\blam\in \cP_n^\Theta \atop \hgt(\blam) = 1} \chi^{\blam} \qquad\text{and}\qquad
\tilde\Gamma_{(n)}=\sum_{\blam\in \cP_n^{\tilde\Theta} \atop \hgt(\blam) = 1} \chi^{\blam}.$$
\end{theorem}

\subsection{The character values of the Gelfand-Graev character}

A \emph{unipotent conjugacy class} $K^{\bmu}$ of $U_n$ or $G_n$ is a conjugacy class that satisfies
$$\bmu_u=\bmu.$$
The unipotent conjugacy classes of $U_n$ and $G_n$ are thus parameterized by partitions $\mu$ of $n$.

\begin{lemma}\label{GGGreenFormula}  \hfill
\begin{enumerate}
\item[(a)] Let $\bmu \in \cP_n^{\Theta}$, $\bmu_u^{(\{1\})} = \mu$.  Then 
$$\dd{\Gamma_{(n)}(\bmu)=\left\{\begin{array}{ll} \sum_{\nu\in \cP_n}\frac{(-1)^{n+\lfloor n/2\rfloor-\ell(\nu)}}{z_\nu} |T_\nu|  Q_{\nu}^{\mu}(-q) & \text{if $\bmu$ is unipotent, } \\ 0 & \text{otherwise.}\end{array}\right.}$$
\item[(b)] Let $\bmu \in \cP_n^{\tilde\Theta}$, $\bmu_u^{(\{1\})} = \mu$. Then $$\dd{\tilde\Gamma_{(n)}(\bmu)=\left\{\begin{array}{ll} \sum_{\nu\in \cP_n}\frac{(-1)^{n-\ell(\nu)}}{z_\nu} |\tilde{T}_\nu|  Q_{\nu}^{\mu}(q) & \text{if $\bmu$ is unipotent, } \\ 0 & \text{otherwise.}\end{array}\right.}$$
\end{enumerate}
\end{lemma}

\begin{proof}  
Note that if $\hgt(\blam)\leq 1$, then $n(\blam)=0$.  Thus, by applying the characteristic map and (\ref{SCHURtoPOWER}) to Theorem \ref{GGDecomposition},
\begin{equation*}
\ch(\Gamma_{(n)}) 
=(-1)^{\lfloor n/2\rfloor} \sum_{\blam\in \cP_n^\Theta\atop \hgt(\blam)\leq 1} \sum_{\bnu\in \cP_n^\Theta\atop \bnu_s=\blam_s} \biggl(\prod_{\vphi\in \Theta} \frac{\omega^{\blam^{(\vphi)}}(\bnu^{(\vphi)})}{z_{\bnu^{(\vphi)}}}\biggr) p_{\bnu}.
\end{equation*}
Since $\hgt(\blam)\leq 1$, $\omega^{\blam^{(\vphi)}}$ is the trivial character for all $\vphi\in \Theta$.   Thus, the summand is independent of $\blam$, and
$$\ch(\Gamma_{(n)}) =(-1)^{\lfloor n/2\rfloor} \sum_{\bnu\in \cP_n^\Theta} \biggl(\prod_{\vphi\in \Theta} z_{\bnu^{(\vphi)}}^{-1}\biggr) p_{\bnu}.$$
By (\ref{TorusMultiplicities}),
\begin{align*}
&\ch(\Gamma_{(n)})=\\
&=(-1)^{\lfloor n/2\rfloor} \sum_{\bnu\in \cP_n^\Theta} \Biggl(\prod_{\vphi \in \Theta} |\vphi|^{\ell(\bnu^{(\vphi)})} \prod_{i \ge 1} \frac{\big(m_i(\bnu_u^{(\{1\})}) \big)!}{\prod_{\vphi \in \Theta}(m_{i/|\vphi|}(\bnu^{(\vphi)}))!}\Biggl)^{-1}  \sum_{\theta\in \Hom(T_{\bnu_u},\CC^\times) \atop \btau_{\Theta}(\theta) = \bnu}  \biggl(\prod_{\vphi\in \Theta} z_{\bnu^{(\vphi)}}^{-1}\biggr) p_{\bnu}\\
&=(-1)^{\lfloor n/2\rfloor} \sum_{\bnu\in \cP_n^\Theta}  \sum_{\theta\in \Hom(T_{\bnu_u},\CC^\times) \atop \btau_{\Theta}(\theta) = \bnu}  z_{\bnu_u}^{-1} p_{\bnu}\\
&=(-1)^{\lfloor n/2\rfloor} \sum_{\nu\in \cP_n}  \sum_{\theta\in \Hom(T_{\nu},\CC^\times)}  z_{\nu}^{-1} p_{\btau_\Theta(\theta)}.
\end{align*}
The change of basis (\ref{POWERtoHALL}) gives  
\begin{align*}
\ch(\Gamma_{(n)}) &= (-1)^{\lfloor n/2\rfloor} \sum_{\nu\in \cP_n}  \sum_{\theta\in \Hom(T_{\nu},\CC^\times)} \frac{(-1)^{n-\ell(\nu)}}{z_\nu}  \sum_{\bmu\in \P^\Phi_n}\sum_{t\in T_{\nu}\atop \btau_\Phi(t)_s=\bmu_s} \theta(t) Q_{\btau_\Phi(t)}^{\bmu}(-q) P_{\bmu}\\
&= (-1)^{\lfloor n/2\rfloor} \sum_{\bmu\in \P^\Phi_n} \sum_{\nu\in \cP_n}  \frac{(-1)^{n-\ell(\nu)}}{z_\nu} \sum_{t\in T_{\nu}\atop \btau_\Phi(t)_s=\bmu_s} \sum_{\theta\in \Hom(T_{\nu},\CC^\times)} \theta(t) Q_{\btau_\Phi(t)}^{\bmu}(-q) P_{\bmu}.
\end{align*}
By the orthogonality of characters of $T_\nu$, the inner-most sum is equal to zero for all $t\neq 1$.  If $t=1$, then $\btau_\Phi(1,1,\ldots,1)^{(f)}=\emptyset$ for $f\neq \{1\}$ and $\btau_\Phi(1,1,\ldots,1)^{(\{1\})}=\nu$.  Thus, 
$$\ch(\Gamma_{(n)}) = (-1)^{\lfloor n/2\rfloor}\sum_{\bmu\in \cP^\Phi_n\atop \bmu_u=\bmu}\sum_{\nu\in \cP_n}\frac{(-1)^{n-\ell(\nu)}}{z_\nu} |T_\nu|  Q_{\nu}^{\bmu_u}(-q) P_{\bmu},$$
and in particular, if $\bmu_u^{(\{1\})} = \mu$, 
$$\Gamma_{(n)}(\bmu)=\left\{\begin{array}{ll} \sum_{\nu\in \cP_n}\frac{(-1)^{n+\lfloor n/2\rfloor-\ell(\nu)}}{z_\nu} |T_\nu|  Q_{\nu}^{\mu}(-q) & \text{if $\bmu$ is unipotent,} \\ 0 & \text{otherwise.}\end{array}\right.$$

(b) The proof is similar to (a), just using the $G_n$ characteristic map.
\end{proof}

The values of the Gelfand-Graev character of the finite general linear group are well-known.  An elementary proof of the following Theorem is given in \cite{howzwor}.

\begin{theorem} \label{GLGGFormula} Let $\bmu\in \cP^{\tilde\Phi}_n$ with $\mu=\bmu_u^{(\{1\})}$.  Then
$$\tilde\Gamma_{(n)}(\bmu)=\left\{\begin{array}{ll} (-1)^{n- \ell(\mu)} \prod_{i = 1}^{\ell(\mu)} \bigl( q^i - 1 \bigr) & \text{if $\bmu$ is unipotent, } \\ 0 & \text{otherwise.}\end{array}\right.$$
\end{theorem}

We may now apply Theorem \ref{GLGGFormula} and Lemma \ref{GGGreenFormula} to give the values of the Gelfand-Graev character of $U_n$.

\begin{corollary} \label{UGGFormula}
Let $\bmu\in \cP^\Phi_n$ with $\mu=\bmu_u^{(\{1\})}$.  Then
$$\Gamma_{(n)}(\bmu)=\left\{\begin{array}{ll} (-1)^{\lfloor n/2 \rfloor - \ell(\mu)} \prod_{i = 1}^{\ell(\mu)} \bigl( (-q)^i - 1 \bigr) & \text{if $\bmu$ is unipotent,}\\ 0 & \text{otherwise.}\end{array}\right.$$
\end{corollary}

\begin{proof}
Combine Lemma \ref{GGGreenFormula} (b) with Theorem \ref{GLGGFormula} to get
$$(-1)^{\ell(\mu)} \prod_{i = 1}^{\ell(\mu)} \bigl(q^i - 1 \bigr)= \sum_{\nu\in \cP_n}\frac{(-1)^{\ell(\nu)}}{z_\nu} |\tilde{T}_\nu|  Q_{\nu}^{\mu}(q),$$
which implies
$$ \prod_{i = 1}^{\ell(\mu)} \bigl(1-q^i \bigr)= \sum_{\nu\in \cP_n}\frac{1}{z_\nu}\prod_{i=1}^{\ell(\nu)} (1-q^{\nu_i})  Q_{\nu}^\mu(q).$$
Make the substitution $q\mapsto -q$ to get
$$ \prod_{i = 1}^{\ell(\mu)} \bigl(1-(-q)^i \bigr)= \sum_{\nu\in \cP_n}\frac{1}{z_\nu}\prod_{i=1}^{\ell(\nu)} (1-(-q)^{\nu_i})  Q_{\nu}^\mu(-q),$$
which implies
$$(-1)^{\lfloor n/2\rfloor+\ell(\mu)}\prod_{i = 1}^{\ell(\mu)} \bigl((-q)^i-1 \bigr)= \sum_{\nu\in \cP_n}\frac{(-1)^{\lfloor n/2\rfloor+|\nu|-\ell(\nu)}}{z_\nu}|T_\nu|  Q_{\nu}^\mu(-q).$$
Apply this last identity to Lemma \ref{GGGreenFormula} (a) to obtain the desired result.
\end{proof}

\section{Degenerate Gelfand-Graev characters}\label{SectionDegenerateGelfandGraev}

\subsection{$G_n=\GL(n,\FF_{q^2})$ notation (different from Section \ref{SectionGelfandGraev})}

In this Section \ref{SectionDegenerateGelfandGraev}, let $G_n={\rm GL}(n, \FF_{q^2})$, 
$$\tilde\Phi=\{F^2\text{-orbits of $\bar{\FF}_q^\times$}\}, \qquad \text{and}\qquad  \tilde\Theta  = \{F^2\text{-orbits of $\bar{\FF}_q^*$}\}.$$
Note that now $G_n^{F} = U_n$ and $G_n^{F^{\prime}} = U_n^{\prime}$.
The same set-up of sections \ref{multipartitions} and \ref{UnitaryCharacteristic} gives a characteristic map for $G_{n}$ by replacing $\Phi$ by $\tilde\Phi$, $\Theta$ by $\tilde\Theta$, $-q$ by $q$, $T_{(k)}$ by $T_{(2k)}$, and $(-1)^{\lfloor n/2\rfloor+n(\lambda)} s_{\blam}$ by $s_{\blam}$.

\subsection{Degenerate Gelfand-Graev characters}

Let $(k,\nu)$ be a pair such that $\nu\vdash \frac{n-k}{2}\in \ZZ_{\geq 0}$, and let
$$\nu_{\leq}=(\nu_{\leq 1}, \nu_{\leq 2},\ldots, \nu_{\leq \ell}), \quad \text{where}\quad \nu_{\leq j}=\nu_1+\nu_2+\cdots +\nu_j.$$
Then the map $\psi_{(k,\nu)} :B^<_n\rightarrow \CC^\times$, given by
$$\psi_{(k,\nu)}\bigg|_{\X_\alpha/[\X_\alpha,\X_\alpha]}=\left\{\begin{array}{ll} \psi & \text{if $\alpha=(i,i+1)$, $1\leq i< \lfloor n/2\rfloor$, and  $i\notin\nu_{\leq}$},\\ 
 \psi & \text{if $\alpha=\lfloor n/2\rfloor$ and  $\lfloor n/2\rfloor \notin\nu_{\leq}$},\\ 
1 & \text{otherwise,}\end{array}\right.$$
is a linear character of $U_n^\prime$.  Note that $\psi_{(\lceil n/2\rceil-\lfloor n/2\rfloor, (1^{\lfloor n/2\rfloor}))}$ is the trivial character and $\psi_{(n,\emptyset)}=\psi_{(n)}$ of Section \ref{SectionGelfandGraev}.

The \emph{degenerate Gelfand-Graev character} $\Gamma_{(k,\nu)}$ is
$$\Gamma_{(k,\nu)}=\Ind_{B^<_n}^{U_n^\prime}(\psi_{(k,\nu)})\cong\Ind_{yB^<_ny^{-1}}^{U_n}(y\psi_{(k,\nu)} y^{-1}),$$
where $y$ is an element of $\bar{G}_n$ such that $yU_n^{\prime} y^{-1} = U_n$.
In particular, the Gelfand-Graev character is $\Gamma_{(n,\emptyset)}$.

Let
$$L^\prime_{(k,\nu)}=\langle L_{k},  L_\nu^{(1)}, L_\nu^{(2)}, \cdots, L_\nu^{(\ell)} \rangle,$$
where
\begin{align*} 
 L_k&=\langle \X_{ij},\X_{i,n+1-j},\X_r\ \mid\ \, |\nu| < i<j\leq |\nu|+k,|\nu|\leq r\leq |\nu|+k\rangle\cong {\rm U}(k, \FF_{q^2})\\
L_{\nu}^{(r)}&=\langle \X_{ij}\ \mid\ \nu_{\leq r-1} \leq i<j\leq \nu_{\leq r}\rangle\cong {\rm GL}(\nu_r, \FF_{q^2}).
\end{align*}
Then
$$L^\prime_{(k,\nu)}\cong {\rm U}(k, \FF_{q^2})\oplus {\rm GL}(\nu_1, \FF_{q^2})\oplus\cdots \oplus {\rm GL}(\nu_\ell, \FF_{q^2})$$
is a maximally split Levi subgroup of $U_n^\prime$.  For example, if $n=9$, $k=3$, and $\nu=(2,1)$, then 
$$L_{(k,\nu)}'=\left\{\left(\begin{array}{ccccc} A & 0 & 0 & 0 & 0\\  0 & B & 0 & 0 & 0\\ 0 & 0 & C & 0 & 0\\ 0 & 0 & 0 & F'(B) & 0\\ 0 & 0 & 0 & 0 & F'(C)\end{array}\right) \bigg| A\in \GL(2,\FF_{q^2}), B\in \GL(1,\FF_{q^2}), C\in U(3,\FF_{q^2})\right\}.$$

Note that since $L_\nu^{(i)}\subseteq U_{2\nu_i}^{\prime}\cong U_{2\nu_i}$, the Levi subgroup
$$U_{(k,\nu)}=U_{k}\oplus U_{2\nu_1}\oplus U_{2\nu_2}\oplus\cdots\oplus U_{2\nu_\ell}\subseteq U_n$$
contains a Levi subgroup $L=U_k\oplus L_1\oplus\cdots \oplus L_\ell$ with $L_i\subseteq U_{2\nu_i}$ such that $L\cong L^\prime_{(k,\nu)}$.

\begin{proposition} \label{DGGProductDecomposition} Let $(k,\nu)$ be such that $\nu\vdash \frac{n-k}{2}\in \ZZ_{\geq 0}$.  Then 
$$\ch(\Gamma_{(k,\nu)})=\ch\bigl(\Gamma_{(k)}\bigr)\ch\biggl(R_{G_{\nu_1}}^{U_{2\nu_1}}(\tilde\Gamma_{(\nu_1)})\biggr)\ch\biggl(R_{G_{\nu_2}}^{U_{2\nu_2}}(\tilde\Gamma_{(\nu_2)})\biggr)\cdots \ch\biggl(R_{G_{\nu_\ell}}^{U_{2\nu_\ell}}(\tilde\Gamma_{(\nu_\ell)})\biggr),$$
where $\tilde\Gamma_{(m)}$ is the Gelfand-Graev character of $G_m=\GL(m, \FF_{q^2})$.
\end{proposition}

This proposition is a consequence of Theorem \ref{CharacteristicMap} and the following lemma.

\begin{lemma} Let $(k,\nu)$ be such that $\nu\vdash \frac{n-k}{2}\in \ZZ_{\geq 0}$.  Then 
$$\Gamma_{(k,\nu)}\cong R_{U_{(k,\nu)}}^{U_n}\big(\Gamma_{(k)}\otimes R_{L_1}^{U_{2\nu_1}}(\tilde\Gamma_{(\nu_1)})\otimes\cdots \otimes R_{L_\ell}^{U_{2\nu_\ell}}(\tilde\Gamma_{(\nu_\ell)})\big).$$
\end{lemma}

\begin{proof}
Since $L^\prime_{(k,\nu)}$ is maximally split, 
$$\Ind_{yB^<_ny^{-1}}^{U_n}(y\psi_{(k,\nu)}y^{-1})\cong\Ind_{B^<_n}^{U_n^{\prime}}(\psi_{(k,\nu)})\cong \Indf_{L^\prime_{(k,\nu)}}^{U_n^{\prime}}(\Gamma_{(k)}\otimes\tilde\Gamma_{(\nu_1)}\otimes\cdots \otimes \tilde\Gamma_{(\nu_\ell)}).$$
where $\Indf_{L}^{G}$ is Harish-Chandra induction.  However,
\begin{align*}
\Indf_{L^\prime_{(k,\nu)}}^{U_n^{\prime}}(\Gamma_{(k)}\otimes\tilde\Gamma_{(\nu_1)}\otimes\cdots \otimes \tilde\Gamma_{(\nu_\ell)})&=R_{L^\prime_{(k,\nu)}}^{U_n^{\prime}}(\Gamma_{(k)}\otimes\tilde\Gamma_{(\nu_1)}\otimes\cdots \otimes \tilde\Gamma_{(\nu_\ell)}),\\
\cong R_{L}^{U_n}(\Gamma_{(k)}\otimes\tilde\Gamma_{(\nu_1)}\otimes\cdots \otimes \tilde\Gamma_{(\nu_\ell)}).
\end{align*}
By transitivity of Deligne-Lusztig induction, we now have
$$\Ind_{yB^<_ny^{-1}}^{U_n}(y\psi_{(k,\nu)}y^{-1})\cong R_{U_{(k,\nu)}}^{U_n}\big(\Gamma_{(k)}\otimes R_{L_1}^{U_{2\nu_1}}(\tilde\Gamma_{(\nu_1)})\otimes\cdots \otimes R_{L_\ell}^{U_{2\nu_\ell}}(\tilde\Gamma_{(\nu_\ell)})\big). \qedhere $$ 
\end{proof}

\subsection{Symplectic tableaux and domino tableaux combinatorics}

Augment the nonnegative integers by symbols $\{\bar{i}\ \mid\ i\in \ZZ_{>0}\}$, so that we have
$$\{0,\bar{1},1,\bar{2},2,\bar{3},3, \ldots\},$$
and order this set by  $i-1<\bar{i}<i<\overline{i+1}$.  Alternatively, one  could identify this augmented set with $\frac{1}{2}\ZZ_{\geq 0}$ by  $\bar{i}=i-\frac{1}{2}$.

Let $\lambda=(\lambda_1,\lambda_2,\ldots, \lambda_r)$ be a partition of $n$ and $(m_0, m_1,m_2,\ldots,m_\ell)$ be a sequence of nonnegative integers that sum to $n$ with $m_0\leq \lambda_1$.   A \emph{symplectic tableau} $Q$ of \emph{shape} $\lambda/(m_0)$ and \emph{weight} $(m_0,m_1,\ldots, m_\ell)$ is a column strict filling of the boxes of $\lambda$ by symbols 
$$\{0,\bar{1},1,\bar{2},2,\ldots, \bar{\ell},\ell\},$$
such that 
$$m_i=\left\{\begin{array}{ll} \text{number of $0$'s in $Q$} &\text{if $i=0$,}\\
 \text{number of $\bar{i}$'s $+$ number of $i$'s in $Q$} & \text{if $i>0$.}\end{array}\right.$$
We write $\sh(Q)=\lambda/(m_0)$ and $\wt(Q)=(m_0,m_1,\ldots,m_\ell)$.  
For example, if 
$$Q
=
\xy<0cm,.5cm>\xymatrix@R=.35cm@C=.35cm{
*={} & *={} \ar @{-} [l] & *={} \ar @{-} [l] & *={} \ar @{-} [l] & *={} \ar @{-} [l] & *={} \ar @{-} [l]  \\
*={} \ar @{-} [u] &   *={} \ar @{-} [l] \ar@{-} [u]  \ar @{} [ul]|{\ss 0} &   *={} \ar @{-} [l] \ar@{-} [u]  \ar @{} [ul]|{\ss 0}  &   *={} \ar @{-} [l] \ar@{-} [u]  \ar @{} [ul]|{\ss \bar{1}} &   *={} \ar @{-} [l] \ar@{-} [u]  \ar @{} [ul]|{\ss 1} &   *={} \ar @{-} [l] \ar@{-} [u]  \ar @{} [ul]|{\ss \bar{4}} \\ 
*={} \ar @{-} [u] &   *={} \ar @{-} [l] \ar@{-} [u]  \ar @{} [ul]|{\ss 1} &   *={} \ar @{-} [l] \ar@{-} [u]  \ar @{} [ul]|{\ss \bar{2}}&   *={} \ar @{-} [l] \ar@{-} [u]  \ar @{} [ul]|{\ss \bar{2}} \\
*={} \ar @{-} [u] &   *={} \ar @{-} [l] \ar@{-} [u]  \ar @{} [ul]|{\ss \bar{3}}\\
*={} \ar @{-} [u]  &   *={} \ar @{-} [l] \ar@{-} [u] \ar @{} [ul]|{\ss 3} }\endxy,\quad\text{then}\quad
\sh(Q)
=
\xy<0cm,.5cm>\xymatrix@R=.3cm@C=.3cm{
*={} & *={}   & *={} & *={} \ar @{-} [l] & *={} \ar @{-} [l] & *={} \ar @{-} [l]  \\
*={}   &   *={} \ar @{-} [l]   &   *={} \ar @{-} [l] \ar@{-} [u]  &   *={} \ar @{-} [l] \ar@{-} [u] &   *={} \ar @{-} [l] \ar@{-} [u]  &   *={} \ar @{-} [l] \ar@{-} [u] \\ 
*={} \ar @{-} [u] &   *={} \ar @{-} [l] \ar@{-} [u] &   *={} \ar @{-} [l] \ar@{-} [u] &   *={} \ar @{-} [l] \ar@{-} [u] \\
*={} \ar @{-} [u] &   *={} \ar @{-} [l] \ar@{-} [u] \\
*={} \ar @{-} [u] &   *={} \ar @{-} [l] \ar@{-} [u] }\endxy
\quad\text{and}\quad
\wt(Q)=(2,3,2,2,1)
$$

Let
\begin{equation}\label{SymplecticTableauxSet}
\cT_{(m_0,m_1,\ldots,m_\ell)}^\lambda=\left\{\begin{array}{c}\text{symplectic tableaux of shape $\lambda/(m_0)$}\\ \text{and weight $(m_0,m_1,\ldots, m_\ell)$}\end{array}\right\}.
\end{equation}

A \emph{tiling of $\lambda$ by dominoes} is a partition of the boxes of $\lambda$ into pairs of adjacent boxes. For example, if 
$$\lambda=   
\xy<0cm,.5cm>\xymatrix@R=.3cm@C=.3cm{
*={} & *={} \ar @{-} [l] & *={} \ar @{-} [l] & *={} \ar @{-} [l] & *={} \ar @{-} [l] & *={} \ar @{-} [l]  \\
*={} \ar @{-} [u] &   *={} \ar @{-} [l] \ar@{-} [u] &   *={} \ar @{-} [l] \ar@{-} [u]  &   *={} \ar @{-} [l] \ar@{-} [u] &   *={} \ar @{-} [l] \ar@{-} [u]  &   *={} \ar @{-} [l] \ar@{-} [u] \\ 
*={} \ar @{-} [u] &   *={} \ar @{-} [l] \ar@{-} [u] &   *={} \ar @{-} [l] \ar@{-} [u] &   *={} \ar @{-} [l] \ar@{-} [u] \\
*={} \ar @{-} [u] &   *={} \ar @{-} [l] \ar@{-} [u] \\
*={} \ar @{-} [u] &   *={} \ar @{-} [l] \ar@{-} [u] }\endxy,
\quad \text{then}\quad 
\xy<0cm,.5cm>\xymatrix@R=.3cm@C=.3cm{
*={} & *={} \ar @{-} [l] & *={} \ar @{-} [l] & *={} \ar @{-} [l] & *={} \ar @{-} [l] & *={} \ar @{-} [l]  \\
*={} \ar @{-} [u] &   *={} \ar @{-} [l] &   *={} \ar @{-} [l] \ar@{-} [u]  &   *={} \ar@{-} [u] &   *={} \ar @{-} [l]  &   *={} \ar @{-} [l] \ar@{-} [u] \\ 
*={} \ar @{-} [u] &   *={} \ar @{-} [l]  &   *={} \ar @{-} [l] \ar@{-} [u] &   *={} \ar @{-} [l] \ar@{-} [u] \\
*={} \ar @{-} [u] &   *={} \ar@{-} [u] \\
*={} \ar @{-} [u] &   *={} \ar @{-} [l] \ar@{-} [u] }\endxy$$
is a tiling of $\lambda$ by dominoes.

Let $(m_0,m_1,\ldots, m_\ell)$ be a sequence of nonnegative integers such that $m_0\leq \lambda_1$ and $|\lambda|=m_0+2(m_1+\cdots+m_\ell)$.  A \emph{domino tableau} $Q$ of \emph{shape} $\lambda/(m_0)=\sh(Q)$ and \emph{weight} $(m_0,m_1,\ldots, m_\ell)= \wt(Q)$ is a column strict filling of a tiling of the shape $\lambda/(m_0)$ by dominoes, where we put $0$'s in the non-tiled boxes of $\lambda$, and where $m_i$ is the number of $i$'s which appear.  For example, if
$$Q=\xy<0cm,.5cm>\xymatrix@R=.3cm@C=.3cm{
*={} & *={} \ar @{-} [l] & *={} \ar @{-} [l] & *={} \ar @{-} [l] & *={} \ar @{-} [l] & *={} \ar @{-} [l]  \\
*={} \ar @{-} [u] &   *={} \ar @{-} [l] \ar @{-} [u]  \ar @{} [ul]|{0}  &   *={} \ar @{-} [l] \ar@{-} [u] \ar @{} [ul]|{0} &   *={} \ar@{-} [u] &   *={} \ar @{-} [l]  &   *={} \ar @{-} [l] \ar@{-} [u] \ar @{} [ull]|{3} \\ 
*={} \ar @{-} [u] &   *={} \ar @{-} [l]  &   *={} \ar @{-} [l] \ar@{-} [u] \ar @{} [ull]|{1}  &   *={} \ar @{-} [l] \ar@{-} [u] \ar @{} [uul]|{3} \\
*={} \ar @{-} [u] &   *={} \ar@{-} [u] \\
*={} \ar @{-} [u] &   *={} \ar @{-} [l] \ar@{-} [u]\ar @{} [uul]|{2}  }\endxy, 
\quad\text{then} \quad 
\sh(Q)= \xy<0cm,.5cm>\xymatrix@R=.3cm@C=.3cm{
*={} & *={} & *={}& *={} \ar @{-} [l] & *={} \ar @{-} [l] & *={} \ar @{-} [l]  \\
*={} &   *={} \ar @{-} [l] &   *={} \ar @{-} [l] \ar@{-} [u]  &   *={} \ar @{-} [l] \ar@{-} [u] &   *={} \ar @{-} [l] \ar@{-} [u]  &   *={} \ar @{-} [l] \ar@{-} [u] \\ 
*={} \ar @{-} [u] &   *={} \ar @{-} [l] \ar@{-} [u] &   *={} \ar @{-} [l] \ar@{-} [u] &   *={} \ar @{-} [l] \ar@{-} [u] \\
*={} \ar @{-} [u] &   *={} \ar @{-} [l] \ar@{-} [u] \\
*={} \ar @{-} [u] &   *={} \ar @{-} [l] \ar@{-} [u] }\endxy
\quad\text{and} \quad
\wt(Q)=(2,1,1,2). 
$$

Let
\begin{equation}\label{DominoTableauxSet}
\cD_{(m_0,m_1,\ldots,m_\ell)}^\lambda=\left\{\begin{array}{c}\text{domino tableaux of shape $\lambda/(m_0)$}\\ \text{and weight $(m_0,m_1,\ldots, m_\ell)$}\end{array}\right\}.
\end{equation}
In the following Lemma, (a) is a straightforward use of the usual Pieri rule, and (b) is both similar to (and perhaps a special case of) \cite[Theorem 6.3]{LT00}, and also related to a Pieri formula in \cite{La05}.

\begin{lemma} \label{PieriLemma}  Let $(m_0,m_1,\ldots, m_\ell)$ be an $\ell+1$-tuple of nonnegative integers which sum to $n$.   Then
\begin{enumerate}
\item[(a)] 
${\dd s_{(m_0)}\prod_{r=1}^\ell \sum_{i=0}^{m_r} s_{(i)}s_{(m_r-i)}=\sum_{\lambda\in \cP_n}|\cT_{(m_0,m_1,\ldots,m_\ell)}^\lambda| s_\lambda},$
\item[(b)]
${\dd s_{(m_0)}\prod_{r=1}^\ell \sum_{i=0}^{2m_r}(-1)^i s_{(i)}s_{(2m_r-i)}=\sum_{\lambda\in \cP_{2n-m_0}}(-1)^{n(\lambda)} |\cD_{(m_0,m_1,\ldots,m_\ell)}^\lambda| s_\lambda.}$
\end{enumerate}
\end{lemma}

\begin{proof}
(a) Note that
$$s_{(m_0)}\prod_{r=1}^\ell \sum_{i=0}^{m_r} s_{(i)}s_{(m_r-i)}=\sum_{0\leq i_r\leq m_r\atop 1\leq r\leq \ell} s_{(m_0)}\prod_{r=1}^\ell s_{(i_r)}s_{(m_r-i_r)}.$$
Now repeated applications of Pieri's rule implies the result.

(b) Note that
$$s_{(m_0)}\prod_{r=1}^\ell \sum_{i=0}^{2m_r}(-1)^i s_{(i)}s_{(2m_r-i)}=\sum_{0\leq i_r\leq 2m_r\atop 1\leq r\leq \ell} (-1)^{i_1+\cdots+i_\ell} s_{(m_0)}\prod_{r=1}^{\ell} s_{(i_r)}s_{(2m_r-i_r)}.$$
By Pieri's rule,
$$\sum_{0\leq i_r\leq 2m_r\atop 1\leq r\leq \ell} s_{(m_0)}\prod_{r=1}^{\ell} s_{(i_r)}s_{(2m_r-i_r)}=\sum_{\lambda\in \cP_{2n-m_0}} \left(\begin{array}{c}  \text{Number of column strict fillings of $\lambda$ }\\  \text{using $m_0$ $0$'s, and for $r=1,2,\ldots,\ell$,}\\ \text{ using  $i_r$ $\bar{r}$'s and $(2m_r-i_r)$ $r$'s.}\end{array}\right)s_\lambda.$$
By observing that the sign counts the number of barred entries,
\begin{equation}\label{SignedEquation} s_{(m_0)}\prod_{r=1}^\ell \sum_{i=0}^{2m_r}(-1)^i s_{(i)}s_{(2m_r-i)}=\sum_{\lambda\in \cP_{2n-m_0}}\biggl(  \sum_{Q\in \cT_{(m_0,2m_1,\ldots,2m_\ell)}^\lambda} \hspace{-.7cm}(-1)^{\text{Number of barred entries in $Q$}}\biggr)s_\lambda.\end{equation}
We therefore need to determine the cancellations for a given shape $\lambda$.

Fix $r\in \{1,2,\ldots, \ell\}$ and $\lambda\in \cP$ such that $\cT_{(m_0,2m_1,\ldots,2m_\ell)}^\lambda\neq \emptyset$.  For a tableau $Q\in \cT_{(m_0,2m_1,\ldots,2m_\ell)}^\lambda$, let
\begin{align*}
Q_r&=\text{skew tableaux consisting of the boxes in $Q$ containing $\bar{r}$ or $r$},\\
\cS_Q^{(r)}&=\{\text{column strict fillings of $\sh(Q_r)$ by elements in $\{\bar{r},r\}$}\}.
\end{align*}
For example, if
$$Q=
\xy<0cm,.5cm>\xymatrix@R=.35cm@C=.35cm{
*={} & *={} \ar @{-} [l] & *={} \ar @{-} [l] & *={} \ar @{-} [l] & *={} \ar @{-} [l] & *={} \ar @{-} [l]  \\
*={} \ar @{-} [u] &   *={} \ar @{-} [l] \ar@{-} [u]  \ar @{} [ul]|{\ss 0} &   *={} \ar @{-} [l] \ar@{-} [u]  \ar @{} [ul]|{\ss 0}  &   *={} \ar @{-} [l] \ar@{-} [u]  \ar @{} [ul]|{\ss \bar{1}} &   *={} \ar @{-} [l] \ar@{-} [u]  \ar @{} [ul]|{\ss 1} &   *={} \ar @{-} [l] \ar@{-} [u]  \ar @{} [ul]|{\ss 1} \\ 
*={} \ar @{-} [u] &   *={} \ar @{-} [l] \ar@{-} [u]  \ar @{} [ul]|{\ss 1} &   *={} \ar @{-} [l] \ar@{-} [u]  \ar @{} [ul]|{\ss \bar{2}}&   *={} \ar @{-} [l] \ar@{-} [u]  \ar @{} [ul]|{\ss \bar{2}} \\
*={} \ar @{-} [u] &   *={} \ar @{-} [l] \ar@{-} [u]  \ar @{} [ul]|{\ss \bar{3}}\\
*={} \ar @{-} [u]  &   *={} \ar @{-} [l] \ar@{-} [u] \ar @{} [ul]|{\ss 3} }
\endxy
\quad\text{then}\quad Q_1=\xy<0cm,.5cm>\xymatrix@R=.35cm@C=.35cm{
*={} & *={} \ar @{.} [l] & *={} \ar @{.} [l] & *={} \ar @{-} [l] & *={} \ar @{-} [l] & *={} \ar @{-} [l]  \\
*={} \ar @{.} [u] &   *={} \ar @{-} [l] \ar@{.} [u]  &   *={} \ar @{.} [l] \ar@{-} [u]   &   *={} \ar @{-} [l] \ar@{-} [u]  \ar @{} [ul]|{\ss \bar{1}} &   *={} \ar @{-} [l] \ar@{-} [u]  \ar @{} [ul]|{\ss 1} &   *={} \ar @{-} [l] \ar@{-} [u]  \ar @{} [ul]|{\ss 1}   \\ 
*={} \ar @{-} [u] &   *={} \ar @{-} [l] \ar@{-} [u]  \ar @{} [ul]|{\ss 1} &   *={} \ar @{.} [l] \ar@{.} [u] &   *={} \ar @{.} [l] \ar@{.} [u]  \\
*={} \ar @{.} [u] &   *={} \ar @{.} [l] \ar@{.} [u]  \\
*={} \ar @{.} [u]  &   *={} \ar @{.} [l] \ar@{.} [u]  }
\endxy\quad \text{and}\quad 
\xy<0cm,.5cm>\xymatrix@R=.35cm@C=.35cm{
*={} & *={} \ar @{.} [l] & *={} \ar @{.} [l] & *={} \ar @{-} [l] & *={} \ar @{-} [l] & *={} \ar @{-} [l]  \\
*={} \ar @{.} [u] &   *={} \ar @{-} [l] \ar@{.} [u]  &   *={} \ar @{.} [l] \ar@{-} [u]   &   *={} \ar @{-} [l] \ar@{-} [u]  \ar @{} [ul]|{\ss \bar{1}} &   *={} \ar @{-} [l] \ar@{-} [u]  \ar @{} [ul]|{\ss \bar{1}} &   *={} \ar @{-} [l] \ar@{-} [u]  \ar @{} [ul]|{\ss 1}    \\ 
*={} \ar @{-} [u] &   *={} \ar @{-} [l] \ar@{-} [u]  \ar @{} [ul]|{\ss \bar{1}} &   *={} \ar @{.} [l] \ar@{.} [u] &   *={} \ar @{.} [l] \ar@{.} [u]  \\
*={} \ar @{.} [u] &   *={} \ar @{.} [l] \ar@{.} [u]  \\
*={} \ar @{.} [u]  &   *={} \ar @{.} [l] \ar@{.} [u]  }\endxy,\xy<0cm,.5cm>\xymatrix@R=.35cm@C=.35cm{
*={} & *={} \ar @{.} [l] & *={} \ar @{.} [l] & *={} \ar @{-} [l] & *={} \ar @{-} [l] & *={} \ar @{-} [l]  \\
*={} \ar @{.} [u] &   *={} \ar @{-} [l] \ar@{.} [u]  &   *={} \ar @{.} [l] \ar@{-} [u]   &   *={} \ar @{-} [l] \ar@{-} [u]  \ar @{} [ul]|{\ss 1} &   *={} \ar @{-} [l] \ar@{-} [u]  \ar @{} [ul]|{\ss 1} &   *={}  \ar @{-} [l] \ar@{-} [u]  \ar @{} [ul]|{\ss 1}  \\ 
*={} \ar @{-} [u] &   *={} \ar @{-} [l] \ar@{-} [u]  \ar @{} [ul]|{\ss \bar{1}} &   *={} \ar @{.} [l] \ar@{.} [u] &   *={} \ar @{.} [l] \ar@{.} [u]  \\
*={} \ar @{.} [u] &   *={} \ar @{.} [l] \ar@{.} [u]  \\
*={} \ar @{.} [u]  &   *={} \ar @{.} [l] \ar@{.} [u]  }\endxy \in \cS_Q^{(1)}. $$
(In fact, $|\cS_Q^{(1)}|=8$).

In light of (\ref{SignedEquation}), (b) is equivalent to
$$\sum_{Q'\in \cS_Q^{(r)}} (-1)^{\text{Number of $\bar{r}$'s in $Q'$}}=\left\{\begin{array}{ll} (-1)^{n(\sh(Q_r))}& \text{if $\sh(Q_r)$ has a domino tiling,}\\ 0 &\text{otherwise.}\end{array}\right.$$
Note that in row $j$, $Q'\in \cS_Q^{(r)}$ is of the form
$$\xy<0cm,.5cm>\xymatrix@R=.35cm@C=.35cm{ *={} & *={} & *={} &  *={} & *={} & *={} & *={} & *={} &  *={} & *={} & *={} & *={} & *={} & *={} & *={} & *={} & *={} \ar @{-} [l]  & *={} & *={} &  *={} \ar @{-} [lll] & *={} \ar @{-} [l] & *={} \ar @{-} [l]  & *={} \save []+<.13cm,.3cm> *{\overbrace{\hspace{1.6cm}}^{d_{j-1}}} \restore & *={} &  *={} \ar @{-} [lll] & *={} \ar @{-} [l]\\
*={} & *={} & *={} &  *={} & *={} & *={} & *={} \ar @{-} [l]  & *={} & *={} &  *={} \ar @{-} [lll] & *={} \ar @{-} [l] & *={} \ar @{-} [l]  & *={} \save []+<.13cm,.3cm> *{\overbrace{\hspace{1.6cm}}^{d_j}} \restore & *={} &  *={} \ar @{-} [lll] & *={} \ar @{-} [l] \ar @{-} [u] & *={} \ar @{-} [l]  \ar @{-} [u] \ar @{} [ul]|{\bar{r}}  & *={}  & *={} &  *={} \ar @{-} [lll] \ar @{-} [u] \ar @{} [ulll]|{\cdots} & *={} \ar @{-} [l] \ar @{-} [u] \ar @{} [ul]|{\bar{r}}  & *={} \ar @{-} [l]  \ar @{-} [u] \ar @{} [ul]|{?}  & *={} & *={} &  *={} \ar @{-} [lll] \ar @{-} [u] \ar @{} [ulll]|{\cdots} & *={} \ar @{-} [l] \ar @{-} [u] \ar @{} [ul]|{?}\\
 *={} & *={} & *={}  \save []+<.13cm,.3cm> *{\overbrace{\hspace{1.6cm}}^{d_{j+1}}} \restore &  *={} & *={} & *={} \ar @{-} [u] \ar @{-} [lllll] & *={} \ar @{-} [u] \ar @{-} [l] \ar @{} [ul]|{\bar{r}} & *={} & *={} &  *={} \ar @{-} [lll] \ar @{} [ulll]|{\cdots} \ar @{-} [u] & *={} \ar @{-} [l] \ar @{-} [u] \ar @{} [ul]|{\bar{r}} & *={} \ar @{-} [l]  \ar @{-} [u] \ar @{} [ul]|{?}  & *={} & *={} &  *={} \ar @{-} [lll] \ar @{-} [u] \ar @{} [ulll]|{\cdots} & *={} \ar @{-} [l] \ar @{-} [u] \ar @{} [ul]|{?}  & *={} \ar @{-} [l]  \ar @{-} [u] \ar @{} [ul]|{r}  & *={} & *={} &  *={} \ar @{-} [lll] \ar @{-} [u] \ar @{} [ulll]|{\cdots} & *={} \ar @{-} [l] \ar @{-} [u] \ar @{} [ul]|{r}\\
 *={} \ar @{-} [u] & *={} \ar @{-} [u] \ar @{-} [l] \ar @{} [ul]|{?} & *={} & *={} &  *={} \ar @{-} [lll] \ar @{} [ulll]|{\cdots} \ar @{-} [u] & *={} \ar @{-} [l] \ar @{-} [u] \ar @{} [ul]|{?} & *={} \ar @{-} [l]  \ar @{-} [u] \ar @{} [ul]|{r}  & *={} & *={} &  *={} \ar @{-} [lll] \ar @{-} [u] \ar @{} [ulll]|{\cdots} & *={} \ar @{-} [l] \ar @{-} [u] \ar @{} [ul]|{r} }\endxy \begin{array}{l} \longleftarrow{\ss\text{row $j-1$}}\vspace{-.1cm}\\ \longleftarrow{\ss\text{row $j$}}\vspace{-.1cm}\\ \longleftarrow{\ss\text{row $j+1$}}\end{array}$$
Thus, we have $d_j+1$ choices for the values in row $j$.  If the total number of choices is even, then exactly half of these choices give a positive sign and half give a negative sign.  So we have 
$$\sum_{Q'\in \cS_Q^{(r)}} (-1)^{\text{Number of $\bar{r}$'s in $Q'$}}=0,$$
unless $d_j$ is even for all rows $j$ (that is, $\sh(Q_r)$ has a tiling by dominoes).  In this case, the signs of all but one of the possible tableaux will cancel each other out, so the only tableau that we have to count has row $j$ of the form
$$\xy<0cm,.5cm>\xymatrix@R=.35cm@C=.35cm{ *={} & *={} & *={} &  *={} & *={} & *={} & *={} & *={} &  *={} & *={} & *={} & *={} & *={} & *={} & *={} & *={} & *={} \ar @{-} [l]  & *={} & *={} &  *={} \ar @{-} [lll] & *={} \ar @{-} [l] & *={} \ar @{-} [l]  & *={} \save []+<.13cm,.3cm> *{\overbrace{\hspace{1.6cm}}^{d_{j-1}}} \restore & *={} &  *={} \ar @{-} [lll] & *={} \ar @{-} [l]\\
 *={} & *={} & *={} &  *={} & *={} & *={} & *={} \ar @{-} [l]  & *={} & *={} &  *={} \ar @{-} [lll] & *={} \ar @{-} [l] & *={} \ar @{-} [l]  & *={} \save []+<.13cm,.3cm> *{\overbrace{\hspace{1.6cm}}^{d_j}} \restore & *={} &  *={} \ar @{-} [lll] & *={} \ar @{-} [l] \ar @{-} [u] & *={} \ar @{-} [l]  \ar @{-} [u] \ar @{} [ul]|{\bar{r}}  & *={}  & *={} &  *={} \ar @{-} [lll] \ar @{-} [u] \ar @{} [ulll]|{\cdots} & *={} \ar @{-} [l] \ar @{-} [u] \ar @{} [ul]|{\bar{r}}  & *={} \ar @{-} [l]  \ar @{-} [u] \ar @{} [ul]|{r}  & *={} & *={} &  *={} \ar @{-} [lll] \ar @{-} [u] \ar @{} [ulll]|{\cdots} & *={} \ar @{-} [l] \ar @{-} [u] \ar @{} [ul]|{r}\\
  *={} & *={} & *={}  \save []+<.13cm,.3cm> *{\overbrace{\hspace{1.6cm}}^{d_{j+1}}} \restore &  *={} & *={} & *={} \ar @{-} [u] \ar @{-} [lllll] & *={} \ar @{-} [u] \ar @{-} [l] \ar @{} [ul]|{\bar{r}} & *={} & *={} &  *={} \ar @{-} [lll] \ar @{} [ulll]|{\cdots} \ar @{-} [u] & *={} \ar @{-} [l] \ar @{-} [u] \ar @{} [ul]|{\bar{r}} & *={} \ar @{-} [l]  \ar @{-} [u] \ar @{} [ul]|{r}  & *={} & *={} &  *={} \ar @{-} [lll] \ar @{-} [u] \ar @{} [ulll]|{\cdots} & *={} \ar @{-} [l] \ar @{-} [u] \ar @{} [ul]|{r}  & *={} \ar @{-} [l]  \ar @{-} [u] \ar @{} [ul]|{r}  & *={} & *={} &  *={} \ar @{-} [lll] \ar @{-} [u] \ar @{} [ulll]|{\cdots} & *={} \ar @{-} [l] \ar @{-} [u] \ar @{} [ul]|{r}\\
 *={} \ar @{-} [u] & *={} \ar @{-} [u] \ar @{-} [l] \ar @{} [ul]|{r} & *={} & *={} &  *={} \ar @{-} [lll] \ar @{} [ulll]|{\cdots} \ar @{-} [u] & *={} \ar @{-} [l] \ar @{-} [u] \ar @{} [ul]|{r} & *={} \ar @{-} [l]  \ar @{-} [u] \ar @{} [ul]|{r}  & *={} & *={} &  *={} \ar @{-} [lll] \ar @{-} [u] \ar @{} [ulll]|{\cdots} & *={} \ar @{-} [l] \ar @{-} [u] \ar @{} [ul]|{r} }\endxy \begin{array}{l} \longleftarrow{\ss\text{row $j-1$}}\vspace{-.1cm}\\ \longleftarrow{\ss\text{row $j$}}\vspace{-.1cm}\\ \longleftarrow{\ss\text{row $j+1$}}\end{array}$$
 which can clearly be tiled by dominoes of the form $\xy<0cm,.25cm>\xymatrix@R=.25cm@C=.25cm{*={} & *={} \ar @{-} [l]\\ *={} \ar @{-} [u] & *={} \ar @{-} [u] \ar @{} [l]|{\sss r}\\ *={} \ar @{-} [u] &*={} \ar @{-} [u] \ar @{-} [l]}\endxy$ and $\xy<0cm,.12cm>\xymatrix@R=.25cm@C=.25cm{*={} & *={} \ar @{-} [l] & *={} \ar @{-} [l]\\ *={} \ar @{-} [u] & *={} \ar @{} [u]|{\sss r} \ar @{-} [l]  & *={} \ar @{-} [u] \ar @{-} [l]}\endxy$.  For this tableau, we have
 $$(-1)^{\text{Number of $\bar{r}$'s}}=(-1)^{n(\sh(Q_r))}.$$
 Thus, 
 $$s_{(m_0)}\prod_{r=1}^\ell \sum_{i=0}^{2m_r}(-1)^i s_{(i)}s_{(2m_r-i)}=\sum_{\lambda\in \cP_{2n-m_0}}(-1)^{n(\lambda)} |\cD_{(m_0,m_1,\ldots,m_\ell)}^\lambda| s_\lambda,$$
 as desired.
\end{proof}

Let $\blam\in \cP^\Theta$ and $\bgamma\in \cP^\Theta$ be such that $\hgt(\bgamma)\leq 1$ and $|\bgamma^{(\vphi)}|\leq \blam^{(\vphi)}_1$ for all $\vphi\in \Theta$.   A \emph{battery $\Theta$-tableau ${\mathbf Q}$ of shape} $\blam/\bgamma$ is a sequence of tableaux indexed by $\Theta$ such that
$${\mathbf Q}^{(\vphi)} =\left\{\begin{array}{ll} \text{a domino tableau of shape $\blam^{(\vphi)}/\bgamma^{(\vphi)}$} & \text{if $|\vphi|$ is odd,}\\  \text{a symplectic tableau of shape $\blam^{(\vphi)}/\bgamma^{(\vphi)}$} & \text{if $|\vphi|$ is even.}\end{array}\right.$$
The \emph{weight of} ${\mathbf Q}$ is  $\wt({\mathbf Q})=(\wt({\mathbf Q})_1,\wt({\mathbf Q})_2,\ldots)$, where
$$\wt({\mathbf Q})_i=\sum_{\vphi\in \Theta\atop |\vphi|\text{ odd}}  |\vphi| \wt({\mathbf Q}^{(\vphi)})_i +\sum_{\vphi\in \Theta\atop |\vphi|\text{ even}}\frac{|\vphi|}{2} \wt({\mathbf Q}^{(\vphi)})_i .$$
Let
\begin{equation}
\cB_{(k,\nu)}^{\blam}=\{{\mathbf Q} \text{ battery tableaux}\ \mid\ \sh({\mathbf Q})=\blam/\bgamma,\bgamma\in \cP^\Theta_k, \hgt(\bgamma)\leq 1,\wt({\mathbf Q})=\nu\}.
\end{equation}

\noindent\textbf{Example.} If 
$$\blam=\left(\xy<0cm,.3cm>\xymatrix@R=.3cm@C=.3cm{
	*={} & *={} \ar @{-} [l] & *={} \ar @{-} [l]  \\
	*={} \ar @{-} [u] &   *={} \ar @{-} [l] \ar@{-} [u] &   *={} \ar @{-} [l] \ar@{-} [u]   \\
*={} \ar @{-} [u] &   *={} \ar @{-} [l] \ar@{-} [u]  & *={} \ar @{-} [l] \ar@{-} [u]  }\endxy^{(\vphi_1)}, \  \xy<0cm,.3cm>\xymatrix@R=.3cm@C=.3cm{
	*={} & *={} \ar @{-} [l] & *={} \ar @{-} [l] \\
	*={} \ar @{-} [u] &   *={} \ar @{-} [l] \ar@{-} [u]  &  *={} \ar @{-} [l] \ar@{-} [u] }\endxy^{(\vphi_2)},\ 
\xy<0cm,.3cm>\xymatrix@R=.3cm@C=.3cm{
	*={} & *={} \ar @{-} [l] & *={} \ar @{-} [l] & *={} \ar @{-} [l]  \\
	*={} \ar @{-} [u] &   *={} \ar @{-} [l] \ar@{-} [u] &   *={} \ar @{-} [l] \ar@{-} [u]  &   *={} \ar @{-} [l] \ar@{-} [u]  \\ *={} \ar @{-} [u] &   *={} \ar @{-} [l] \ar@{-} [u] }\endxy^{(\vphi_3)}\ \right)\qquad\text{where $|\vphi_i|=i$},$$
then $\cB_{(2,(5,4))}^{\blam}$ contains
$$\left(\xy<0cm,.4cm>\xymatrix@R=.35cm@C=.35cm{
	*={} & *={} \ar @{-} [l] & *={} \ar @{-} [l]  \\
	*={} \ar @{-} [u] &   *={} \ar @{-} [l] \ar@{-} [u] \ar @{} [ul]|{0} &   *={} \ar @{-} [l] \ar@{-} [u] \ar @{} [ul]|{0}  \\
*={} \ar @{-} [u] &   *={} \ar @{-} [l] \ar@{} [u]|{1}  & *={} \ar @{-} [l] \ar@{-} [u]  }\endxy^{(\vphi_1)}, \xy<0cm,.4cm>\xymatrix@R=.35cm@C=.35cm{
	*={} & *={} \ar @{-} [l] & *={} \ar @{-} [l] \\
	*={} \ar @{-} [u] &   *={} \ar @{-} [l] \ar@{-} [u] \ar @{} [ul]|{1}  &  *={} \ar @{-} [l] \ar@{-} [u] \ar @{} [ul]|{2} }\endxy^{(\vphi_2)},
\xy<0cm,.4cm>\xymatrix@R=.35cm@C=.35cm{
	*={} & *={} \ar @{-} [l] & *={} \ar @{-} [l] & *={} \ar @{-} [l]  \\
	*={} \ar @{-} [u] &   *={} \ar @{} [l]|{1} \ar@{-} [u] &   *={} \ar @{-} [l] \ar@{} [u]|{2}  &   *={} \ar @{-} [l] \ar@{-} [u]  \\ *={} \ar @{-} [u] &   *={} \ar @{-} [l] \ar@{-} [u] }\endxy^{(\vphi_3)} \right),
\left(\xy<0cm,.4cm>\xymatrix@R=.35cm@C=.35cm{
	*={} & *={} \ar @{-} [l] & *={} \ar @{-} [l]  \\
	*={} \ar @{-} [u] &   *={} \ar @{-} [l] \ar@{-} [u] \ar @{} [ul]|{0} &   *={} \ar @{-} [l] \ar@{-} [u] \ar @{} [ul]|{0}  \\
*={} \ar @{-} [u] &   *={} \ar @{-} [l] \ar@{} [u]|{1}  & *={} \ar @{-} [l] \ar@{-} [u]  }\endxy^{(\vphi_1)}, \xy<0cm,.4cm>\xymatrix@R=.35cm@C=.35cm{
	*={} & *={} \ar @{-} [l] & *={} \ar @{-} [l] \\
	*={} \ar @{-} [u] &   *={} \ar @{-} [l] \ar@{-} [u] \ar @{} [ul]|{\bar{1}}  &  *={} \ar @{-} [l] \ar@{-} [u] \ar @{} [ul]|{\bar{2}} }\endxy^{(\vphi_2)},
\xy<0cm,.4cm>\xymatrix@R=.35cm@C=.35cm{
	*={} & *={} \ar @{-} [l] & *={} \ar @{-} [l] & *={} \ar @{-} [l]  \\
	*={} \ar @{-} [u] &   *={} \ar @{} [l]|{1} \ar@{-} [u] &   *={} \ar @{-} [l] \ar@{} [u]|{2}  &   *={} \ar @{-} [l] \ar@{-} [u]  \\ *={} \ar @{-} [u] &   *={} \ar @{-} [l] \ar@{-} [u] }\endxy^{(\vphi_3)}\right),
\left(\xy<0cm,.4cm>\xymatrix@R=.35cm@C=.35cm{
	*={} & *={} \ar @{-} [l] & *={} \ar @{-} [l]  \\
	*={} \ar @{-} [u] &   *={} \ar @{-} [l] \ar@{-} [u] \ar @{} [ul]|{0} &   *={} \ar @{-} [l] \ar@{-} [u] \ar @{} [ul]|{0}  \\
*={} \ar @{-} [u] &   *={} \ar @{-} [l] \ar@{} [u]|{1}  & *={} \ar @{-} [l] \ar@{-} [u]  }\endxy^{(\vphi_1)},  \xy<0cm,.4cm>\xymatrix@R=.35cm@C=.35cm{
	*={} & *={} \ar @{-} [l] & *={} \ar @{-} [l] \\
	*={} \ar @{-} [u] &   *={} \ar @{-} [l] \ar@{-} [u] \ar @{} [ul]|{\bar{1}}  &  *={} \ar @{-} [l] \ar@{-} [u] \ar @{} [ul]|{2} }\endxy^{(\vphi_2)},
\xy<0cm,.4cm>\xymatrix@R=.35cm@C=.35cm{
	*={} & *={} \ar @{-} [l] & *={} \ar @{-} [l] & *={} \ar @{-} [l]  \\
	*={} \ar @{-} [u] &   *={} \ar @{} [l]|{1} \ar@{-} [u] &   *={} \ar @{-} [l] \ar@{} [u]|{2}  &   *={} \ar @{-} [l] \ar@{-} [u]  \\ *={} \ar @{-} [u] &   *={} \ar @{-} [l] \ar@{-} [u] }\endxy^{(\vphi_3)} \right),
$$
$$
\left(\xy<0cm,.4cm>\xymatrix@R=.35cm@C=.35cm{
	*={} & *={} \ar @{-} [l] & *={} \ar @{-} [l]  \\
	*={} \ar @{-} [u] &   *={} \ar @{-} [l] \ar@{-} [u] \ar @{} [ul]|{0} &   *={} \ar @{-} [l] \ar@{-} [u] \ar @{} [ul]|{0}  \\
*={} \ar @{-} [u] &   *={} \ar @{-} [l] \ar@{} [u]|{1}  & *={} \ar @{-} [l] \ar@{-} [u]  }\endxy^{(\vphi_1)}, \xy<0cm,.4cm>\xymatrix@R=.35cm@C=.35cm{
	*={} & *={} \ar @{-} [l] & *={} \ar @{-} [l] \\
	*={} \ar @{-} [u] &   *={} \ar @{-} [l] \ar@{-} [u] \ar @{} [ul]|{1}  &  *={} \ar @{-} [l] \ar@{-} [u] \ar @{} [ul]|{\bar{2}} }\endxy^{(\vphi_2)},
\xy<0cm,.4cm>\xymatrix@R=.35cm@C=.35cm{
	*={} & *={} \ar @{-} [l] & *={} \ar @{-} [l] & *={} \ar @{-} [l]  \\
	*={} \ar @{-} [u] &   *={} \ar @{} [l]|{1} \ar@{-} [u] &   *={} \ar @{-} [l] \ar@{} [u]|{2}  &   *={} \ar @{-} [l] \ar@{-} [u]  \\ *={} \ar @{-} [u] &   *={} \ar @{-} [l] \ar@{-} [u] }\endxy^{(\vphi_3)} \right),
\left(\xy<0cm,.4cm>\xymatrix@R=.35cm@C=.35cm{
	*={} & *={} \ar @{-} [l] & *={} \ar @{-} [l]  \\
	*={} \ar @{-} [u] &   *={} \ar @{-} [l] \ar@{-} [u] \ar @{} [ul]|{0} &   *={} \ar @{-} [l] \ar@{-} [u] \ar @{} [ul]|{0}  \\
*={} \ar @{-} [u] &   *={} \ar @{-} [l] \ar@{} [u]|{2}  & *={} \ar @{-} [l] \ar@{-} [u]  }\endxy^{(\vphi_1)},  \xy<0cm,.4cm>\xymatrix@R=.35cm@C=.35cm{
	*={} & *={} \ar @{-} [l] & *={} \ar @{-} [l] \\
	*={} \ar @{-} [u] &   *={} \ar @{-} [l] \ar@{-} [u] \ar @{} [ul]|{\bar{1}}  &  *={} \ar @{-} [l] \ar@{-} [u] \ar @{} [ul]|{\bar{1}} }\endxy^{(\vphi_2)},
\xy<0cm,.4cm>\xymatrix@R=.35cm@C=.35cm{
	*={} & *={} \ar @{-} [l] & *={} \ar @{-} [l] & *={} \ar @{-} [l]  \\
	*={} \ar @{-} [u] &   *={} \ar @{} [l]|{1} \ar@{-} [u] &   *={} \ar @{-} [l] \ar@{} [u]|{2}  &   *={} \ar @{-} [l] \ar@{-} [u]  \\ *={} \ar @{-} [u] &   *={} \ar @{-} [l] \ar@{-} [u] }\endxy^{(\vphi_3)}\right),
\left(\xy<0cm,.4cm>\xymatrix@R=.35cm@C=.35cm{
	*={} & *={} \ar @{-} [l] & *={} \ar @{-} [l]  \\
	*={} \ar @{-} [u] &   *={} \ar @{-} [l] \ar@{-} [u] \ar @{} [ul]|{0} &   *={} \ar @{-} [l] \ar@{-} [u] \ar @{} [ul]|{0}  \\
*={} \ar @{-} [u] &   *={} \ar @{-} [l] \ar@{} [u]|{2}  & *={} \ar @{-} [l] \ar@{-} [u]  }\endxy^{(\vphi_1)}, \xy<0cm,.4cm>\xymatrix@R=.35cm@C=.35cm{
	*={} & *={} \ar @{-} [l] & *={} \ar @{-} [l] \\
	*={} \ar @{-} [u] &   *={} \ar @{-} [l] \ar@{-} [u] \ar @{} [ul]|{\bar{1}}  &  *={} \ar @{-} [l] \ar@{-} [u] \ar @{} [ul]|{1} }\endxy^{(\vphi_2)},
\xy<0cm,.4cm>\xymatrix@R=.35cm@C=.35cm{
	*={} & *={} \ar @{-} [l] & *={} \ar @{-} [l] & *={} \ar @{-} [l]  \\
	*={} \ar @{-} [u] &   *={} \ar @{} [l]|{1} \ar@{-} [u] &   *={} \ar @{-} [l] \ar@{} [u]|{2}  &   *={} \ar @{-} [l] \ar@{-} [u]  \\ *={} \ar @{-} [u] &   *={} \ar @{-} [l] \ar@{-} [u] }\endxy^{(\vphi_3)} \right),
$$
$$	
\left(\xy<0cm,.4cm>\xymatrix@R=.35cm@C=.35cm{
*={} & *={} \ar @{-} [l] & *={} \ar @{-} [l]  \\
	*={} \ar @{-} [u] &   *={} \ar @{-} [l] \ar@{-} [u] \ar @{} [ul]|{0} &   *={} \ar @{-} [l] \ar@{-} [u] \ar @{} [ul]|{0}  \\
*={} \ar @{-} [u] &   *={} \ar @{-} [l] \ar@{} [u]|{2}  & *={} \ar @{-} [l] \ar@{-} [u]  }\endxy^{(\vphi_1)},   \xy<0cm,.4cm>\xymatrix@R=.35cm@C=.35cm{
	*={} & *={} \ar @{-} [l] & *={} \ar @{-} [l] \\
	*={} \ar @{-} [u] &   *={} \ar @{-} [l] \ar@{-} [u] \ar @{} [ul]|{1}  &  *={} \ar @{-} [l] \ar@{-} [u] \ar @{} [ul]|{1} }\endxy^{(\vphi_2)},
\xy<0cm,.4cm>\xymatrix@R=.35cm@C=.35cm{
	*={} & *={} \ar @{-} [l] & *={} \ar @{-} [l] & *={} \ar @{-} [l]  \\
	*={} \ar @{-} [u] &   *={} \ar @{} [l]|{1} \ar@{-} [u] &   *={} \ar @{-} [l] \ar@{} [u]|{2}  &   *={} \ar @{-} [l] \ar@{-} [u]  \\ *={} \ar @{-} [u] &   *={} \ar @{-} [l] \ar@{-} [u] }\endxy^{(\vphi_3)} \right).
$$

\vspace{.25cm}

\noindent\textbf{Some intuition.}  If $\blam \in \cP^\Theta$, we can think of the boxes in $\blam^{(\vphi)}$ as being $|\vphi|$ deep, so in the above example,  
$$\blam=\left(\xy<0cm,.5cm>\xymatrix@R=.15cm@C=.15cm{*={} & *={} \ar @{-} [rrrrrr] \ar @{.} [dddddd]  & *={}  & *={} & *={} \ar @{.} [dddddd] & *={} & *={} & *={} \ar @{-} [dddddd] \\
*={}\ar @{-} [ur] \ar @{-} [dddddd] \ar @{-} [rrrrrr]  & *={} & *={} & *={} \ar @{-} [dddddd] \ar @{-} [ur] & *={} & *={} & *={} \ar @{-} [ur] \ar @{-} [dddddd]\\
*={}\\
*={} & *={} \ar @{.} [rrrrrr]  & *={} & *={} & *={} & *={} & *={} & *={} \\
*={} \ar @{.} [ur] \ar @{-} [rrrrrr] & *={} & *={} & *={} \ar @{.} [ur] & *={} & *={} & *={} \ar @{-} [ur]\\
*={}\\
*={} & *={} \ar @{.} [rrrrrr]  & *={} & *={} & *={} & *={} & *={} & *={} \\
*={} \ar @{.} [ur] \ar @{-} [rrrrrr] & *={} & *={} & *={} \ar @{.} [ur] & *={} & *={} & *={} \ar @{-} [ur]   }\endxy^{(\vphi_1)}, \  
\xy<0cm,.3cm>\xymatrix@R=.15cm@C=.15cm{*={} & *={} & *={} \ar @{-} [rrrrrr] \ar @{.} [ddd] & *={}  & *={} & *={} \ar @{.} [ddd] & *={} & *={} & *={} \ar @{-} [ddd] \\
*={} & *={} \ar @{-} [rrrrrr] \ar @{.} [ddd] & *={}  & *={} & *={} \ar @{.} [ddd] & *={} & *={} & *={} \ar @{-} [ddd] \\
*={}\ar @{-} [uurr] \ar @{-} [ddd] \ar @{-} [rrrrrr]  & *={} & *={} & *={} \ar @{-} [ddd] \ar @{-} [uurr] & *={} & *={} & *={} \ar @{-} [uurr] \ar @{-} [ddd]\\
*={} & *={} & *={} \ar @{.} [rrrrrr]  & *={} & *={} & *={} & *={} & *={} & *={}\\
*={} & *={} \ar @{.} [rrrrrr]  & *={} & *={} & *={} & *={} & *={} & *={} \\
*={} \ar @{.} [uurr] \ar @{-} [rrrrrr] & *={} & *={} & *={} \ar @{.} [uurr] & *={} & *={} & *={} \ar @{-} [uurr]  }\endxy^{(\vphi_2)},\ 
\xy<0cm,.6cm>\xymatrix@R=.15cm@C=.15cm{*={} & *={} & *={} & *={} \ar @{-} [rrrrrrrrr] \ar @{.} [dddddd] & *={}  & *={} & *={} \ar @{.} [dddddd] & *={} & *={} & *={} \ar @{.} [ddd] & *={} & *={} & *={} \ar @{-} [ddd] \\
*={} & *={} & *={} \ar @{-} [rrrrrrrrr] \ar @{.} [dddddd] & *={}  & *={} & *={} \ar @{.} [dddddd] & *={} & *={} & *={} \ar @{.} [ddd] & *={} & *={} & *={} \ar @{-} [ddd] \\
*={} & *={} \ar @{-} [rrrrrrrrr] \ar @{.} [dddddd] & *={}  & *={} & *={} \ar @{.} [dddddd] & *={} & *={} & *={} \ar @{.} [ddd] & *={} & *={} & *={} \ar @{-} [ddd] \\
*={}\ar @{-} [uuurrr] \ar @{-} [dddddd] \ar @{-} [rrrrrrrrr]  & *={} & *={} & *={} \ar @{-} [dddddd] \ar @{-} [uuurrr] \ar @{.} [rrrrrrrrr] & *={} & *={} & *={} \ar @{-} [uuurrr] \ar @{-} [ddd] & *={} & *={} & *={} \ar @{-} [ddd] \ar @{-} [uuurrr] & *={} & *={} & *={}\\
*={} & *={} & *={} \ar @{.} [rrrrrrrrr]  & *={} & *={} & *={} & *={} & *={} & *={} & *={} & *={} & *={}\\
*={} & *={} \ar @{.} [rrrrrrrrr]  & *={} & *={} & *={} & *={} & *={} & *={} & *={} & *={} & *={} \\
*={} \ar @{.} [uuurrr] \ar @{-} [rrrrrrrrr] & *={} & *={} & *={} \ar @{.} [uuurrr] & *={} & *={} & *={} \ar @{.} [uuurrr] & *={} & *={} & *={} \ar @{-} [uuurrr] \\
*={} & *={} & *={} \ar @{.} [rrr]  & *={} & *={} & *={}  \ar @{-} [u] \\
*={} & *={} \ar @{.} [rrr]  & *={} & *={} & *={}   \ar @{-} [uu] \\
*={} \ar @{.} [uuurrr] \ar @{-} [rrr] & *={} & *={} & *={} \ar @{-} [uuurrr]  \\
 }\endxy^{(\vphi_3)}\right).$$
A battery $\Theta$-tableau is a way of stuffing the slots by numbered ``batteries" where front and back are distinguished by $i$ and $\bar{i}$, but the sides look generically like $i$, so
$$
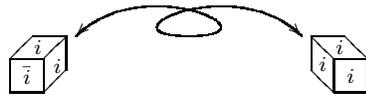
\xymatrix@R=.15cm@C=.15cm{*={} & *={} & *={} \ar @{-} [rrr] & *={}  & *={} & *={} \ar @{-} [ddd]   \\
*={} & *={} & *={}  & *={} & *={} & *={}    \\
*={}\ar @{-} [uurr] \ar @{-} [ddd]\ar @{-} [rrr] \ar @{} [dddrrr]|{\bar{i}} \ar @{} [uurrrrr]|{i} & *={} & *={} & *={} \ar @{-} [ddd] \ar @{-} [uurr]   & *={} & *={} \\
*={} & *={} & *={}  & *={} & *={} & *={}  \\
*={} & *={}  & *={} & *={} & *={} & *={}   \\
*={}  \ar @{-} [rrr] & *={} & *={} & *={} \ar @{-}[uurr] \ar @{} [uuuurr]|(.6){i} & *={} & *={}}\ 
\xy<0cm,.2cm>\xymatrix@R=.1cm@C=1.5cm{*={}  & *={} \ar @(r,ur) [l]\ar @(l,ul) [r] & *={}}\endxy\ 
\xymatrix@R=.15cm@C=.15cm{*={} \ar @{-} [ddd] \ar@{-} [rrr] \ar @{-} [ddrr] \ar @{} [dddddrr]|{i} \ar @{} [ddrrrrr]|{i}  & *={} & *={} & *={}  \ar @{-} [ddrr]   \\
*={} & *={} & *={}  & *={} & *={} & *={} & *={} & *={}  \\
*={} & *={}  & *={} \ar @{-} [rrr] & *={} & *={} & *={} & *={} & *={} \\
*={}\ar @{-} [ddrr]  & *={} & *={} & *={}  & *={} & *={} & *={} \\
*={} & *={} & *={}  & *={} & *={} & *={} & *={} & *={} & *={}\\
*={}  & *={} & *={} \ar @{-} [rrr] \ar @{-} [uuu] \ar @{} [uuurrr]|{i} & *={} & *={} & *={} \ar @{-} [uuu] & *={}}  
$$
Then a battery $\Theta$-tableau might look like:
$$\left(\xy<0cm,.5cm>\xymatrix@R=.15cm@C=.15cm{*={} & *={} \ar @{-} [rrrrrr]  & *={}  & *={} & *={} & *={} & *={} & *={} \ar @{-} [dddddd] \\
*={}\ar @{-} [ur] \ar @{-} [dddddd] \ar @{-} [rrrrrr] \ar @{} [dddrrr]|{0} & *={} & *={} & *={} \ar @{-} [ddd] \ar @{-} [ur]  \ar @{} [dddrrr]|{0} & *={} & *={} & *={} \ar @{-} [ur] \ar @{-} [dddddd]\\
*={}\\
*={} & *={}  & *={} & *={} & *={} & *={} & *={} & *={} \\
*={} \ar @{-} [rrrrrr] & *={} & *={} & *={} & *={} & *={} & *={} \ar @{-} [ur]\\
*={}\\
*={} & *={}  & *={} & *={} & *={} & *={} & *={} & *={} \\
*={} \ar @{-} [rrrrrr] \ar @{} [uuurrrrrr]|{1} & *={} & *={} & *={} & *={} & *={} & *={} \ar @{-} [ur]   }\endxy^{(\vphi_1)}, \  
\xy<0cm,.3cm>\xymatrix@R=.15cm@C=.15cm{*={} & *={} & *={} \ar @{-} [rrrrrr] & *={}  & *={} & *={} & *={} & *={} & *={} \ar @{-} [ddd] \\
*={} & *={} & *={}  & *={} & *={} & *={} & *={} & *={}  \\
*={}\ar @{-} [uurr] \ar @{-} [ddd]\ar @{-} [rrrrrr] \ar @{} [dddrrr]|{1} & *={} & *={} & *={} \ar @{-} [ddd] \ar @{-} [uurr]  \ar @{} [dddrrr]|{\bar{2}} & *={} & *={} & *={} \ar @{-} [uurr] \ar @{-} [ddd]\\
*={} & *={} & *={}  & *={} & *={} & *={} & *={} & *={} & *={}\\
*={} & *={}  & *={} & *={} & *={} & *={} & *={} & *={} \\
*={}  \ar @{-} [rrrrrr] & *={} & *={} & *={} & *={} & *={} & *={} \ar @{-} [uurr]  }\endxy^{(\vphi_2)},\ 
\xy<0cm,.6cm>\xymatrix@R=.15cm@C=.15cm{*={} & *={} & *={} & *={} \ar @{-} [rrrrrrrrr] \ar @{} [ddr]|{\sss 1} & *={}  & *={} & *={} \ar @{} [ddr]|{\sss 2} & *={} & *={} & *={} \ar @{} [ddr]|{\sss 2}  & *={} & *={} & *={} \ar @{-} [ddd] \\
*={} & *={} & *={}  & *={}  & *={} & *={} & *={} & *={} & *={} & *={} & *={} & *={} \\
*={} & *={} \ar @{-} [rrrrrrrrr] & *={}  & *={} & *={} & *={} & *={} & *={}  \ar @{-} [uurr] & *={} & *={} & *={} \ar @{-} [ddd] \\
*={}\ar @{-} [uuurrr] \ar @{-} [dddddd] \ar @{-} [rrrrrrrrr]  & *={} & *={} & *={} \ar @{-} [dddddd] \ar @{-} [uuurrr]  & *={} & *={} & *={} \ar @{} [ddd]|{2} & *={} & *={} & *={} \ar @{-} [ddd] \ar @{-} [uuurrr] & *={} & *={} & *={}\\
*={} & *={} & *={}  & *={} & *={} & *={} & *={} & *={} & *={} & *={} & *={} & *={}\\
*={} & *={}  & *={} & *={} & *={} & *={} & *={} & *={} & *={} & *={} & *={} \\
*={} \ar @{} [rrr]|{1} & *={} & *={} & *={}  \ar @{-} [rrrrrr] & *={} & *={} & *={} & *={} & *={} & *={} \ar @{-} [uuurrr] \\
*={} & *={} & *={}   & *={} & *={} & *={} \save[]+<.25cm,-.2cm>*{\ss 1} \ar []-<.05cm,-.09cm> \restore\\
*={} & *={}   & *={} & *={} & *={} \ar @{-} [uu]  \\
*={}  \ar @{-} [rrr] & *={} & *={} & *={} \ar @{-} [uuurrr]  \\
 }\endxy^{(\vphi_3)}\right),$$
 so the weight of the tableau counts the number of batteries of a given type get used, regardless of the cardinality of $\vphi$.

\subsection{Inducing from $G_n$ to $U_{2n}$}\label{SectionBaseCase}

Note that any maximal torus $\tilde{T}_\nu$ of $G_n\subseteq U_{2n}$ becomes the maximal torus $T_{2\nu}$ of $U_{2n}$, which gives rise to the map
$$\begin{array}{rccc} i:& \left\{\begin{array}{c}\text{Pairs $(\tilde{T}_\nu,\tilde{\theta}_\nu)$ with $\tilde{T}_\nu$ a}\\ \text{ maximal torus of $G_n$,}\\ \tilde\theta_\nu\in \Hom(\tilde{T}_\nu,\CC^\times)\end{array}\right\} &\longrightarrow &\left\{\begin{array}{c}\text{Pairs $(T_{2\nu},\theta_\nu)$ with $T_{2\nu}$ a}\\ \text{maximal torus of $U_{2n}$,}\\ \theta_\nu\in \Hom(T_{2\nu},\CC^\times)\end{array}\right\} \\
& (\tilde{T}_\nu,\tilde\theta_\nu) & \mapsto & (T_{2\nu},\tilde\theta_\nu).\end{array}$$
To translate the combinatorics between $G_n$ and $U_{2n}$, we define the map 
$$\begin{array}{rcc}\iota:\P_n^{\tilde\Theta} & \longrightarrow & \P_{2n}^\Theta\\ \tilde\blam & \mapsto & \iota\tilde\blam\end{array}\qquad \text{where for $\vphi\in \Theta$,} \qquad \iota\tilde\blam^{(\vphi)}=\left\{\begin{array}{ll} 2\tilde\blam^{(\tilde{\vphi})} & \text{if $\vphi=\tilde\vphi$}, \\ \tilde\blam^{(\tilde\vphi)}\cup\tilde{\blam}^{(F(\tilde\vphi))} &\text{if $\vphi=\tilde\vphi\cup F\tilde\vphi$,}\end{array}\right.$$
which has the property that $\tau_{\Theta} \circ i = \iota \circ \tau_{\tilde\Theta}$ (see (\ref{TorusToTheta})).  The map $\iota$ is neither surjective nor injective.  We note that $F\tilde\vphi=\tilde\vphi$ implies that $|\tilde\vphi|$ is odd, and if $F\tilde\vphi\neq \tilde\vphi$, then $\vphi=\tilde\vphi\cup F\tilde\vphi$ implies $|\vphi|$ is even (see \cite{ennolaconj}).  Thus, the image of $\iota$ is the set of \emph{even} $\Theta$-partitions,
$$\mathrm{Image}(\iota)=\{\blam\in \P_n^\Theta\ \mid\ |\vphi|\blam^{(\vphi)} \text{ is even for } \vphi\in \Theta\}.$$

 \begin{theorem} \label{DGGBaseCase}
$$R_{G_n}^{U_{2n}}(\tilde\Gamma_{(n)})=\sum_{\blam\in \cP_{2n}^\Theta\atop \hgt(\blam)\leq 2} |\cB_{(0,(n))}^{\blam}| \chi^{\blam}.$$
\end{theorem}

\begin{proof}   
Note that by Theorem \ref{GGDecomposition}, (\ref{SCHURtoPOWER}), and the characteristic map for $G_n$,
\begin{equation*}
\tilde\Gamma_{(n)}=\sum_{\tilde\blam\in \cP_n^{\tilde\Theta}\atop \hgt(\tilde\blam)= 1} \chi^{\tilde\blam}=\sum_{\tilde\blam\in \cP_n^{\tilde\Theta}\atop \hgt(\blam)= 1}\sum_{\tilde\bnu\in \cP^{\tilde\blam}_s} \frac{(-1)^{n-\ell(\tilde\bnu)}}{z_{\tilde\bnu}}R_{\tilde\bnu}^{G_n}.
\end{equation*}
By transitivity of induction, and the fact that $\tau_{\Theta} \circ i = \iota \circ \tau_{\tilde\Theta}$, we have $R_{G_n}^{U_{2n}}(R_{\tilde\bnu}^{G_n})=R_{\iota\tilde\bnu}^{U_{2n}}$, and so
$$R_{G_n}^{U_{2n}}(\tilde\Gamma_{(n)})=\sum_{\tilde\blam\in \cP_n^{\tilde\Theta}\atop \hgt(\tilde\blam)= 1}\sum_{\tilde\bnu\in \cP^{\tilde\blam}_s}  \frac{(-1)^{n - \ell(\tilde\bnu)}}{z_{\tilde\bnu}}R_{\iota\tilde\bnu}^{U_{2n}}.$$
We now change the second sum to a sum over $\bnu = \iota\tilde\bnu \in \cP_s^{\iota\tilde\blam}$, and we obtain
\begin{align*}
R_{G_n}^{U_{2n}}(\tilde\Gamma_{(n)}) &=\sum_{\tilde\blam\in \cP_n^{\tilde\Theta}\atop \hgt(\tilde\blam)= 1} \sum_{\bnu\in \cP^{\iota\tilde\blam}_s} \biggl(\sum_{\tilde\bnu\in \P_s^{\tilde\blam}\atop \iota\tilde\bnu= \bnu} \frac{1}{z_{\tilde\bnu}}\biggr) (-1)^{n-\ell(\bnu)} R_{\bnu}^{U_{2n}}\\
&=\sum_{\bnu\in \cP^\Theta_{2n}\atop \bnu\text{ even}} \biggl(\sum_{\tilde\bnu\in \P_n^{\tilde\Theta}\atop \iota\tilde\bnu= \bnu} \frac{1}{z_{\tilde\bnu}}\biggr) (-1)^{n-\ell(\bnu)} R_{\bnu}^{U_{2n}}.
\end{align*}
Recall that $F\tilde\vphi=\tilde\vphi$ implies that $|\tilde\vphi|$ is odd, and $F\tilde\vphi\neq \tilde\vphi$ implies that $\vphi=\tilde\vphi\cup F\tilde\vphi$ where $|\vphi|$ is even.  Apply the characteristic map, factor, and then reindex to obtain
\begin{align*}
\ch(R_{G_n}^{U_{2n}}(\tilde\Gamma_{(n)}))&=(-1)^n\sum_{\bnu\in \cP_{2n}^\Theta\atop \bnu\text{ even}} \biggl(\sum_{\tilde\bnu\in \P^{\tilde\Theta}\atop \iota\tilde\bnu= \bnu} \frac{1}{z_{\tilde\bnu}}\biggr) p_{\bnu}\\
&=(-1)^n\sum_{\bnu\in \cP_{2n}^\Theta\atop \bnu\text{ even}}\prod_{\vphi\in \Theta\atop |\vphi|\text{ odd}} \frac{1}{z_{\bnu^{(\vphi)}/2}} p_{\bnu^{(\vphi)}}(Y^{(\vphi)}) \prod_{\vphi\in \Theta\atop |\vphi|\text{ even}} \biggl(\sum_{\nu,\mu\in \cP\atop \eta\cup\mu= \bnu^{(\vphi)}} \frac{1}{z_{\eta}z_\mu}\biggr) p_{\bnu^{(\vphi)}}(Y^{(\vphi)})\\
&=(-1)^n\sum_{{\bgamma\in \cP_{2n}^\Theta\atop \hgt(\bgamma)= 1}\atop \bgamma\text{ even}}\prod_{\vphi\in \Theta\atop |\vphi|\text{ odd}} \biggl(\sum_{{\nu \in \cP\atop |\nu|=|\bgamma^{(\vphi)}|}\atop \nu\text{ even}} \frac{1}{z_{\nu/2}} p_{\nu}(Y^{(\vphi)})\biggr) \hspace{-.25cm} \prod_{\vphi\in \Theta\atop |\vphi|\text{ even}}\hspace{-.25cm} \biggl(\sum_{\eta,\mu\in \cP\atop |\eta|+|\mu|=|\bgamma^{(\vphi)}|} \frac{1}{z_{\eta}z_\mu}p_{\eta\cup\mu}(Y^{(\vphi)})\biggr). 
\end{align*}
Note that by (\ref{SchurToPower}), 
$$\sum_{\eta,\mu\in \cP\atop |\eta|+|\mu|=|\gamma|} \frac{1}{z_{\eta}z_\mu}p_{\eta\cup\mu}=\sum_{i=0}^{|\gamma|}\biggl(\sum_{|\eta|=i}z_{\eta}^{-1}p_\eta\biggr)\biggl(\sum_{|\mu|=|\gamma|-i}z_{\mu}^{-1}p_\mu\biggr)=\sum_{i=0}^{|\gamma|}s_{(i)}s_{(|\gamma|-i)}.$$
A computation similar to \cite[I.2.14]{Mac} shows that
$$\sum_{|\nu|=|\gamma|\atop \nu\text{ even}} \frac{1}{z_{\nu/2}} p_{\nu}=\sum_{i=0}^{|\gamma|} (-1)^i s_{(i)}s_{(|\gamma|-i)}.$$
Thus,
\begin{equation} \label{ProductsOfSchurs}
\ch(R_{G_n}^{U_{2n}}(\tilde\Gamma_{(n)}))=(-1)^n\sum_{{\bgamma\in \cP_n^\Theta\atop \hgt(\bgamma)= 1}\atop \bgamma\text{ even}}\prod_{\vphi\in \Theta} \sum_{i=0}^{|\bgamma^{(\vphi)}|}(-1)^{|\vphi|i}s_{(i)}(Y^{(\vphi)}) s_{(|\bgamma^{(\vphi)}|-i)}(Y^{(\vphi)}).\end{equation}
Lemma \ref{PieriLemma} (a) and (b), respectively, imply that
\begin{align*}
\sum_{i=0}^{k} s_{(i)}s_{(k-i)} &= \sum_{\lambda\in \cP_k}|\cT_{(0,k)}^\lambda| s_\lambda, \text{ and }\\
\sum_{i=0}^{k}(-1)^{i}s_{(i)}s_{(k-i)} &= \sum_{\lambda\in \cP_{k}}(-1)^{n(\lambda)} |\cD_{(0,k/2)}^\lambda| s_\lambda.
\end{align*}
Since $|\cD_{(0,k/2)}^\lambda|=|\cT_{(0,k)}^\lambda|=0$ unless $\hgt(\lambda)\leq 2$, 
\begin{align*}
&\ch(R_{G_n}^{U_{2n}}(\tilde\Gamma_{(n)}))\\
&=(-1)^n\hspace{-.25cm}\sum_{{\bgamma\in \cP_{2n}^\Theta\atop \hgt(\bgamma)= 1}\atop \bgamma\text{ even}}\prod_{\vphi\in \Theta\atop |\vphi|\text{ odd}} \sum_{|\blam^{(\vphi)}|={|\bgamma^{(\vphi)}|}} \hspace{-.35cm}(-1)^{n(\blam^{(\vphi)})} \bigl|\cD_{(0,|\bgamma^{(\vphi)}|/2)}^{\blam^{(\vphi)}}\bigr| s_{\blam^{(\vphi)}}(Y^{(\vphi)})\\
& \hspace*{8cm}\cdot\prod_{\vphi\in \Theta\atop |\vphi|\text{ even}} \sum_{|\blam^{(\vphi)}|=|\bgamma^{(\vphi)}|}\hspace{-.35cm}\bigl|\cT_{(0,|\bgamma^{(\vphi)}|)}^{\blam^{(\vphi)}}\bigr| s_{\blam^{(\vphi)}}(Y^{(\vphi)})\\
&=\sum_{\blam\in \cP_{2n}^\Theta\atop \hgt(\blam)\leq 2} (-1)^{n+n(\blam)} |\cB_{(0,(n))}^{\blam}| s_{\blam}.
\end{align*}
Apply $\ch^{-1}$ to get the result.
\end{proof}

\begin{corollary} \label{DGGBaseCaseCorollary} For $n\in \ZZ_{\geq 1}$, 
$$\ch(R_{G_n}^{U_{2n}}(\tilde\Gamma_{(n)}))=(-1)^n\sum_{{\bnu\in \cP_{2n}^\Theta\atop \hgt(\bnu)= 1}\atop \bnu\text{ even}} \prod_{\vphi\in \Theta} \sum_{i=0}^{|\bnu^{(\vphi)}|} (-1)^{i|\vphi|}s_{(i)}(Y^{(\varphi)})s_{(|\bnu^{(\vphi)}|-i)}(Y^{(\varphi)}).$$
\end{corollary}

\begin{proof} This is (\ref{ProductsOfSchurs}) in the proof of Theorem \ref{DGGBaseCase}.\end{proof}

Using similar techniques, we can prove a result for arbitrary irreducible characters of $G_n$.  For $\blam \in \cP^{\tilde\Theta}$ and $\bgamma\in \cP^{\iota\tilde\blam}_s$, let
$$c_{\tilde\blam}^{\bgamma}=\prod_{\vphi\in \Theta} c_{\tilde\blam}^{\bgamma}(\vphi),\qquad\text{where}\qquad c_{\tilde\blam}^{\bgamma}(\vphi)=\left\{\begin{array}{ll} c_{\blam^{(\tilde\vphi)}}^{\bgamma^{(\vphi)}} &\text{if $\vphi=\tilde\vphi\in \tilde\Theta$,}\\ \\
 c_{\tilde\blam^{(\tilde\vphi)}\tilde\blam^{(F\tilde\vphi)}}^{\bgamma^{(\vphi)}} & \text{if $\vphi=\tilde\vphi\cup F\tilde\vphi$ and $F\tilde\vphi\neq \tilde\vphi\in \tilde\Theta$,}\end{array}\right.$$ 
where $c_{\nu\mu}^\lambda$ is as in (\ref{LRCoefficients}), and $c_\lambda^\gamma$ is as in (\ref{PlethysmCoefficients}).

\begin{theorem} \label{InducedIrreducibles} Let $\tilde\blam\in \cP_n^{\tilde{\Theta}}$.  Then
$$R_{G_n}^{U_{2n}}(\chi^{\tilde\blam})=\sum_{\bgamma\in \cP_s^{\iota\tilde\blam}}(-1)^{n(\bgamma)} c_{\tilde\blam}^{\bgamma} \chi^{\bgamma}.$$
\end{theorem}

\begin{proof}
By (\ref{SCHURtoPOWER}) and the characteristic map for $G_n$, 
\begin{equation*}
\chi^{\tilde\blam}=\sum_{\tilde\bnu\in \cP^{\tilde\blam}_s} \Big(\prod_{\tilde\vphi \in \tilde\Theta} \frac{\omega^{\tilde\blam^{(\tilde\vphi)}}(\tilde\bnu^{(\tilde\vphi)})}{z_{\tilde\bnu^{(\tilde\vphi)}}}\Big)(-1)^{n-\ell(\tilde\bnu)}R_{\tilde\bnu}^{G_n}.
\end{equation*}
By transitivity of induction, and the fact that $\tau_{\Theta} \circ i = \iota \circ \tau_{\tilde\Theta}$, we have $R_{G_n}^{U_{2n}}(R_{\tilde\bnu}^{G_n})=R_{\iota\tilde\bnu}^{U_{2n}}$, and so
$$R_{G_n}^{U_{2n}}(\chi^{\tilde\blam})=\sum_{\tilde\bnu\in \cP^{\tilde\blam}_s} \Big(\prod_{\tilde\vphi \in \tilde\Theta} \frac{\omega^{\tilde\blam^{(\tilde\vphi)}}(\tilde\bnu^{(\tilde\vphi)})}{z_{\tilde\bnu^{(\tilde\vphi)}}}\Big)(-1)^{n-\ell(\tilde\bnu)}R_{\iota\tilde\bnu}^{U_{2n}} .$$
We now change the sum to a sum over $\bnu = \iota\tilde\bnu \in \cP_s^{\iota\tilde\blam}$, and using the image of the map $\iota$, we obtain
$$R_{G_n}^{U_{2n}}(\chi^{\tilde\blam}) = \sum_{\bnu\in \P_s^{\iota\tilde\blam}\atop \bnu\text{ even}} \biggl(\sum_{\tilde\bnu\in \P_s^{\tilde\blam}\atop \iota\tilde\bnu= \bnu}\Big(\prod_{\vphi \in \Theta} \frac{\omega^{\tilde\blam^{(\vphi)}}(\tilde\bnu^{(\vphi)})}{z_{\tilde\bnu^{(\vphi)}}}\Big)\biggr) (-1)^{n-\ell(\bnu)} R_{\bnu}^{U_{2n}}.$$
Apply the characteristic map, and rewrite the inner sum and product, to get
\begin{align*}
\ch\bigl(R_{G_n}^{U_{2n}}(\chi^{\tilde\blam})\bigr) &=(-1)^n \sum_{\bnu\in \P_s^{\iota\tilde\blam}\atop \bnu\text{ even}} \biggl(\sum_{\tilde\bnu\in \P_s^{\tilde\blam}\atop \iota\tilde\bnu= \bnu} \Big(\prod_{\tilde\vphi \in \tilde\Theta} \frac{\omega^{\tilde\blam^{(\tilde\vphi)}}(\tilde\bnu^{(\tilde\vphi)})}{z_{\tilde\bnu^{(\tilde\vphi)}}}\Big) \biggr)p_{\bnu}\\
&= (-1)^n\sum_{\bnu\in \P_s^{\iota\tilde\blam}\atop \bnu\text{ even}}\prod_{\tilde\vphi\in \tilde\Theta\atop F\tilde\vphi=\tilde\vphi}\frac{\omega^{\tilde\blam^{(\tilde\vphi)}}(\bnu^{(\tilde\vphi)}/2)}{z_{\bnu^{(\tilde\vphi)}/2}} \prod_{\tilde\vphi\in \tilde\Theta\atop F\tilde\vphi\neq \tilde\vphi} \biggl(\sum_{|\gamma|=|\tilde\blam^{(\tilde\vphi)}|\atop {|\mu|=|\tilde\blam^{(F\tilde\vphi)}| \atop \gamma\cup \mu = \bnu^{(\tilde\vphi\cup F\tilde\vphi)}} } \frac{\omega^{\tilde\blam^{(\tilde\vphi)}}(\gamma) \omega^{\tilde\blam^{(F\tilde\vphi)}}(\mu) }{z_{\gamma}z_{\mu}}\biggr) p_{\bnu}.\end{align*}
Recall that $F\tilde\vphi=\tilde\vphi$ implies that $|\tilde\vphi|$ is odd, and $F\tilde\vphi\neq \tilde\vphi$ implies that $\vphi=\tilde\vphi\cup F\tilde\vphi$ where $|\vphi|$ is even.  Thus, for every $\vphi\in \Theta$ such that $\bnu^{(\vphi)}\neq \emptyset$, if $|\vphi|$ is odd then $\vphi=\tilde\vphi$ for some $\tilde\vphi\in \tilde\Theta$, and if $|\vphi|$ is even then $\vphi=\tilde\vphi\cup F\tilde\vphi$ for some $\tilde\vphi\in \tilde\Theta$.  Factor our expression accordingly as
\begin{align*}
&\ch\bigl(R_{G_n}^{U_{2n}}(\chi^{\tilde\blam})\bigr) =(-1)^n\sum_{\bnu\in \P_s^{\iota\tilde\blam}\atop \bnu\text{ even}}\prod_{\vphi\in \Theta\atop \vphi=\tilde\vphi}\frac{\omega^{\tilde\blam^{(\tilde\vphi)}}(\bnu^{(\vphi)}/2)}{z_{\bnu^{(\vphi)}/2}}p_{\bnu^{(\vphi)}} \hspace{-.25cm} \prod_{\vphi\in \Theta\atop \vphi=\tilde\vphi\cup F\tilde\vphi} \biggl(\sum_{|\gamma|=|\tilde\blam^{(\tilde\vphi)}|\atop {|\mu|=|\tilde\blam^{(F\tilde\vphi)}| \atop \gamma\cup \mu = \bnu^{(\vphi)}} } \frac{\omega^{\tilde\blam^{(\tilde\vphi)}}(\gamma) \omega^{\tilde\blam^{(F\tilde\vphi)}}(\mu) }{z_{\gamma}z_{\mu}}\biggr) p_{\bnu^{(\vphi)}}\\
& = (-1)^n\hspace{-.15cm} \prod_{\vphi\in \Theta\atop \vphi=\tilde\vphi}\sum_{|\nu|=|\tilde\blam^{(\tilde\vphi)}|}\hspace{-.25cm}  \frac{\omega^{\tilde\blam^{(\tilde\vphi)}}(\nu)}{z_{\nu}}p_{2\nu}(Y^{(\vphi)}) \hspace{-.25cm}\prod_{\vphi\in \Theta\atop \vphi=\tilde\vphi\cup F\tilde\vphi} \sum_{|\nu|=|\tilde\blam^{(\tilde\vphi)}|+|\tilde\blam^{(F\tilde\vphi)}|} \hspace{-.05cm}\biggl(\hspace{-.25cm}\sum_{|\gamma|=|\tilde\blam^{(\tilde\vphi)}|\atop {|\mu|=|\tilde\blam^{(F\tilde\vphi)}| \atop \gamma\cup \mu = \nu} } \hspace{-.35cm} \frac{\omega^{\tilde\blam^{(\tilde\vphi)}}(\gamma) \omega^{\tilde\blam^{(F\tilde\vphi)}}(\mu) }{z_{\gamma}z_{\mu}}\biggr) p_\nu(Y^{(\vphi)}).
\end{align*}
The first product is the case that $|\vphi|$ is odd, and the second product is the case that $|\vphi|$ is even.  For the sum in the first product, note that
\begin{align*}
\sum_{|\nu|=|\lambda|+|\eta|} \biggl(\sum_{|\gamma|=|\lambda|\atop {|\mu|=|\eta| \atop \gamma\cup \mu = \nu} } \frac{\omega^{\lambda}(\gamma) \omega^{\eta}(\mu) }{z_{\gamma}z_{\mu}}\biggr) p_\nu&=\biggl( \sum_{|\gamma|=|\lambda|} \frac{\omega^\lambda(\gamma)}{z_\gamma} p_\gamma\biggr)\biggl(
 \sum_{|\mu|=|\eta|} \frac{\omega^\eta(\mu)}{z_\mu} p_\mu\biggr)\\
&=s_\lambda s_\eta. 
\end{align*}
For the sum in the product for $|\vphi|$ even, we have
\begin{align*}
\sum_{|\nu|=|\lambda|}  \frac{\omega^{\lambda}(\nu)}{z_{\nu}}p_{2\nu}&=\sum_{|\nu|=|\lambda|}  \frac{\omega^{\lambda}(\nu)}{z_{\nu}} p_\nu \circ p_{(2)}\\
&= s_\lambda \circ p_{(2)}
\end{align*}
where $\circ$ is the plethysm product (\ref{PlethysmCoefficients}).  Thus, from the definition of the coefficients $c_{\tilde\blam}^{\bgamma}$, we have
$$\ch\bigl(R_{G_n}^{U_{2n}}(\chi^{\tilde\blam})\bigr) =(-1)^n \sum_{\bgamma\in \cP^{\iota\tilde\blam}_s}c_{\tilde\blam}^{\bgamma} s_{\bgamma},$$
as desired.
\end{proof}

It is perhaps worth noting that since we know the Harish-Chandra induction $R_{G_n}^{U_{2n}}(\chi^{\tilde\blam})$ gives a character, then the sign of the coefficient $c_{\tilde\blam}^{\bgamma}$ must be $(-1)^{n(\bgamma)}$.

\subsection{Degenerate Gelfand-Graev characters}

The following theorem is our main theorem of Section \ref{SectionDegenerateGelfandGraev}.
\begin{theorem} \label{decomp} Let $n\in \ZZ_{\geq 1}$ and let $(k,\nu)$ satisfy $\nu\vdash \frac{n-k}{2}\in \ZZ_{\geq 0}$.  Then
$$\Gamma_{(k,\nu)}=\sum_{\blam\in \cP^\Theta_n} |\cB_{(k,\nu)}^{\blam}| \chi^{\blam}.$$
\end{theorem}

\begin{proof}
Recall that by Proposition \ref{DGGProductDecomposition}, 
$$\ch(\Gamma_{(k,\nu)})=\ch\bigl(\Gamma_{(k)}\bigr)\ch\biggl(R_{G_{\nu_1}}^{U_{2\nu_1}}(\tilde\Gamma_{(\nu_1)})\biggr)\ch\biggl(R_{G_{\nu_2}}^{U_{2\nu_2}}(\tilde\Gamma_{(\nu_2)})\biggr)\cdots \ch\biggl(R_{G_{\nu_\ell}}^{U_{2\nu_\ell}}(\tilde\Gamma_{(\nu_\ell)})\biggr).$$
By Theorem \ref{GGDecomposition} and Corollary \ref{DGGBaseCaseCorollary},
\begin{align*}
\ch(\Gamma_{(k)})&=(-1)^{\lfloor k/2 \rfloor}\sum_{\bgamma\in \cP^\Theta_k\atop \hgt(\bgamma)= 1} \prod_{\vphi\in \Theta} s_{\bgamma^{(\vphi)}}(Y^{(\vphi)}), \text{ and }\\
\ch(R_{G_{r}}^{U_{2r}}(\tilde\Gamma_{(r)}))&=(-1)^r\sum_{{\bgamma\in \cP^\Theta_{2r}\atop \hgt(\bgamma)= 1}\atop \bgamma\text{ even}} \prod_{\vphi\in \Theta} \sum_{i=0}^{|\bgamma^{(\vphi)}|} (-1)^{i|\vphi|} s_{(i)}(Y^{(\vphi)})s_{(|\bgamma^{(\vphi)}|-i)}(Y^{(\vphi)}).
\end{align*}
Thus,
$$\ch(\Gamma_{(k,\nu)})=(-1)^{\lfloor n/2 \rfloor}\sum_{\bgamma_0\in \cP^\Theta_k\atop \hgt(\bgamma_0)= 1}\sum_{ {{1 \leq r \leq \ell(\nu) \atop \bgamma_r\in \cP^\Theta_{2\nu_r}}\atop  \hgt(\bgamma_r)= 1}\atop \bgamma_r\text{ even}}\prod_{\vphi\in \Theta} s_{\bgamma_0^{(\vphi)}}(Y^{(\vphi)})\prod_{r=1}^{\ell(\nu)} \sum_{ i=0}^{|\bgamma_r^{(\vphi)}|} (-1)^{i|\vphi|} s_{(i)}(Y^{(\vphi)})s_{(|\bgamma_r^{(\vphi)}|-i)}(Y^{(\vphi)}).$$
Fix $\vphi\in \Theta$, and let $m_r=|\bgamma_r^{(\vphi)}|$.  If $|\vphi|$ is even, then Lemma \ref{PieriLemma} (a) implies
\begin{equation*}
s_{(m_0)}\prod_{r=1}^{\ell(\nu)}\sum_{i=0}^{m_r}  s_{(i)}s_{(m_r-i)}
=\sum_{|\lambda|=m_0+\cdots+m_{\ell(\nu)}} \bigl|\cT_{(m_0,\ldots,m_{\ell(\nu)})}^\lambda\bigr| s_\lambda.
\end{equation*}
If $|\vphi|$ is odd, then Lemma \ref{PieriLemma} (b) implies
$$s_{(m_0)}\prod_{r=1}^{\ell(\nu)} \sum_{i=0}^{m_r} (-1)^{i} s_{(i)}s_{(m_r-i)}=\sum_{|\lambda|=m_0+\cdots+m_{\ell(\nu)}} (-1)^{n(\lambda)}\bigl|\cD_{(m_0,m_1/2,\ldots,m_{\ell(\nu)}/2)}^\lambda\bigr| s_\lambda.$$
Therefore,
\begin{align*}
&\ch(\Gamma_{(k,\nu)})\\
&=(-1)^{\lfloor n/2 \rfloor} \sum_{\bgamma_0\in \cP^\Theta_k\atop \hgt(\bgamma_0)= 1}\sum_{ {{1 \leq r \leq \ell(\nu) \atop \bgamma_r\in \cP^\Theta_{2\nu_r}}\atop  \hgt(\bgamma_r)= 1}\atop \bgamma_r\text{ even}} \prod_{\vphi\in \Theta\atop|\vphi|\text{ odd}}\sum_{\blam^{(\vphi)}} (-1)^{n(\blam^{(\vphi)})}\bigl|\cD_{(|\bgamma_0^{(\vphi)}|, |\bgamma_1^{(\vphi)}|/2,\ldots)}^{\blam^{(\vphi)}}\bigr| s_{\blam^{(\vphi)}}(Y^{(\vphi)}) \\
&\hspace*{8cm}\cdot \prod_{\vphi\in \Theta\atop |\vphi|\text{ even}} \sum_{\blam^{(\vphi)}} \bigl|\cT_{(|\bgamma_0^{(\vphi)}|, |\bgamma_1^{(\vphi)}|\ldots)}^{\blam^{(\vphi)}}\bigr| s_{\blam^{(\vphi)}}(Y^{(\vphi)})\\
&=(-1)^{\lfloor n/2 \rfloor}\hspace{-.25cm} \sum_{\blam\in \cP_n^\Theta} \biggl(\sum_{\bgamma_0\in \cP^\Theta_k\atop \hgt(\bgamma_0)= 1}\sum_{ {{1 \leq r \leq \ell(\nu) \atop \bgamma_r\in \cP^\Theta_{2\nu_r}}\atop  \hgt(\bgamma_r)= 1}\atop \bgamma_r\text{ even}} \prod_{\vphi\in \Theta\atop|\vphi|\text{ odd}}(-1)^{n(\blam^{(\vphi)})}\bigl|\cD_{(|\bgamma_0^{(\vphi)}|, |\bgamma_1^{(\vphi)}|/2,\ldots)}^{\blam^{(\vphi)}}\bigr| \hspace{-.25cm} \prod_{\vphi\in \Theta\atop |\vphi|\text{ even}} \bigl|\cT_{(|\bgamma_0^{(\vphi)}|, |\bgamma_1^{(\vphi)}|\ldots)}^{\blam^{(\vphi)}}\bigr| \biggr)s_{\blam}\\
&=\sum_{\blam \in \cP_n^\Theta} (-1)^{\lfloor n/2 \rfloor +n(\blam)} |\cB_{(k,\nu)}^{\blam}|s_{\blam}.
\end{align*}

\noindent
The result follows by applying $\ch^{-1}$.
\end{proof}

\section{Some multiplicity consequences}

In this section we explore some of the multiplicity implications of Theorem \ref{decomp}.

\subsection{A bijection between domino tableaux and pairs of column strict tableaux}

The \emph{2-core} of a partition $\lambda\in \cP$, which we denote $\core_2(\lambda)$, is the partition of minimal size such that the skew partition $\lambda / \core_2(\lambda)$ may be tiled by dominoes.  It is not difficult to see that the 2-core of any partition is always of the form $(m, m-1, \ldots, 2, 1)$ for some nonnegative integer $m$ (where $(0)$ is the empty partition).  

The \emph{$2$-quotient} of a partition $\lambda$, $\quot_2(\lambda)$, is a pair of partitions $(\quot_2(\lambda)^{(0)}, \quot_2(\lambda)^{(1)})$ (defined in \cite[I.1, Example 8]{Mac}).   We define
$$\quot_2(\lambda)_i = \quot_2(\lambda)^{(0)}_i + \quot_2(\lambda)^{(1)}_i.$$
Also define the \emph{content} of a box $\Box$ in the $i$th row and $j$th column of a partition $\lambda$ to be $j-i$.

Let $\lambda\in \cP_n$ with $\core_2(\lambda)\in \{(0),(1)\}$.  Consider the bijection
\begin{equation}\label{DominoToTableau}\begin{array}{rccc} & \cD^\lambda_{(|\core_2(\lambda)|,m_1,\ldots,m_\ell)} & \longleftrightarrow & \left\{\begin{array}{c} \text{Pairs of column strict}\\ \text{tableaux of shape $\mathrm{quot}_2(\lambda)$}\\ \text{and weight $(m_1,m_2,\ldots, m_\ell)$}\end{array}\right\}\\  & Q & \leftrightarrow & (Q^{(0)}, Q^{(1)}),\end{array}\end{equation}
given by the following algorithm, which originally appeared in \cite{SW}, and is in a more general form in \cite{LLT95}.
\begin{enumerate}
\item[(1)]  Each domino in $Q$ covers two boxes of $\lambda/\core_2(\lambda)$.  Move the entries in $Q$ to the box that has content $0$ modulo 2. 
$$Q=\xy<0cm,.75cm>\xymatrix@R=.3cm@C=.3cm{
*={} & *={} \ar @{-} [l] & *={} \ar @{-} [l] & *={} \ar @{-} [l] & *={} \ar @{-} [l] & *={} \ar @{-} [l]  \\
*={} \ar @{-} [u] &   *={} \ar @{-} [l] \ar @{-} [u] &   *={} \ar @{} [l]|{\sss 1} \ar@{-} [u]  &   *={}  \ar @{-}[l] \ar@{} [u]|{\sss 1}  &   *={} \ar @{-} [l]  \ar @{-} [u] &   *={} \ar @{} [l]|{\sss 2}  \ar @{-} [u] \\ 
*={} \ar @{-} [u] &   *={} \ar @{} [l]|{\sss 1} \ar @{-} [u]  &   *={} \ar @{-} [l] \ar@{-} [u]   &   *={} \ar @{-} [l] \ar@{} [u]|{\sss 2} &   *={} \ar @{-} [l] \ar@{-} [u]  &   *={} \ar @{-} [l] \ar@{-} [u]\\
*={} \ar @{-} [u] &   *={} \ar@{-} [u] \ar @{-} [l]    &   *={} \ar @{-} [l] \ar@{} [u]|{\sss 3} &   *={} \ar @{-} [l] \ar@{-} [u] &   *={} \ar @{-} [l] \ar@{} [u]|{\sss 3} &   *={} \ar @{-} [l] \ar@{-} [u]  \\
*={} \ar @{-} [u] &   *={} \ar @{-} [l] \ar@{} [u]|{\sss 4}   &   *={} \ar @{-} [l] \ar@{-} [u]  &   *={} \ar @{} [l]|{\sss 6} \ar@{-} [u] \\
*={} \ar @{-} [u] &   *={} \ar @{-} [l] \ar@{} [u]|{\sss 5}    &   *={} \ar @{-} [l] \ar@{-} [u] &   *={} \ar @{-} [l] \ar@{-} [u] }\endxy
\quad\longmapsto\quad 
\xy<0cm,.75cm>\xymatrix@R=.3cm@C=.3cm{
*={} & *={} \ar @{-} [l] & *={} \ar @{-} [l] & *={} \ar @{-} [l] & *={} \ar @{-} [l] & *={} \ar @{-} [l]  \\
*={} \ar @{-} [u] &   *={} \ar @{-} [l] \ar @{-} [u] &   *={} \ar @{.} [l] \ar@{-} [u]  &   *={}  \ar @{-}[l] \ar@{.} [u] \ar @{} [ul]|{\sss 1} &   *={} \ar @{-} [l]  \ar @{-} [u] &   *={} \ar @{.} [l]  \ar @{-} [u] \ar @{} [ul]|{\sss 2} \\ 
*={} \ar @{-} [u] &   *={} \ar @{.} [l] \ar @{-} [u]  &   *={} \ar @{-} [l] \ar@{-} [u] \ar @{} [ul]|{\sss 1}  &   *={} \ar @{-} [l] \ar@{.} [u] &   *={} \ar @{-} [l] \ar@{-} [u] \ar @{} [ul]|{\sss 2}  &   *={} \ar @{-} [l] \ar@{-} [u]\\
*={} \ar @{-} [u] &   *={} \ar@{-} [u] \ar @{-} [l]  \ar @{} [ul]|{\sss 1}  &   *={} \ar @{-} [l] \ar@{.} [u] &   *={} \ar @{-} [l] \ar@{-} [u] \ar @{} [ul]|{\sss 3} &   *={} \ar @{-} [l] \ar@{.} [u] &   *={} \ar @{-} [l] \ar@{-} [u] \ar @{} [ul]|{\sss 3} \\
*={} \ar @{-} [u] &   *={} \ar @{-} [l] \ar@{.} [u]   &   *={} \ar @{-} [l] \ar@{-} [u] \ar @{} [ul]|{\sss 4} &   *={} \ar @{.} [l] \ar@{-} [u] \\
*={} \ar @{-} [u] &   *={} \ar @{-} [l] \ar@{.} [u]   \ar @{} [ul]|{\sss 5} &   *={} \ar @{-} [l] \ar@{-} [u] &   *={} \ar @{-} [l] \ar@{-} [u] \ar @{} [ul]|{\sss 6}}\endxy$$
\item[(2)]  Let $\cS^{(0)}$ denote the set of all dominoes that have the entry in the lower or leftmost box, and $\cS^{(1)}$ be the set of dominoes that have the entry in the upper or rightmost box.
$$\cS^{(0)}=\left\{\xy<0cm,.75cm>\xymatrix@R=.3cm@C=.3cm{
*={} & *={} \ar @{-} [l] & *={} \ar @{-} [l] & *={} \ar @{-} [l] & *={} \ar @{-} [l] & *={} \ar @{-} [l]  \\
*={} \ar @{-} [u] &   *={} \ar @{-} [l] \ar @{-} [u] &   *={} \ar @{.} [l] \ar@{-} [u]  &   *={}  \ar @{-}[l] \ar@{.} [u] \ar @{} [ul]|{\sss 1} &   *={} \ar @{-} [l]  \ar @{-} [u] &   *={} \ar @{.} [l]  \ar @{-} [u]  \\ 
*={} \ar @{-} [u] &   *={} \ar @{.} [l] \ar @{-} [u]  &   *={} \ar @{-} [l] \ar@{-} [u] \ar @{} [ul]|{\sss 1}  &   *={} \ar @{-} [l] \ar@{.} [u] &   *={} \ar @{-} [l] \ar@{-} [u]   &   *={} \ar @{-} [l] \ar@{-} [u]\\
*={} \ar @{-} [u] &   *={} \ar@{-} [u] \ar @{-} [l]  \ar @{} [ul]|{\sss 1}  &   *={} \ar @{-} [l] \ar@{.} [u] &   *={} \ar @{-} [l] \ar@{-} [u] &   *={} \ar @{-} [l] \ar@{.} [u] &   *={} \ar @{-} [l] \ar@{-} [u]  \\
*={} \ar @{-} [u] &   *={} \ar @{-} [l] \ar@{.} [u]   &   *={} \ar @{-} [l] \ar@{-} [u]  &   *={} \ar @{.} [l] \ar@{-} [u] \\
*={} \ar @{-} [u] &   *={} \ar @{-} [l] \ar@{.} [u]   \ar @{} [ul]|{\sss 5} &   *={} \ar @{-} [l] \ar@{-} [u] &   *={} \ar @{-} [l] \ar@{-} [u] \ar @{} [ul]|{\sss 6}}\endxy\right\}
\qquad \text{and}\qquad 
\cS^{(1)}=\left\{\xy<0cm,.75cm>\xymatrix@R=.3cm@C=.3cm{
*={} & *={} \ar @{-} [l] & *={} \ar @{-} [l] & *={} \ar @{-} [l] & *={} \ar @{-} [l] & *={} \ar @{-} [l]  \\
*={} \ar @{-} [u] &   *={} \ar @{-} [l] \ar @{-} [u] &   *={} \ar @{.} [l] \ar@{-} [u]  &   *={}  \ar @{-}[l] \ar@{.} [u]  &   *={} \ar @{-} [l]  \ar @{-} [u] &   *={} \ar @{.} [l]  \ar @{-} [u] \ar @{} [ul]|{\sss 2} \\ 
*={} \ar @{-} [u] &   *={} \ar @{.} [l] \ar @{-} [u]  &   *={} \ar @{-} [l] \ar@{-} [u] &   *={} \ar @{-} [l] \ar@{.} [u] &   *={} \ar @{-} [l] \ar@{-} [u] \ar @{} [ul]|{\sss 2}  &   *={} \ar @{-} [l] \ar@{-} [u]\\
*={} \ar @{-} [u] &   *={} \ar@{-} [u] \ar @{-} [l]    &   *={} \ar @{-} [l] \ar@{.} [u] &   *={} \ar @{-} [l] \ar@{-} [u] \ar @{} [ul]|{\sss 3} &   *={} \ar @{-} [l] \ar@{.} [u] &   *={} \ar @{-} [l] \ar@{-} [u] \ar @{} [ul]|{\sss 3} \\
*={} \ar @{-} [u] &   *={} \ar @{-} [l] \ar@{.} [u]   &   *={} \ar @{-} [l] \ar@{-} [u] \ar @{} [ul]|{\sss 4} &   *={} \ar @{.} [l] \ar@{-} [u] \\
*={} \ar @{-} [u] &   *={} \ar @{-} [l] \ar@{.} [u]   &   *={} \ar @{-} [l] \ar@{-} [u] &   *={} \ar @{-} [l] \ar@{-} [u]}\endxy\right\}$$
\item[(3)]  For even $-\ell(\lambda)<i<\lambda_1$ and $j\in \{0,1\}$, let
$$D_i^{(j)}=\begin{array}{c}\text{The increasing sequence of entries whose}\\  \text{content is $i$ and whose domino is in $\cS^{(j)}$.}\end{array}$$
\begin{align*}
\big(D_{-4}^{(0)},D_{-2}^{(0)},D_{0}^{(0)},D_{2}^{(0)},D_{4}^{(0)}\big)&=\big((5),(1,6),(1),(1),()\big)\\
\big(D_{-4}^{(1)},D_{-2}^{(1)},D_{0}^{(1)},D_{2}^{(1)},D_{4}^{(1)}\big)&=\big((),(4),(3),(2,3),(2)\big)
\end{align*}
\item[(4)] Let $Q^{(j)}$ be the unique tableau that has increasing diagonal sequences given by the $D_i^{(j)}$  for all even  $-\ell(\lambda)<i<\lambda_1$.
$$Q^{(0)}=
\xy<0cm,.3cm>\xymatrix@R=.3cm@C=.3cm{
*={} & *={} \ar @{-} [l] & *={} \ar @{-} [l]  & *={} \ar @{-} [l] \\
*={} \ar @{-} [u] &   *={} \ar @{-} [l] \ar @{-} [u] \ar @{} [ul]|{\sss 1} &   *={} \ar @{-} [l] \ar@{-} [u] \ar @{} [ul]|{\sss 1}&   *={} \ar @{-} [l] \ar@{-} [u] \ar @{} [ul]|{\sss 1} \\
*={} \ar @{-} [u] &   *={} \ar @{-} [l] \ar @{-} [u] \ar @{} [ul]|{\sss 5} &   *={} \ar @{-} [l] \ar@{-} [u] \ar @{} [ul]|{\sss 6}}\endxy
\qquad \text{and}\qquad 
Q^{(1)}=\xy<0cm,.45cm>\xymatrix@R=.3cm@C=.3cm{
*={} & *={} \ar @{-} [l] & *={} \ar @{-} [l]   \\
*={} \ar @{-} [u] &   *={} \ar @{-} [l] \ar @{-} [u] \ar @{} [ul]|{\sss 2} &   *={} \ar @{-} [l] \ar@{-} [u] \ar @{} [ul]|{\sss 2} \\
*={} \ar @{-} [u] &   *={} \ar @{-} [l] \ar @{-} [u] \ar @{} [ul]|{\sss 3} &   *={} \ar @{-} [l] \ar@{-} [u] \ar @{} [ul]|{\sss 3}\\
*={} \ar @{-} [u] &   *={} \ar @{-} [l] \ar @{-} [u] \ar @{} [ul]|{\sss 4} 
}\endxy$$
\end{enumerate}

\noindent\textbf{Remarks.} 
\begin{enumerate}
\item If the shape of the domino tableau $Q$ is $\lambda/\core_2(\lambda)$, then the shape of $(Q^{(0)},Q^{(1)})$ is $\quot_2(\lambda)$.
\item We may apply this algorithm to a domino tableau of shape $\lambda/(m)$ with $m\equiv |\core_2(\lambda)|\mod 2$, by requiring that the tableau of shape $\lambda/\core_2(\lambda)$ has $\lfloor m/2\rfloor$ horizontal dominoes filled with zeroes.  For example, 
$$Q =\xy<0cm,.75cm>\xymatrix@R=.3cm@C=.3cm{
*={} & *={} \ar @{-} [l] & *={} \ar @{-} [l] & *={} \ar @{-} [l] & *={} \ar @{-} [l] & *={} \ar @{-} [l]  \\
*={} \ar @{-} [u] &   *={} \ar @{-} [l] \ar @{-} [u] &   *={} \ar @{-} [l] \ar@{-} [u]  &   *={}  \ar @{-}[l] \ar@{-} [u]  &   *={} \ar @{-} [l]  \ar @{-} [u] &   *={} \ar @{-} [l]  \ar @{-} [u]  \\ 
*={} \ar @{-} [u] &   *={} \ar @{.} [l] \ar @{-} [u]  &   *={} \ar @{.} [l] \ar@{-} [u] \ar @{} [ul]|{\sss 1}  &   *={} \ar @{-} [l] \ar@{.} [u] &   *={} \ar @{-} [l] \ar@{-} [u] \ar @{} [ul]|{\sss 2}  &   *={} \ar @{.} [l] \ar@{-} [u]\\
*={} \ar @{-} [u] &   *={} \ar@{-} [u] \ar @{-} [l]  \ar @{} [ul]|{\sss 1}  &   *={} \ar @{-} [l] \ar@{-} [u] &   *={} \ar @{-} [l] \ar@{.} [u] \ar @{} [ul]|{\sss 3} &   *={} \ar @{-} [l] \ar@{-} [u] &   *={} \ar @{-} [l] \ar@{-} [u] \ar @{} [ul]|{\sss 3} \\
*={} \ar @{-} [u] &   *={} \ar @{-} [l] \ar@{.} [u]   &   *={} \ar @{-} [l] \ar@{-} [u] \ar @{} [ul]|{\sss 4} &   *={} \ar @{.} [l] \ar@{-} [u] \\
*={} \ar @{-} [u] &   *={} \ar @{-} [l] \ar@{.} [u]   \ar @{} [ul]|{\sss 5} &   *={} \ar @{-} [l] \ar@{-} [u] &   *={} \ar @{-} [l] \ar@{-} [u] \ar @{} [ul]|{\sss 6}}\endxy\quad\text{has shape $\lambda/(5)$, so apply the algorithm to}\quad 
\xy<0cm,.75cm>\xymatrix@R=.3cm@C=.3cm{
*={} & *={} \ar @{-} [l] & *={} \ar @{-} [l] & *={} \ar @{-} [l] & *={} \ar @{-} [l] & *={} \ar @{-} [l]  \\
*={} \ar @{-} [u] &   *={} \ar @{-} [l] \ar @{-} [u] &   *={} \ar @{-} [l] \ar@{.} [u]  &   *={}  \ar @{-}[l] \ar@{-} [u] \ar @{} [ul]|{\sss 0}  &   *={} \ar @{-} [l]  \ar @{.} [u] &   *={} \ar @{-} [l]  \ar @{-} [u] \ar @{} [ul]|{\sss 0}  \\ 
*={} \ar @{-} [u] &   *={} \ar @{.} [l] \ar @{-} [u]  &   *={} \ar @{.} [l] \ar@{-} [u] \ar @{} [ul]|{\sss 1}  &   *={} \ar @{-} [l] \ar@{.} [u] &   *={} \ar @{-} [l] \ar@{-} [u] \ar @{} [ul]|{\sss 2}  &   *={} \ar @{.} [l] \ar@{-} [u]\\
*={} \ar @{-} [u] &   *={} \ar@{-} [u] \ar @{-} [l]  \ar @{} [ul]|{\sss 1}  &   *={} \ar @{-} [l] \ar@{-} [u] &   *={} \ar @{-} [l] \ar@{.} [u] \ar @{} [ul]|{\sss 3} &   *={} \ar @{-} [l] \ar@{-} [u] &   *={} \ar @{-} [l] \ar@{-} [u] \ar @{} [ul]|{\sss 3} \\
*={} \ar @{-} [u] &   *={} \ar @{-} [l] \ar@{.} [u]   &   *={} \ar @{-} [l] \ar@{-} [u] \ar @{} [ul]|{\sss 4} &   *={} \ar @{.} [l] \ar@{-} [u] \\
*={} \ar @{-} [u] &   *={} \ar @{-} [l] \ar@{.} [u]   \ar @{} [ul]|{\sss 5} &   *={} \ar @{-} [l] \ar@{-} [u] &   *={} \ar @{-} [l] \ar@{-} [u] \ar @{} [ul]|{\sss 6}}\endxy$$
Note that all of the zero dominoes are in the same set $\cS^{(|\core_2(\lambda)|)}$, so changing $m$ corresponds to adding or subtracting the number of zeroes in the first row of $Q^{(|\core_2(\lambda)|)}$.
\end{enumerate}

We will use the lexicographic total ordering on partitions given by
\begin{equation}\label{LexicographicOrdering} \lambda\leq \mu\text{ if there exists $k\in \ZZ_{\geq 1}$ such that  $\lambda_k<\mu_k$ and $\lambda_i=\mu_i$ for $1\leq i< k$.}\end{equation}

\begin{lemma} \label{dominobij}
Let $\lambda \in \cP_n$ be such that $\core_2(\lambda)\in \{(0),(1)\}$, and let $0\leq m\leq \lambda_1$ be such that $m\equiv |\core_2(\lambda)|\mod 2$.  Then there exists a lexicographically maximal weight  $\mu=(\mu_1,\mu_2,\ldots,\mu_\ell)$ such that there exists exactly one domino tableau of shape $\lambda/(m)$ and weight $(m, \mu_1, \ldots, \mu_{\ell})$.
\end{lemma}

\begin{proof}
First suppose $(\lambda^{(0)}/\gamma^{(0)}, \lambda^{(1)}/\gamma^{(1)})$ is a pair of skew partitions.  Let $\mu_1$ be the maximal number of 1's we can put in a tableau of shape $(\lambda^{(0)}/\gamma^{(0)}, \lambda^{(1)}/\gamma^{(1)})$, $\mu_2$ be the maximal number of 2's we can thereafter fill into $(\lambda^{(0)}/\gamma^{(0)}, \lambda^{(1)}/\gamma^{(1)})$, and $\mu_j$ be the maximal number of $j$'s we can fill given that we have filled in a maximum number at each step up to $j$.  Then
there is exactly one tableau $(Q^{(0)}, Q^{(1)})$ of shape $(\lambda^{(0)}/\gamma^{(0)}, \lambda^{(1)}/\gamma^{(1)})$ and weight $\mu$, and this weight is lexicographically maximal.  The result now follows from pulling back $(Q^{(0)}, Q^{(1)})$ through the bijection (\ref{DominoToTableau}) to get a domino tableau of the same weight, along with the second remark preceding this Lemma.
\end{proof}

\noindent\textbf{Remark.}  If $m=\core_2(\lambda)$, then $\mu$ is given by $\mu_0=|\core_2(\lambda)|$ and $\mu_i=\quot_2(\lambda)_i$ for $i\geq 1$.

\subsection{Multiplicity results}

Our first consequence of Theorem \ref{decomp} characterizes which irreducible characters of ${\rm U}(n,\FF_{q^2})$ appear with nonzero multiplicity in some degenerate Gelfand-Graev character.

\begin{corollary} \label{decompcor} 
The set of all $\blam \in \cP^{\Theta}_n$ such that the character $\chi^{\blam}$ of ${\rm U}(n, \FF_{q^2})$ satisfies $\langle \chi^{\blam}, \Gamma_{(k, \nu)} \rangle \neq 0$ for some degenerate Gelfand-Graev character $\Gamma_{(k,\nu)}$ is
$$ \{ \blam \in \cP^{\Theta}_n \; \big| \; \core_2(\blam^{(\varphi)}) \in \{(0),(1)\} \text{ whenever } |\varphi| \text{ is odd } \}.$$
\end{corollary}

\begin{proof}
By Theorem \ref{decomp}, the irreducible character $\chi^{\blam}$ appears with nonzero multiplicity in some degenerate Gelfand-Graev character if and only if there a battery tableau of shape $\blam/\bgamma$ for some $\bgamma\in \cP^\Theta$ with $\hgt(\bgamma)\leq 1$.

If for some odd $\vphi\in \Theta$, we have $\core_2(\blam^{(\vphi)})\notin\{(0),(1)\}$,  then the 2-core of $\blam^{(\varphi)}$ has at least two parts.   But then there is no choice of $\bgamma^{(\varphi)}$ that allows us to tile $\blam^{(\varphi)} / \bgamma^{(\varphi)}$ by dominoes.  On the other hand, if  $\core_2(\blam^{(\varphi)})=(0)$, we can choose $\bgamma^{(\varphi)}=(0)$, and if $\core_2(\blam^{(\varphi)})=(1)$, we can let $\bgamma^{(\varphi)} = (1)$, and $\blam^{(\varphi)}/\bgamma^{(\varphi)}$ can be tiled by dominoes.   \end{proof}

A character $\chi$ of ${\rm U}(n, \FF_{q^2})$ is called {\em regular} if it is a component of the Gelfand-Graev character $\Gamma_{(n)}$.  The next result follows directly from Corollary \ref{decompcor} and \cite[Corollary 2.6]{kot}.  For the definition of Harish-Chandra series, see \cite{Ca85}.

\begin{corollary} \label{kotlarcor}
The set all $\blam \in \cP^{\Theta}_n$ such that the character $\chi^{\blam}$ of ${\rm U}(n, \FF_{q^2})$ is in the union of all Harish-Chandra series that contain regular characters is
$$ \{ \blam \in \cP^{\Theta}_n \; \big| \; \core_2(\blam^{(\varphi)}) \in \{(0),(1)\} \text{ whenever } |\varphi| \text{ is odd } \}.$$
\end{corollary}

We now specify multiplicities of certain characters $\chi^{\blam}$ in degenerate Gelfand-Graev characters. 

\begin{theorem} \label{bigmult}
Let $\blam\in \cP_n^\Theta$ be such that $\core_2(\blam^{(\varphi)}) \in \{ (0), (1) \}$ whenever $|\varphi|$ is odd.  Then there exists $\nu \vdash \frac{n-k}{2}$ such that 
$$\langle \Gamma_{(k,\nu)},\chi^{\blam}\rangle=\prod_{\vphi\in \Theta\atop |\vphi|\text{ even}} \prod_{i\text{ odd}} \bigl(\blam^{(\vphi)}_i-\blam^{(\vphi)}_{i+1}+1\bigr).$$
\end{theorem}

\begin{proof}
Let $k=\sum_{|\vphi|\text{ odd}}|\vphi||\core_2(\blam^{(\vphi)})|$ and define $\bgamma$ by
$$ \bgamma^{(\vphi)} =\left\{\begin{array}{ll}\core_2(\blam^{(\vphi)}) & \text{if $|\vphi|$ is odd,}\\ \emptyset & \text{otherwise.}\end{array}\right.$$
Since $|\bgamma|=k$, by Theorem \ref{decomp} and Corollary \ref{decompcor}, it suffices to find $\nu\vdash \frac{n-k}{2}$ such that there exist 
$$\prod_{\vphi\in \Theta\atop |\vphi|\text{ even}} \prod_{i\text{ odd}} \bigl(\blam^{(\vphi)}_i-\blam^{(\vphi)}_{i+1}+1\bigr)$$
battery tableaux with shape $\blam/\bgamma$ and weight $(k,\nu)$.   

We construct the battery tableau ${\bf Q}$ as follows.  For odd $\varphi\in \Theta$, let ${\bf Q}^{(\varphi)}$ be the unique domino tableau of shape $\blam^{(\varphi)}/(\core_2(\blam^{(\varphi)}))$ and weight $(|\core_2(\blam^{(\varphi)})|, \quot_2(\blam^{(\varphi)})_1, \quot_2(\blam^{(\varphi)})_2, \ldots)$, obtained from Lemma \ref{dominobij} (see, in particular, the remark after the lemma).

For even $\varphi\in \Theta$ and for each $i \geq 1$, we fill the $(2i-1)$st row of $\blam^{(\varphi)}$ with $\bar{i}$'s, and the $2i$th row with $i$'s.  The resulting symplectic tableau ${\bf Q}^{(\varphi)}$ has weight $(\blam^{(\varphi)}_1 + \blam^{(\varphi)}_2, \blam^{(\varphi)}_3 + \blam^{(\varphi)}_4, \ldots)$.  Note that we may change up to $\blam^{(\varphi)}_{2i-1} - \blam^{(\varphi)}_{2i}$ of the $\bar{i}$'s to $i$'s in row $2i-1$ while leaving the weight unchanged.  We therefore have exactly
$$ \prod_{i \text{ odd}}(\blam^{(\varphi)}_i - \blam^{(\vphi)}_{i-1} +1)$$
symplectic tableaux of shape $\blam^{(\vphi)}$ and weight $(\blam^{(\varphi)}_1 + \blam^{(\varphi)}_2, \blam^{(\varphi)}_3 + \blam^{(\varphi)}_4, \ldots)$.

We combine these to create a battery tableau of shape $\blam/\bgamma$ and weight $\nu$, where
$$\nu_i=\sum_{\vphi\in \Theta\atop |\vphi|\text{ odd}} |\vphi|\quot_2(\blam^{(\vphi)})_i+\sum_{\vphi\in \Theta\atop |\vphi|\text{ even}} \frac{|\vphi|}{2} \bigl(\blam_{2i}^{(\vphi)}+\blam_{2i-1}^{(\vphi)}\bigr).$$ 
Note that from this construction, $\nu$ is the maximal weight under the lexicographical ordering  (\ref{LexicographicOrdering}) of a battery tableau of shape $\blam/\bgamma$, while each $\bgamma^{(\varphi)}$ is chosen minimally.  It follows that the weight $\nu$ will change if we change the weight of any ${\bf Q}^{(\vphi)}$.

For example, if
$$\blam=\left(\xy<0cm,.5cm>\xymatrix@R=.3cm@C=.3cm{
	*={} & *={} \ar @{-} [l] & *={} \ar @{-} [l]  \\
	*={} \ar @{-} [u] &   *={} \ar @{-} [l] \ar@{-} [u] &   *={} \ar @{-} [l] \ar@{-} [u]   \\
*={} \ar @{-} [u] &   *={} \ar @{-} [l] \ar@{-} [u]  & *={} \ar @{-} [l] \ar@{-} [u] \\ 
*={} \ar @{-} [u] &   *={} \ar @{-} [l] \ar@{-} [u]  }\endxy^{(\vphi_1)},\ 
\xy<0cm,.5cm>\xymatrix@R=.3cm@C=.3cm{
	*={} & *={} \ar @{-} [l] & *={} \ar @{-} [l] & *={} \ar @{-} [l]  \\
	*={} \ar @{-} [u] &   *={} \ar @{-} [l] \ar@{-} [u] &   *={} \ar @{-} [l] \ar@{-} [u]  &   *={} \ar @{-} [l] \ar@{-} [u]  \\ *={} \ar @{-} [u] &   *={} \ar @{-} [l] \ar@{-} [u] \\ *={} \ar @{-} [u] &   *={} \ar @{-} [l] \ar@{-} [u]  }\endxy^{(\vphi_2)}\  ,\  \xy<0cm,.5cm>\xymatrix@R=.3cm@C=.3cm{
	*={} & *={} \ar @{-} [l] & *={} \ar @{-} [l] \\
	*={} \ar @{-} [u] &   *={} \ar @{-} [l] \ar@{-} [u]  &  *={} \ar @{-} [l] \ar@{-} [u] \\ *={} \ar @{-} [u] &   *={} \ar @{-} [l] \ar@{-} [u] \\ *={} \ar @{-} [u] &   *={} \ar @{-} [l] \ar@{-} [u] }\endxy^{(\vphi_3)}\ \right)\qquad\text{where $|\vphi_i|=i$},$$
then $k=1$, $\nu=(2+4+3, 0+1+3)=(9,4)$, and every battery tableau of shape $\blam$ and weight $(k,\nu)$ must be of the form 
$$\left(\xy<0cm,.5cm>\xymatrix@R=.3cm@C=.3cm{
	*={} & *={} \ar @{-} [l] & *={} \ar @{-} [l]  \\
	*={} \ar @{-} [u] &   *={} \ar @{-} [l] \ar@{-} [u] \ar @{} [ul]|{\ss 0} &   *={} \ar @{} [l]|{\ss 1} \ar@{-} [u]   \\
*={} \ar @{-} [u] &   *={} \ar @{} [l]|{1} \ar@{-} [u]  & *={} \ar @{-} [l] \ar@{-} [u] \\ 
*={} \ar @{-} [u] &   *={} \ar @{-} [l] \ar@{-} [u]  }\endxy^{(\vphi_1)},\ 
\xy<0cm,.5cm>\xymatrix@R=.3cm@C=.3cm{
	*={} & *={} \ar @{-} [l] & *={} \ar @{-} [l] & *={} \ar @{-} [l]  \\
	*={} \ar @{-} [u] &   *={} \ar @{-} [l] \ar@{-} [u]\ar @{} [ul]|{\ss \bar{1}} &   *={} \ar @{-} [l] \ar@{-} [u]  \ar @{} [ul]|{\ss \ddot 1} &   *={} \ar @{-} [l] \ar@{-} [u] \ar @{} [ul]|{\ss \ddot 1} \\ *={} \ar @{-} [u] &   *={} \ar @{-} [l] \ar@{-} [u] \ar @{} [ul]|{\ss 1}\\ *={} \ar @{-} [u] &   *={} \ar @{-} [l] \ar@{-} [u] \ar @{} [ul]|{\ss \ddot 2} }\endxy^{(\vphi_2)}\  ,\  \xy<0cm,.5cm>\xymatrix@R=.3cm@C=.3cm{
	*={} & *={} \ar @{-} [l] & *={} \ar @{-} [l] \\
	*={} \ar @{-} [u] &   *={} \ar @{-} [l] \ar@{} [u]|{\ss 1}  &  *={} \ar @{-} [l] \ar@{-} [u] \\ *={} \ar @{-} [u] &   *={} \ar @{} [l]|{\ss 2} \ar@{-} [u] \\ *={} \ar @{-} [u] &   *={} \ar @{-} [l] \ar@{-} [u] }\endxy^{(\vphi_3)}\ \right),\qquad \text{where $\ddot i\in \{\bar{i},i\}$.}$$
There are $3\cdot 2=6$ of such tableaux.  
\end{proof}

Theorem \ref{bigmult} and its proof give the following multiplicity one result.

\begin{corollary}\label{OurMultiplicity}
Let $\blam\in \cP^\Theta_n$.  Suppose $\core_2(\blam^{(\vphi)})\in \{(0),(1)\}$ whenever $|\vphi|$ is odd, and for all $i> 0$, the multiplicities $m_i(\blam^{(\vphi)})$ are even for all  even $\vphi\in\Theta$.  Then $\langle \Gamma_{(k,\nu)},\chi^{\blam}\rangle=1$, where 
$$k=\sum_{\vphi\in \Theta\atop |\vphi|\text{ odd}} |\vphi||\core_2(\blam^{(\vphi)})| \quad \text{and}\quad \nu_i=\sum_{\vphi\in \Theta\atop |\vphi|\text{ odd}} |\vphi|\quot_2(\blam^{(\vphi)})_i+\sum_{\vphi\in \Theta\atop |\vphi|\text{ even}} \frac{|\vphi|}{2} \bigl(\blam_{2i}^{(\vphi)}+\blam_{2i-1}^{(\vphi)}\bigr).$$
\end{corollary}

The next theorem generalizes Corollary \ref{OurMultiplicity}.

\begin{theorem}\label{MultiplicityOne} Suppose $\blam\in \cP^\Theta_n$ satisfies the following:
\begin{enumerate}
\item[(a)] Whenever $|\vphi|$ is odd, $\core_2(\blam^{(\vphi)}) \in \{ (0), (1) \}$,
\item[(b)] Whenever $|\vphi|$ is even, the partition $\blam^{(\vphi)}$ has at most one nonzero part with odd multiplicity,
\item[(c)] There exists an $r >0$ such that for every $\vphi\in \Theta$ with $|\vphi|$ even, either $\ell(\blam^{(\vphi)})< r$ and $\blam^{(\vphi)}$ has no nonzero part of odd multiplicity, or $\blam^{(\vphi)}_r$ has odd multiplicity and $\blam^{(\vphi)}_r < \blam^{(\vphi)}_{r-1}$.
\end{enumerate}
Then the irreducible character $\chi^{\blam}$ of $U_n$ appears with multiplicity one in some degenerate Gelfand-Graev character of $U_n$.
\end{theorem}

\begin{proof} 
Suppose $\blam\in \cP_n^\Theta$ satisfies (a), (b), and (c).   Let $\bgamma\in \cP^\Theta$ be given by
$$\bgamma^{(\vphi)}=\left\{\begin{array}{ll} (i) & \text{if $|\vphi|$ is even and $m_i(\blam^{(\vphi)})$ is odd,}\\ \Big(|\core_2(\lambda)|+2\quot_2(\lambda)^{(|\core_2(\lambda)|)}_{\lceil r/2\rceil}\Big) & \text{if $|\vphi|$ is odd and $\lambda=\blam^{(\vphi)}$.}\end{array}\right.$$
For example, the $\Theta$-partition 
$$\blam=\left(\xy<0cm,.75cm>\xymatrix@R=.3cm@C=.3cm{
*={} & *={} \ar @{-} [l] & *={} \ar @{-} [l] & *={} \ar @{-} [l] & *={} \ar @{-} [l] & *={} \ar @{-} [l]  \\
*={} \ar @{-} [u] &   *={} \ar @{-} [l] \ar @{-} [u] &   *={} \ar @{-} [l] \ar@{-} [u]  &   *={}  \ar @{-}[l] \ar@{-} [u]  &   *={} \ar @{-} [l]  \ar @{-} [u] &   *={} \ar @{-} [l]   \ar @{-} [u]\\ 
*={} \ar @{-} [u] &   *={} \ar @{-} [l] \ar @{-} [u]  &   *={} \ar @{-} [l] \ar@{-} [u] &   *={} \ar @{-} [l] \ar@{-} [u] &   *={} \ar @{-} [l] \ar@{-} [u]  &   *={} \ar @{-} [l] \ar@{-} [u]\\
*={} \ar @{-} [u] &   *={} \ar@{-} [u] \ar @{-} [l]    &   *={} \ar @{-} [l] \ar@{-} [u] &   *={} \ar @{-} [l] \ar@{-} [u] &   *={} \ar @{-} [l] \ar@{-} [u] &   *={} \ar @{-} [l] \ar@{-} [u] \\
*={} \ar @{-} [u] &   *={} \ar @{-} [l] \ar@{-} [u]   &   *={} \ar @{-} [l] \ar@{-} [u] &   *={} \ar @{-} [l] \ar@{-} [u] \\
*={} \ar @{-} [u] &   *={} \ar @{-} [l] \ar@{-} [u]   &   *={} \ar @{-} [l] \ar@{-} [u] &   *={} \ar @{-} [l] \ar@{-} [u]}\endxy^{(\vphi_1)},\xy<0cm,.75cm>\xymatrix@R=.3cm@C=.3cm{
*={} & *={} \ar @{-} [l] & *={} \ar @{-} [l]  \\
*={} \ar @{-} [u] &   *={} \ar @{-} [l] \ar @{-} [u] &   *={} \ar @{-} [l] \ar@{-} [u]    \\ 
*={} \ar @{-} [u] &   *={} \ar @{-} [l] \ar @{-} [u]  &   *={} \ar @{-} [l] \ar@{-} [u] \\
*={} \ar @{-} [u] &   *={} \ar@{-} [u] \ar @{-} [l] &   *={} \ar @{-} [l] \ar@{-} [u]  \\
*={} \ar @{-} [u] &   *={} \ar @{-} [l] \ar@{-} [u] &   *={} \ar @{-} [l] \ar@{-} [u] \\
*={} \ar @{-} [u] &   *={} \ar @{-} [l] \ar@{-} [u] }\endxy^{(\vphi_2)}, 
\xy<0cm,.3cm>\xymatrix@R=.3cm@C=.3cm{
*={} &   *={} \ar @{-} [l]   \\
*={} \ar @{-} [u] &   *={} \ar @{-} [l] \ar@{-} [u]  \\
*={} \ar @{-} [u] &   *={} \ar @{-} [l] \ar@{-} [u] }\endxy^{(\vphi_4)}\right), \qquad \text{with $|\vphi_i|=i$,}$$
satisfies (a), (b), and (c) with $r=5$.  Since 
$$\quot_2(\blam^{(\vphi_1)})=\left(\xy<0cm,.3cm>\xymatrix@R=.3cm@C=.3cm{
*={} & *={} \ar @{-} [l] & *={} \ar @{-} [l]  & *={} \ar @{-} [l] \\
*={} \ar @{-} [u] &   *={} \ar @{-} [l] \ar @{-} [u] &   *={} \ar @{-} [l] \ar@{-} [u] &   *={} \ar @{-} [l] \ar@{-} [u] \\
*={} \ar @{-} [u] &   *={} \ar @{-} [l] \ar @{-} [u] &   *={} \ar @{-} [l] \ar@{-} [u] }\endxy^{(0)}, 
\xy<0cm,.45cm>\xymatrix@R=.3cm@C=.3cm{
*={} & *={} \ar @{-} [l] & *={} \ar @{-} [l]   \\
*={} \ar @{-} [u] &   *={} \ar @{-} [l] \ar @{-} [u] &   *={} \ar @{-} [l] \ar@{-} [u]  \\
*={} \ar @{-} [u] &   *={} \ar @{-} [l] \ar @{-} [u]  &   *={} \ar @{-} [l] \ar@{-} [u]\\
*={} \ar @{-} [u] &   *={} \ar @{-} [l] \ar @{-} [u] 
}\endxy^{(1)}\right)\qquad \text{and} \qquad \core_2(\blam^{(\vphi_1)})=(1)$$
(as in the example for (\ref{DominoToTableau})), we have that 
$$ |\core_2(\blam^{(\vphi_1)})|+2\quot_2(\blam^{(\vphi_1)})^{(1)}_{\lceil r/2\rceil}=3\qquad \text{and}\qquad \bgamma=\left(\xy<0cm,.15cm>\xymatrix@R=.3cm@C=.3cm{
*={} & *={} \ar @{-} [l] & *={} \ar @{-} [l]  & *={} \ar @{-} [l]\\
*={} \ar @{-} [u] &   *={} \ar @{-} [l] \ar @{-} [u] &   *={} \ar @{-} [l] \ar@{-} [u]&   *={} \ar @{-} [l] \ar@{-} [u] }\endxy^{(\vphi_1)}, \xy<0cm,.15cm>\xymatrix@R=.3cm@C=.3cm{
*={} & *={} \ar @{-} [l] \\
*={} \ar @{-} [u] &   *={} \ar @{-} [l] \ar @{-} [u] }\endxy^{(\vphi_2)}\right).$$

Consider the following battery tableau ${\bf Q}$ of shape $\blam/\bgamma$.   
\begin{enumerate}
\item For $\vphi\in \Theta$ such that $|\vphi|$ is even, fill $\blam^{(\vphi)}/\bgamma^{(\vphi)}$ with $\blam^{(\vphi)}_{2j-1}$ $j$'s and $\blam^{(\vphi)}_{2j}$ $\bar{j}$'s for $2j \leq \ell(\blam^{(\vphi)})$ such that $2j < r$, and $\blam^{(\vphi)}_{2j}$ $j$'s and $\blam^{(\vphi)}_{2j+1}$ $\bar{j}$'s for $2j+1 \leq \ell(\blam^{(\vphi)})$ such that $2j > r$.  Then all of the nonzero entries come in pairs $\xy<0cm,.3cm>\xymatrix@R=.3cm@C=.3cm{*={} & *={} \ar @{-} [l]\\  *={} \ar @{-} [u] &  *={} \ar @{-} [u] \ar @{-} [l] \ar @{} [ul]|{\sss \bar{j}}\\ *={} \ar @{-} [u] &  *={} \ar @{-} [u] \ar @{-} [l] \ar @{} [ul]|{\sss j}}\endxy$ and the resulting weight is lexicographically maximal.
\item For $\vphi\in \Theta$ such that $|\vphi|$ is odd, use Lemma \ref{dominobij} to fill $\blam^{(\vphi)}/\bgamma^{(\vphi)}$ in a lexicographically maximal way.
\end{enumerate}
In our running example, we have
$${\bf Q}=\left(\xy<0cm,.75cm>\xymatrix@R=.3cm@C=.3cm{
*={} & *={} \ar @{-} [l] & *={} \ar @{-} [l] & *={} \ar @{-} [l] & *={} \ar @{-} [l] & *={} \ar @{-} [l]  \\
*={} \ar @{-} [u] &   *={} \ar @{-} [l] \ar @{-} [u] \ar @{} [ul]|{\sss 0} &   *={} \ar @{-} [l] \ar@{-} [u]  \ar@{} [ul]|{\sss 0}  &   *={}  \ar @{-}[l] \ar @{-} [u] \ar@{} [ul]|{\sss 0}  &   *={} \ar @{} [l]|{\sss 1}  \ar @{-} [u]  &   *={} \ar @{} [l]|{\sss 1}   \ar @{-} [u]\\ 
*={} \ar @{-} [u] &   *={} \ar @{} [l]|{\sss 1} \ar @{-} [u]  &   *={} \ar @{} [l]|{\sss 1} \ar@{-} [u] &   *={} \ar @{} [l]|{\sss 1} \ar@{-} [u] &   *={} \ar @{-} [l] \ar@{-} [u]  &   *={} \ar @{-} [l] \ar@{-} [u]\\
*={} \ar @{-} [u] &   *={} \ar@{-} [u] \ar @{-} [l]    &   *={} \ar @{-} [l] \ar@{-} [u] &   *={} \ar @{-} [l] \ar@{-} [u] &   *={} \ar @{-} [l] \ar@{} [u]|{\sss 2} &   *={} \ar @{-} [l] \ar@{-} [u] \\
*={} \ar @{-} [u] &   *={} \ar @{} [l]|{\sss 2} \ar@{-} [u]   &   *={} \ar @{} [l]|{\sss 2} \ar@{-} [u] &   *={} \ar @{} [l]|{\sss 2} \ar@{-} [u] \\
*={} \ar @{-} [u] &   *={} \ar @{-} [l] \ar@{-} [u]   &   *={} \ar @{-} [l] \ar@{-} [u] &   *={} \ar @{-} [l] \ar@{-} [u]}\endxy^{(\vphi_1)},\xy<0cm,.75cm>\xymatrix@R=.3cm@C=.3cm{
*={} & *={} \ar @{-} [l] & *={} \ar @{-} [l]  \\
*={} \ar @{-} [u] &   *={} \ar @{-} [l] \ar @{-} [u] \ar @{} [ul]|{\sss 0}&   *={} \ar @{-} [l] \ar@{-} [u]  \ar@{} [ul]|{\sss \bar{1}}   \\ 
*={} \ar @{-} [u] &   *={} \ar @{-} [l] \ar @{-} [u] \ar@{} [ul]|{\sss \bar{1}}  &   *={} \ar @{-} [l] \ar@{-} [u] \ar@{} [ul]|{\sss 1} \\
*={} \ar @{-} [u] &   *={} \ar@{-} [u] \ar @{-} [l] \ar@{} [ul]|{\sss 1}  &   *={} \ar @{-} [l] \ar@{-} [u]  \ar@{} [ul]|{\sss \bar{2}}  \\
*={} \ar @{-} [u] &   *={} \ar @{-} [l] \ar@{-} [u] \ar@{} [ul]|{\sss \bar{2}}  &   *={} \ar @{-} [l] \ar@{-} [u] \ar@{} [ul]|{\sss 2} \\
*={} \ar @{-} [u] &   *={} \ar @{-} [l] \ar@{-} [u]\ar@{} [ul]|{\sss 2}  }\endxy^{(\vphi_2)}, 
\xy<0cm,.3cm>\xymatrix@R=.3cm@C=.3cm{
*={} &   *={} \ar @{-} [l]   \\
*={} \ar @{-} [u] &   *={} \ar @{-} [l] \ar@{-} [u] \ar@{} [ul]|{\sss \bar{1}}  \\
*={} \ar @{-} [u] &   *={} \ar @{-} [l] \ar@{-} [u] \ar@{} [ul]|{\sss 1} }\endxy^{(\vphi_4)}\right)$$

Note that by Lemma \ref{dominobij}, ${\bf Q}$ is the only battery tableau of shape $\blam/\bgamma$ and weight $\wt({\bf Q})$.  Thus, it suffices to show that there is no $\bnu\subseteq \blam$ with $|\bnu|=|\bgamma|$ and $\hgt(\bnu)\leq 1$ such that there exists a battery tableau ${\bf P}$ of shape $\blam/\bnu$ and weight $\wt({\bf Q})$.  

Since $|\bnu|=|\bgamma|$, we may think of moving from ${\bf Q}$ to ${\bf P}$ by shifting zero entries between $\vphi\in \Theta$ in ${\bf Q}$.  If $|\vphi|$ is even, it is clear from the construction of ${\bf Q}^{(\vphi)}$ that if we add a zero, an entry $<r/2$ is lost, while if we remove a zero, an entry $>r/2$ is gained.  Now consider when $|\vphi|$ is odd, with $Q= {\bf Q}^{(\vphi)}$ and $\lambda = \blam^{(\vphi)}$.  Apply the bijection (\ref{DominoToTableau}) to the domino tableau $Q$, and notice that from Remark 2 preceding Lemma \ref{dominobij}, our choice of $\bgamma^{(\vphi)}$ forces $Q^{(|\core_2(\lambda)|)}$ to have exactly $\quot_2(\lambda)^{(\core_2(\lambda))}_{\lceil r/2 \rceil}$ $0$'s.  Now, adding a pair of zero entries to or removing a pair of zero entries from $Q$ is the same as adding a zero to or removing a zero from $Q^{(|\core_2(\lambda)|)}$.  It is clear that adding a zero to $Q^{(|\core_2(\lambda)|)}$ results in losing an entry $<r/2$, which removes a domino with entry $<r/2$ in $Q$, and removing a zero from $Q^{(|\core_2(\lambda)|)}$ results in gaining an entry $>r/2$, which adds a domino with entry $>r/2$ in $Q$.  Thus, no matter how we change $\bgamma$ to $\bnu$, we are forced to change the weight of the full battery tableau to a lexicographically smaller weight.  So, there is no such $\bnu$ which leaves the weight unchanged, and uniqueness follows.
\end{proof}

\noindent\textbf{Remarks.}  Corollary \ref{OurMultiplicity} follows from Theorem \ref{MultiplicityOne}, since (a) and (b) are easily satisfied, and 
$$r=\max\{\ell(\blam^{(\vphi)})\ \mid\ \vphi\in \Theta, |\vphi| \text{ even}\}+1.$$

Another consequence of Theorem \ref{MultiplicityOne} is a result by Ohmori \cite{ohm}.

\begin{corollary}[Ohmori] \label{OhmoriMultiplicity}
Let $\blam \in \cP^{\Theta}_n$, and define the partition $\mu$ to have parts
$$\mu_j = \sum_{\varphi \in \Theta} |\varphi| \blam^{(\varphi)}_j. $$
Suppose that $\mu = (1^{m_1} 2^{m_2} \ldots)$ is such that $m_i$ is even for all $i$ except for the one value $i = k$, or that $m_i$ is always even, in which case we let $k=0$.  Define the partition $\nu$ to be $\nu = (1^{m_1/2} 2^{m_2/2} \cdots k^{(m_k - 1)/2} \cdots)$.  Then the irreducible character $\chi^{\blam}$ appears with multiplicity one in the degenerate Gelfand-Graev character $\Gamma_{(k, \nu)}$.
\end{corollary}

\begin{proof}
Note that if $\blam\in \cP^\Theta_n$ satisfies the hypotheses of the corollary, then for any $\vphi\in \Theta$ the partition $\blam^{(\vphi)}$ has at most one nonzero part size with odd multiplicity, otherwise $\mu$ would have more parts with odd multiplicity.  Thus, $\blam$ satisfies condition (b) of Theorem \ref{MultiplicityOne}.  Moreover, the fact that $\mu$ has at most one part with odd multiplicity implies that there must be an $r>0$ such that for every $\vphi \in \Theta$, either $\ell(\blam^{(\vphi)}) <r$ or $\blam^{(\vphi)}_r$ has odd multiplicity in $\blam^{(\vphi)}$ and $\blam^{(\vphi)}_r < \blam^{(\vphi)}_{r-1}$.  In particular, this holds when $|\vphi|$ is even, and so $\blam$ satisfies condition (c) of Theorem \ref{MultiplicityOne}.  If $\vphi\in \Theta$ is odd, and $\blam^{(\vphi)}$ has a part size $i$ with odd multiplicity, where $i=0$ if $\ell(\blam^{(\vphi)}) <r$, then 
$$\core_2(\blam^{(\vphi)})= \left\{\begin{array}{ll} (1) & \text{if $i$ is odd},\\ (0) & \text{if $i$ is even}.\end{array}\right.$$
Thus, $\blam$ satisfies (a) of Theorem \ref{MultiplicityOne}.

Now define $\bgamma$ by $\bgamma^{(\vphi)} = i$ if $m_i(\blam^{(\vphi)})$ is odd, where $i=0$ if $\ell(\blam^{(\vphi)}) <r$.  Then $|\bgamma| = k$.  When $|\vphi|$ is even, fill $\blam^{(\vphi)}/\bgamma^{(\vphi)}$ just as in the proof of Theorem \ref{MultiplicityOne}.  When $|\vphi|$ is odd, fill $\blam^{(\vphi)}/\bgamma^{(\vphi)}$ with all vertical dominoes such that there are $\blam^{(\vphi)}_{2j-1}$ $j$'s for $2j-1 \leq \ell(\blam^{(\vphi)})$ such that $2j-1 < r$, and $\blam^{(\vphi)}_{2j}$ $j$'s for $2j \leq \ell(\blam^{(\vphi)})$ such that $2j >r$.  This gives a battery tableau of shape $\blam/\bgamma$, where $|\bgamma|=k$, and weight $\nu$ as defined above.  

Fix a $\vphi$ such that $|\vphi|$ is odd, let $\lambda = \blam^{(\vphi)}$, and let ${\bf Q}^{(\vphi)}= Q$ be the domino tableau just defined, and apply the bijection (\ref{DominoToTableau}) to $Q$.  Since $Q$ has been filled with all vertical dominoes, the resulting weight is lexicographically maximal, and so by Lemma \ref{dominobij}, the tableaux $Q^{(0)}$ and $Q^{(1)}$ obtained from the bijection (\ref{DominoToTableau}) also have lexicographically maximal weights.  Let $j = |\core_2(\lambda)|$, and consider $Q^{(j)}$.  From our choice of $\bgamma^{(\vphi)}$, the tableau $Q^{(j)}$ has exactly $\lfloor \lambda_r/2 \rfloor$ $0$'s.  By the bijection (\ref{DominoToTableau}), we also have row $\lceil r/2 \rceil$ of $Q^{(j)}$, which is $\quot_2(\lambda^{(\vphi)})_{\lceil r/2 \rceil}^{(j)}$, is exactly $\lfloor \lambda_r/2 \rfloor$.  This means the domino tableau ${\bf Q}^{(\vphi)}$ is exactly what is constructed in the proof of Theorem \ref{MultiplicityOne}.  Therefore, ${\bf Q}$ is exactly the battery tableau obtained in the proof of Theorem \ref{MultiplicityOne}, and so we have $\langle \chi^{\blam}, \Gamma_{(k,\nu)} \rangle = 1$.
  \end{proof}

  Note that by Corollary \ref{decompcor}, condition (a) of Theorem \ref{MultiplicityOne} is a necessary condition.  The following proposition shows that condition (b) is also necessary.

\begin{proposition} Let $\blam\in \cP^\Theta$.  If there exists a $\vphi \in \Theta$ such that $|\vphi|$ is even and $\blam^{(\vphi)}$ has at least two distinct part sizes with odd multiplicity, then 
$$\langle\chi^{\blam},\Gamma_{(k,\nu)}\rangle\neq 1$$
 for all degenerate Gelfand-Graev characters $\Gamma_{(k,\nu)}$.
\end{proposition}

\begin{proof}
Suppose $\blam\in \cP^\Theta$ and $|\vphi|$ is even, such that $\lambda=\blam^{(\vphi)}$ has part sizes $x<y$ with odd multiplicity.  Let $Q$ be a symplectic tableau of shape $\lambda/(m)$ for some $m\leq \lambda_1$, and suppose $\wt(Q)=\mu$.  If there exists an $\bar{i}$ such that there is no $i$ directly south of $\bar{i}$ in $Q$, then there is a second symplectic tableau $P$ of shape $\lambda/(m)$ and weight $\mu$ obtained by changing this $\bar{i}$ to an $i$ in $Q$.  Similarly, if there is an $i$ with no $\bar{i}$ directly north of it, then there is a second tableau $P$ with the same weight and shape as $Q$.  Thus, the only way $Q$ is the only tableau of shape $\lambda/(m)$ and weight $\mu$, is if $\lambda/(m)$ can be tiled by vertical dominoes.

If $m<y$, then the $y$th column of $\lambda/(m)$ has an odd number of boxes, and therefore cannot be tiled by vertical dominoes.  If $m>y$, then the $m$th column of $\lambda/(m)$ has an odd number of boxes.  If $m=y$, then the $x$th column of $\lambda/(m)$ has an odd number of boxes.  In all cases, $\lambda/(m)$ cannot be tiled by dominoes, and the result follows.
\end{proof}

\noindent\textbf{Remarks.}  
\begin{enumerate}
\item While conditions (a) and (b) of Theorem \ref{MultiplicityOne} are necessary, condition (c) is not. For example, the only battery tableau of weight $(2,(8))$ for the $\Theta$-partition
$$\blam=\left(\xy<0cm,.3cm>\xymatrix@R=.3cm@C=.3cm{
	*={} & *={} \ar @{-} [l]  \\
	*={} \ar @{-} [u] &   *={} \ar @{-} [l] \ar@{-} [u]   \\
*={} \ar @{-} [u] &   *={} \ar @{-} [l] \ar@{-} [u]   }\endxy^{(\alpha)},\ \xy<0cm,.2cm>
\xymatrix@R=.3cm@C=.3cm{
	*={} & *={} \ar @{-} [l] \\
	*={} \ar @{-} [u] &   *={} \ar @{-} [l] \ar@{-} [u]}\endxy^{(\beta)}\right)
\quad\text{with $|\alpha|=4, |\beta|=2$, is}\quad
{\bf Q}=\left( \xy<0cm,.3cm>\xymatrix@R=.3cm@C=.3cm{
	*={} & *={} \ar @{-} [l]  \\
	*={} \ar @{-} [u] &   *={} \ar @{-} [l] \ar@{-} [u]\ar @{} [ul]|{\sss \bar{1}}   \\
*={} \ar @{-} [u] &   *={} \ar @{-} [l] \ar@{-} [u] \ar @{} [ul]|{\sss 1}  }\endxy^{(\alpha)},\ \xy<0cm,.2cm>\xymatrix@R=.3cm@C=.3cm{
	*={} & *={} \ar @{-} [l] \\
	*={} \ar @{-} [u] &   *={} \ar @{-} [l] \ar@{-} [u] \ar @{} [ul]|{\sss 0}}\endxy^{(\beta)} \right).
	$$
\item At the same time, conditions (a) and (b) of Theorem \ref{MultiplicityOne} are not alone sufficient.  For example, 
$$\blam=\left(\xy<0cm,.3cm>\xymatrix@R=.3cm@C=.3cm{
	*={} & *={} \ar @{-} [l]  \\
	*={} \ar @{-} [u] &   *={} \ar @{-} [l] \ar@{-} [u]   \\
*={} \ar @{-} [u] &   *={} \ar @{-} [l] \ar@{-} [u]   }\endxy^{(\alpha)}
,\  \xy<0cm,.3cm>\xymatrix@R=.3cm@C=.3cm{
	*={} & *={} \ar @{-} [l] \\
	*={} \ar @{-} [u] &   *={} \ar @{-} [l] \ar@{-} [u]}\endxy^{(\beta)}\ \right),
\quad\text{with $|\alpha|=|\beta|=2$},$$
satisfies (a) and (b).  The possible weights and two of their battery tableaux are 
\begin{align*}
(0,(6)):\qquad 
&\left(\xy<0cm,.3cm>\xymatrix@R=.3cm@C=.3cm{
	*={} & *={} \ar @{-} [l]  \\
	*={} \ar @{-} [u] &   *={} \ar @{-} [l] \ar@{-} [u]  \ar @{} [ul]|{\sss \bar{1}} \\
*={} \ar @{-} [u] &   *={} \ar @{-} [l] \ar@{-} [u] \ar @{} [ul]|{\sss 1}  }\endxy^{(\alpha)}
,\  \xy<0cm,.3cm>\xymatrix@R=.3cm@C=.3cm{
	*={} & *={} \ar @{-} [l] \\
	*={} \ar @{-} [u] &   *={} \ar @{-} [l] \ar@{-} [u]\ar @{} [ul]|{\sss \bar{1}}}\endxy^{(\beta)}\ \right),\qquad
	\left(\xy<0cm,.3cm>\xymatrix@R=.3cm@C=.3cm{
	*={} & *={} \ar @{-} [l]  \\
	*={} \ar @{-} [u] &   *={} \ar @{-} [l] \ar@{-} [u]  \ar @{} [ul]|{\sss \bar{1}} \\
*={} \ar @{-} [u] &   *={} \ar @{-} [l] \ar@{-} [u] \ar @{} [ul]|{\sss 1}  }\endxy^{(\alpha)}
,\  \xy<0cm,.3cm>\xymatrix@R=.3cm@C=.3cm{
	*={} & *={} \ar @{-} [l] \\
	*={} \ar @{-} [u] &   *={} \ar @{-} [l] \ar@{-} [u]\ar @{} [ul]|{\sss 1}}\endxy^{(\beta)}\ \right)\qquad (\text{2 total}), \\
(0,(4,2)): \qquad &\left(\xy<0cm,.3cm>\xymatrix@R=.3cm@C=.3cm{
	*={} & *={} \ar @{-} [l]  \\
	*={} \ar @{-} [u] &   *={} \ar @{-} [l] \ar@{-} [u]  \ar @{} [ul]|{\sss \bar{1}} \\
*={} \ar @{-} [u] &   *={} \ar @{-} [l] \ar@{-} [u] \ar @{} [ul]|{\sss 1}  }\endxy^{(\alpha)}
,\  \xy<0cm,.3cm>\xymatrix@R=.3cm@C=.3cm{
	*={} & *={} \ar @{-} [l] \\
	*={} \ar @{-} [u] &   *={} \ar @{-} [l] \ar@{-} [u]\ar @{} [ul]|{\sss \bar{2}}}\endxy^{(\beta)}\ \right),\qquad 
	\left(\xy<0cm,.3cm>\xymatrix@R=.3cm@C=.3cm{
	*={} & *={} \ar @{-} [l]  \\
	*={} \ar @{-} [u] &   *={} \ar @{-} [l] \ar@{-} [u]  \ar @{} [ul]|{\sss \bar{1}} \\
*={} \ar @{-} [u] &   *={} \ar @{-} [l] \ar@{-} [u] \ar @{} [ul]|{\sss 1}  }\endxy^{(\alpha)}
,\  \xy<0cm,.3cm>\xymatrix@R=.3cm@C=.3cm{
	*={} & *={} \ar @{-} [l] \\
	*={} \ar @{-} [u] &   *={} \ar @{-} [l] \ar@{-} [u]\ar @{} [ul]|{\sss 2}}\endxy^{(\beta)}\ \right)\qquad (\text{10 total}),\\
(0,(2^3)):\qquad &\left(\xy<0cm,.3cm>\xymatrix@R=.3cm@C=.3cm{
	*={} & *={} \ar @{-} [l]  \\
	*={} \ar @{-} [u] &   *={} \ar @{-} [l] \ar@{-} [u]  \ar @{} [ul]|{\sss \bar{1}} \\
*={} \ar @{-} [u] &   *={} \ar @{-} [l] \ar@{-} [u] \ar @{} [ul]|{\sss \bar{2}}  }\endxy^{(\alpha)}
,\  \xy<0cm,.3cm>\xymatrix@R=.3cm@C=.3cm{
	*={} & *={} \ar @{-} [l] \\
	*={} \ar @{-} [u] &   *={} \ar @{-} [l] \ar@{-} [u]\ar @{} [ul]|{\sss \bar{3}}}\endxy^{(\beta)}\ \right),\qquad 
	\left(\xy<0cm,.3cm>\xymatrix@R=.3cm@C=.3cm{
	*={} & *={} \ar @{-} [l]  \\
	*={} \ar @{-} [u] &   *={} \ar @{-} [l] \ar@{-} [u]  \ar @{} [ul]|{\sss 1} \\
*={} \ar @{-} [u] &   *={} \ar @{-} [l] \ar@{-} [u] \ar @{} [ul]|{\sss 2}  }\endxy^{(\alpha)}
,\  \xy<0cm,.3cm>\xymatrix@R=.3cm@C=.3cm{
	*={} & *={} \ar @{-} [l] \\
	*={} \ar @{-} [u] &   *={} \ar @{-} [l] \ar@{-} [u]\ar @{} [ul]|{\sss 3}}\endxy^{(\beta)}\ \right)\qquad (\text{24 total}),\\
(2,(4)):\qquad &\left(\xy<0cm,.3cm>\xymatrix@R=.3cm@C=.3cm{
	*={} & *={} \ar @{-} [l]  \\
	*={} \ar @{-} [u] &   *={} \ar @{-} [l] \ar@{-} [u]  \ar @{} [ul]|{\sss 0} \\
*={} \ar @{-} [u] &   *={} \ar @{-} [l] \ar@{-} [u] \ar @{} [ul]|{\sss \bar{1}}  }\endxy^{(\alpha)}
,\  \xy<0cm,.3cm>\xymatrix@R=.3cm@C=.3cm{
	*={} & *={} \ar @{-} [l] \\
	*={} \ar @{-} [u] &   *={} \ar @{-} [l] \ar@{-} [u]\ar @{} [ul]|{\sss \bar{1}}}\endxy^{(\beta)}\ \right),\qquad 
	\left(\xy<0cm,.3cm>\xymatrix@R=.3cm@C=.3cm{
	*={} & *={} \ar @{-} [l]  \\
	*={} \ar @{-} [u] &   *={} \ar @{-} [l] \ar@{-} [u]  \ar @{} [ul]|{\sss 0} \\
*={} \ar @{-} [u] &   *={} \ar @{-} [l] \ar@{-} [u] \ar @{} [ul]|{\sss 1}  }\endxy^{(\alpha)}
,\  \xy<0cm,.3cm>\xymatrix@R=.3cm@C=.3cm{
	*={} & *={} \ar @{-} [l] \\
	*={} \ar @{-} [u] &   *={} \ar @{-} [l] \ar@{-} [u]\ar @{} [ul]|{\sss 1}}\endxy^{(\beta)}\ \right)\qquad (\text{5 total}),\\
(2,(2^2)):\qquad &\left(\xy<0cm,.3cm>\xymatrix@R=.3cm@C=.3cm{
	*={} & *={} \ar @{-} [l]  \\
	*={} \ar @{-} [u] &   *={} \ar @{-} [l] \ar@{-} [u]  \ar @{} [ul]|{\sss 0} \\
*={} \ar @{-} [u] &   *={} \ar @{-} [l] \ar@{-} [u] \ar @{} [ul]|{\sss \bar{1}}  }\endxy^{(\alpha)}
,\  \xy<0cm,.3cm>\xymatrix@R=.3cm@C=.3cm{
	*={} & *={} \ar @{-} [l] \\
	*={} \ar @{-} [u] &   *={} \ar @{-} [l] \ar@{-} [u]\ar @{} [ul]|{\sss \bar{2}}}\endxy^{(\beta)}\ \right),\qquad 
	\left(\xy<0cm,.3cm>\xymatrix@R=.3cm@C=.3cm{
	*={} & *={} \ar @{-} [l]  \\
	*={} \ar @{-} [u] &   *={} \ar @{-} [l] \ar@{-} [u]  \ar @{} [ul]|{\sss 0} \\
*={} \ar @{-} [u] &   *={} \ar @{-} [l] \ar@{-} [u] \ar @{} [ul]|{\sss 1}  }\endxy^{(\alpha)}
,\  \xy<0cm,.3cm>\xymatrix@R=.3cm@C=.3cm{
	*={} & *={} \ar @{-} [l] \\
	*={} \ar @{-} [u] &   *={} \ar @{-} [l] \ar@{-} [u]\ar @{} [ul]|{\sss 2}}\endxy^{(\beta)}\ \right)\qquad (\text{12 total}),\\
(4,(2)):\qquad &\left(\xy<0cm,.3cm>\xymatrix@R=.3cm@C=.3cm{
	*={} & *={} \ar @{-} [l]  \\
	*={} \ar @{-} [u] &   *={} \ar @{-} [l] \ar@{-} [u]  \ar @{} [ul]|{\sss 0} \\
*={} \ar @{-} [u] &   *={} \ar @{-} [l] \ar@{-} [u] \ar @{} [ul]|{\sss \bar{1}}  }\endxy^{(\alpha)}
,\  \xy<0cm,.3cm>\xymatrix@R=.3cm@C=.3cm{
	*={} & *={} \ar @{-} [l] \\
	*={} \ar @{-} [u] &   *={} \ar @{-} [l] \ar@{-} [u]\ar @{} [ul]|{\sss 0}}\endxy^{(\beta)}\ \right),\qquad 
	\left(\xy<0cm,.3cm>\xymatrix@R=.3cm@C=.3cm{
	*={} & *={} \ar @{-} [l]  \\
	*={} \ar @{-} [u] &   *={} \ar @{-} [l] \ar@{-} [u]  \ar @{} [ul]|{\sss 0} \\
*={} \ar @{-} [u] &   *={} \ar @{-} [l] \ar@{-} [u] \ar @{} [ul]|{\sss 1}  }\endxy^{(\alpha)}
,\  \xy<0cm,.3cm>\xymatrix@R=.3cm@C=.3cm{
	*={} & *={} \ar @{-} [l] \\
	*={} \ar @{-} [u] &   *={} \ar @{-} [l] \ar@{-} [u]\ar @{} [ul]|{\sss 0}}\endxy^{(\beta)}\ \right)\qquad (\text{2 total}).
\end{align*}
\end{enumerate}

\end{document}